\documentclass{article}
\usepackage{authblk}
\usepackage{graphicx,ulem} 
\usepackage{color}
\usepackage{amsfonts,amsmath}
\usepackage{amsthm}
\usepackage{empheq}
\usepackage{mathtools}
\usepackage{multirow}
\usepackage{caption}
\usepackage{graphicx}
\usepackage[toc,page]{appendix}
\usepackage{subcaption}
\usepackage[dvipsnames]{xcolor}

\usepackage{hyperref}
\hypersetup{
    colorlinks=true,
    linkcolor=blue,
    filecolor=magenta,      
    urlcolor=cyan,
    citecolor=red
}

\graphicspath{{./Images/HT1v3/}}

\newcommand{\lfd}{\lambda_{F, D}}
\newcommand{\lfw}{\lambda_{F}}

\newcommand{\cI}{\mathcal{I}}
\newcommand{\mW}{\tilde{m}_W}
\newcommand{\mD}{\tilde{m}_D}
\newcommand{\R}{\mathbb{R}}
\newcommand{\N}{\mathcal{N}}

\newcommand{\YD}[1]{\textcolor{black}{#1}}

\newcommand{\vtrois}[1]{\textcolor{OliveGreen}{#1}}
\newcommand{\vquatre}[1]{\textcolor{black}{#1}}

\DeclareMathOperator{\Tr}{Tr}
\newtheorem{theorem}{Theorem}
\newtheorem{prop}[theorem]{Proposition}
\theoremstyle{definition}
\newtheorem{definition}[theorem]{Definition}
\theoremstyle{remark}
\newtheorem{remark}[theorem]{Remark}

\title{Human-Wildlife interactions in a tropical forest context: modeling, analysis and simulations}
\author[1,2,3]{Yves Dumont \footnote{Corresponding author: yves.dumont@cirad.fr}} 
\author[2]{Marc H\'etier}
\author[2,4]{Valaire Yatat-Djeumen}

\affil[1]{CIRAD, UMR AMAP, 3P, F-97410 Saint Pierre, France}
\affil[2]{AMAP, Univ Montpellier, CIRAD, CNRS, INRAe, IRD, Montpellier, France}
\affil[3]{University of Pretoria, Department of Mathematics and Applied Mathematics, Pretoria, South Africa}
\affil[4]{ENSPY, University of Yaound\'e 1, Cameroon}
\begin{document}

\maketitle

\section*{Abstract}
Anthropisation and excessive hunting in tropical forests threaten biodiversity, ecosystem maintenance and human food security. In this article, we focus on the issue of coexistence between humans and wildlife in an anthropised environment. Assuming that the human population moves between its residential area and the surrounding forest to hunt, we study a resource-consumer model with consumer migration. A comprehensive analysis of the system is carried out using classical theory and monotone systems theory. We show that the possibilities for long-term coexistence between human populations and wildlife populations are determined by hunting rate thresholds. Depending on the level of anthropisation and the hunting rate, the system may converge towards a limit cycle or a co-existence equilibrium. However, the conditions for coexistence become more difficult as anthropisation increases. Numerical simulations are provided to illustrate the theoretical results.

\vspace{0.5cm}
\textbf{Keywords}: Human-Wildlife system; Anthropisation; Tropical Forest; Resource-consumer model; Tikhonov theorem; Monotone systems theory; Periodic cycle; Numerical simulations.
\section{Introduction}

Tropical forests are particularly rich ecosystems, in term of plants and animals diversity \cite{rajpar_tropical_2018}. They also provide resources for forest-based people, and other human populations living nearby. These resources include food \cite{avila_martin_food_2024}, medicine
and energy sources \cite{mangula_energy_2019} for example.

However, tropical forests are increasingly degraded and fragmented by human settlement and activities, as the development of infrastructures (roadways, harbors, dams, ...), agriculture, and industrial complexes (industrial palm groves, industrial logging,...), etc. 
Today, only 20\% of the remaining area are considered as intact, see \cite{benitez-lopez_intact_2019}. It is not only the vegetation that is endangered, but also the wildlife that it shelters. Of course, the destruction of their habitat has indirect consequences on animal populations, but they are also directly threatened by human activities, and especially by over-hunting \cite{benitez-lopez_intact_2019, wilkie_empty_2011}. In 2019, it was estimated in \cite{benitez-lopez_intact_2019} that the abundance of tropical mammal species had declined by an average of 13\%. The authors also predicted a decline of over 70\% for mammals in West Africa.

Defaunation is not harmless. It impacts both vegetation and human populations. In \cite{rajpar_tropical_2018, ripple_bushmeat_2016}, the authors point out that even partial defaunation has consequences on the environment and human populations. Indeed, extinction of mammals affects the forest regeneration, as, for instance they play a role in seed dispersal \cite{estrada-villegas_animal_2023,peres_dispersal_2016, wright_bushmeat_2007}. Defaunation is also associated with the (re)emergence of zoonotic diseases, see for example \cite{dobigny_zoonotic_2022, white_emerging_2020}. 

Moreover, bushmeat remains an important source of food and income for certain populations, see \cite{jones_incentives_2019}. It is the case in South Cameroon, where some of the populations (as the Baka Pygmies) still rely on hunting and agriculture for their livelihood \cite{avila_martin_food_2024}. These populations do not threaten the vegetation: consumption of plant for food or medical uses is low compared to forest regeneration \cite{koppert_consommation_1996}. Furthermore, we cannot ignore the fact that  the forest is also impacted by industrial development: a deep-water harbor has been built near the town of Kribi, and industrial palm groves are expanding rapidly in the surrounding area. These new constructions are supported by the development of roads, which fragmented the environment (see for example \cite{foonde_change_2018, romain_deforestation_2017} about landscape and biodiversity changes around Kribi).

It can be difficult to determine whether a specific hunting, as practised in a particular place, is sustainable or not. Indeed, acquiring knowledge about wildlife and hunting practices can not be done using remote-sensing \cite{peres_detecting_2006}. It requires specific field studies such as hunter recall interviews or continuous monitoring. These surveys are costly and susceptible to biases, see \cite{jones_consequences_2020}.

However, we believe that mathematical modeling can provide some insights on these issues. Indeed, mathematical models  can not only synthesise and generalise a complex reality, determine parameters and thresholds of interest, but also facilitate future field studies and help decision-makers make informed decisions \cite{deangelis_towards_2021}.

Previous studies have modeled human-environment (including both vegetation and wildlife) interactions in different ways.
Some of them consider that the environment provides necessary services to agriculture \cite{bengochea-paz_agricultural_2020, roman_dynamics_2018}, and dealt with how and when social-ecological systems collapse. Other studies (see for example \cite{bulte_habitat_2003, nlom_bio-economic_2021}), in the bio-economic fields, reason in terms of economic costs and profits, seeking to determine the most effective wildlife management strategy. In \cite{fanuel_modelling_2023}, the authors review articles modeling the impact of human activity and human pressure on forest biomass and wildlife. The dynamical models presented rarely focus on over-hunting or on the coexistence between human populations and wildlife. Instead, they focused on the impact of specific types of anthropisation, such as industrialisation \cite{agarwal_depletion_2010}, mining activities or human pressure \cite{dubey_modelling_2009} on  forest biomass. Moreover, the majority of these models involves many (at least 4) variables, which complicates the theoretical analysis and the understanding of the models.

This work aims to investigate the impact of overhunting and the anthropisation of habitat on the possibility of coexistence between a forest based human population and wildlife. It also highlights how the presence of a constant human immigration affects this possibility.
It takes place within the I-CARE project (Impact of the Anthropisation on the Risk of Emergence of Zoonotic Arboviroses in Central Africa). This project, which involves field experts, focuses on South Cameroon. As stated above, the vegetation in this area is not directly threatened but is quickly and abruptly impacted by industrial projects. Moreover, due to the interdisciplinary aspect of the project, there is a need to keep the model as simple as possible, in order, at least, to preserve a deep mathematical tractability property and interpretability, but also to be able to include an epidemiological aspect at a later stage.

The manuscript is organized as follows. In section \ref{sec:model}, a dynamical model describing human-wildlife interactions is introduced. In sections \ref{sec:qssa} and \ref{sec:competitive}, we propose a mathematical analysis of this model. First, we present a simple but inconclusive quasi-steady state approach. Then, we study a more complicated but complete analysis using an equivalent competitive system. In section \ref{sec:ecological}, we provide an ecological interpretation of the results we obtained, as well as numerical simulations. Finally, we end the paper with some conclusions and perspectives in section \ref{sec:conclusion}.

\section{A Model of Human-Wildlife Interactions} \label{sec:model}
We model the interactions between a hunter community and wilderness using a consumer-resource model. We consider two areas, one corresponding to a domestic area (a village for example), the second to a wild area (a tropical forest for example). Humans in the domestic area, $H_D$, are modeled using a consumer equation based on available resource. The dynamic of $H_D$ includes several terms. First there is a constant growth term, denoted by $\cI$, representing the immigration from other inhabited areas. For instance, the development of an industrial complex attracts new workers. A second term in the equation is the natural death-rate of humans, noted $\mu_D$. Furthermore, we assume that the villagers are able to produce a certain amount of food at the rate $f_D$. Note that $f_D$ can be understand as $f_D = e \lfd F_D$ where $F_D$ is a constant amount of cattle or agricultural resources, $\lfd$ a harvest rate and $e$ the proportion of these resources in human diet. 

The villagers come and go between the domestic and wild areas, mainly for hunting activities. These displacements are modeled by the function $M(H_D, H_W) = m_W H_W - m_D H_D$, where $H_W$ is the human population presents in the wild area. We note $m_D$ (respectively $m_W$) the migration rate from the domestic (resp. wild) area to the wild (resp. domestic) area, and we define $m = \dfrac{m_D}{m_W}$. The dynamic of $H_W$ simply corresponds to the migration process. 

Wild animals, $F_W$, hunted in the wild area are partially used to feed the villagers. Therefore, we consider a last growth term in the $H_D$ dynamic, $e \Lambda_F(F_W) H_W$, where $\Lambda_F(F_W)$ is the number of prey hunted by one hunter, and $e$ is the proportion of wild prey in diet. This parameter $e$ encompasses both the direct consumption of the meat by the villagers and the fact that a large proportion of bush-meat is sold on city market, allowing the villagers to purchase other supply, see \cite{wilkie_bushmeat_1998}.

The dynamic of $F_W$ follows a logistic equation, whose carrying capacity, $K_F$, depends on the surrounding vegetation. To take into account the level of anthropisation, we introduce the non-negative parameter $\alpha \in [0, 1)$. When $\alpha > 0$, the carrying capacity of the vegetation is reduced of $\alpha \%$ from its original value. anthropisation may also have a negative impact on the animal's growth rate $r_F$. For the sake of simplicity, we model this impact in the same way, by multiplying the growth rate by $(1-\alpha)$.

We consider two different interactions between wild fauna and the population located in the wild area. Firstly, as mentioned above, the wild fauna is hunted by humans present in the wild area. We model this predation by a Holling-Type 1 functional $\Lambda_F(F_W) =  \lfw F_W $, where $\lfw$ is the hunting rate. We choose an unbounded functional response to take into account the possibility of over-hunt. \vquatre{On the other hand, according to discussion with field experts of the I-CARE project, the presence of human population may benefit the reproduction of certain wild species, particularly rodents. See also \cite{dobigny_zoonotic_2022, dounias_foraging_2011}.}
We take this effect into account by multiplying the wild animal's growth rate, $r_F$, by the functional $(1+\beta H_W)$, where $\beta \vquatre{\geq 0}$ is the positive influence rate that human activities may have on wild animal's growth.

\medskip
Finally, the model is given by the following equations:
\begin{equation}
\def\arraystretch{2}
\left\{ 
\begin{array}{l}
\dfrac{dH_D}{dt}= \cI + e\lfw H_W F_W + (f_D - \mu_D) H_D - m_D H_D + m_W H_W, \\
\dfrac{dF_W}{dt} = r_F(1- \alpha) (1+ \beta H_W) \left(1 - \dfrac{F_W}{K_F(1-\alpha)} \right) F_W - \lfw F_W H_W,\\
\dfrac{dH_W}{dt}= m_D H_D - m_W H_W.
\end{array} \right.
\label{equation:HDFWHW}
\end{equation}

For the rest of the study, we will note $y$ the vector of state variables, that is $y = \Big(H_D, F_W, H_W \Big)$ and $f(y)$ the right hand side of the system.

\begin{figure}[!ht]
\centering
\includegraphics[width=0.5\textwidth]{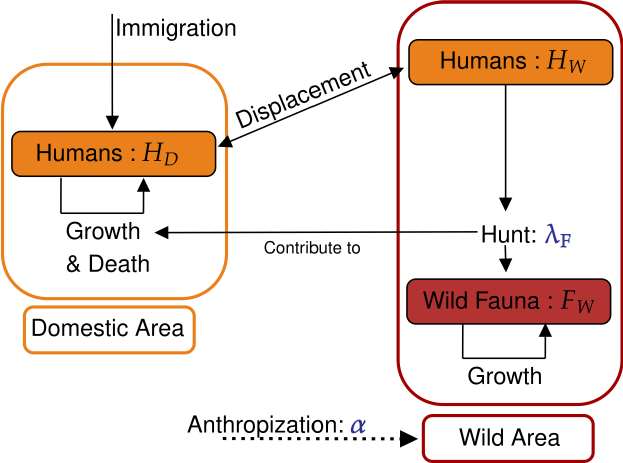}
\caption{System flow chart.}
\label{fig:flow chart}
\end{figure}

The whole dynamic is represented through the flow chart in figure \ref{fig:flow chart}. Table \ref{table:symbol} summarizes the meaning of the state variables and parameters.

\begin{table}[ht]
\center
\begin{tabular}{|c|c|c|}
\hline 
Symbol & Description & Unit \\ 
\hline \hline
$H_D$ & Humans in the domestic area & Ind \\
$H_W$ & Humans in the wild area & Ind \\
$F_W$ & Wild fauna & Ind \\
\hline \hline
$t$ & Time & Year \\
$e$ & Proportion of wild meat in the diet & - \\
$f_D$ & Food produced by human population & Year$^{-1}$ \\
$\mu_D$ & Human mortality rate  & Year$^{-1}$ \\
$m_D$ & Migration from domestic area to wild area & Year$^{-1}$ \\
$m_W$ & Migration from wild area to domestic area & Year$^{-1}$ \\
$r_F$ & Wild animal growth rate & Year$^{-1}$ \\
$K_F$ & Carrying capacity for wild fauna, fixed by the environment& Ind \\
$\alpha$ & Proportion of anthropized environment & - \\
$\beta$ & Positive impact from human activities \YD{on wild} animal's growth rate & Ind$^{-1}$  \\
$\lfw$ & Hunting rate & Ind$^{-1}\times$ Year$^{-1}$\\
$\mathcal{I}$ & Immigration rate &Ind $\times$ Year$^{-1}$\\
\hline
\end{tabular}
\caption{State variables and parameters of the model.}
\label{table:symbol}
\end{table}

To analyze system \eqref{equation:HDFWHW} and its long term behavior, we will use divers notions and theorems; they are recalled in Appendix \ref{sec:litterature theorems}, page \pageref{sec:litterature theorems}. We start by showing that the problem is well posed.

\subsection{Existence and Uniqueness of Global Solutions}
In this section, we state general results on system \eqref{equation:HDFWHW}:  existence of an invariant region, existence and uniqueness of global solutions.

\begin{remark}
To avoid infinite growth of human population, and ensure that the system is well defined, we need to add a constraint on the sign of $f_D - \mu_D$: on the following, we will assume that $f_D - \mu_D < 0$. This means that the food produced by the human population living in the domestic area is not sufficient to ensure the permanence of the population. Hunt, or immigration, is necessary. 
\end{remark}

The following proposition indicates a compact and invariant subset of $\mathbb{R}_+^3$ on which the solutions of system \eqref{equation:HDFWHW} are bounded.

\begin{prop}\label{prop:invariantRegion} 

Assume \YD{$0 \leq \beta <\beta^*$, with }
\begin{equation}
\label{betastar}
\beta^* := \dfrac{4(\mu_D - f_D)}{m e r_F (1-\alpha)^2 K_F}.
\end{equation}
Then, the region
$$\Omega = \Big\{\Big(H_D, F_W, H_W \Big) \in \mathbb{R}_+^3  \Big|H_D + eF_W \leq S^{max}, F_W \leq F_W^{max}, H_W \leq H_W^{max} \Big\},$$
is a compact and invariant set for system \eqref{equation:HDFWHW}, 
where
$$
S^{max} = \dfrac{\cI + \left( \mu_D - f_D   +  \dfrac{(1-\alpha)r_F}{4}\right)e(1-\alpha)K_F}{  \Big(\mu_D -f_D - \dfrac{er_F (1-\alpha)^2K_F m}{4} \beta\Big)} ,
\quad
F_W^{max} = (1-\alpha)K_F,
\quad
H_W^{max} = \dfrac{m_D}{m_W} S^{max}.
$$
In particular, this means that any solutions of equations \eqref{equation:HDFWHW} with initial condition in $\Omega$ are bounded and remain in $\Omega$.
\end{prop}
\begin{proof} 
To prove this proposition we will use the notion of invariant region, see \cite{smoller_shock_1994}. Before, we introduce the variable $S = H_D + e F_W$. 
With this new variable, the model writes:

\begin{equation}
\left\{ \begin{array}{l}
\dfrac{dS}{dt} = \cI + (f_D - \mu_D - m_D) \Big(S - eF_W \Big) + m_WH_W + e (1-\alpha)(1+\beta H_W)r_F  \left(1 - \dfrac{F_W}{(1-\alpha)K_F} \right) F_W, \\
\dfrac{dF_W}{dt} = (1-\alpha)(1+\beta H_W) r_F \left(1 - \dfrac{F_W}{K_F(1-\alpha)} \right) F_W - \lfw F_W H_W, \\
\dfrac{dH_W}{dt}= m_D \left(S - eF_W\right) - m_W H_W.
\end{array} \right.
\label{equation:SFWHW}
\end{equation}

We define the function $g(z)$, with $z=\Big(S, F_W, H_W \Big)$, as the right hand side of system \eqref{equation:SFWHW}. We also introduce the following functions:
$$
G_1(z) = S - S^{max},
\quad
G_2(z) = F_W - F_W^{max},
\quad
G_3(z) = H_W - H_W^{max}.
$$

Following \cite{smoller_shock_1994}, we will show that quantities $(\nabla G_1 \cdot g)|_{S = S^{max}}$, $(\nabla G_2 \cdot g)|_{F_W = F_W^{max}}$ and $(\nabla G_3 \cdot g)|_{H_W = H_W^{max}}$ are non-positive for $z \in \Omega_S = \Big\{ \Big(S, F_W, H_W \Big) \in (\mathbb{R})^3  \Big|S \leq S^{max}, F_W \leq F_W^{max}, H_W \leq H_W^{max} \Big\}$.

We have:
\begin{multline*}
(\nabla G_1 \cdot g)|_{S^{max}} = \cI + (f_D - \mu_D -m_D) \Big(S^{max} - eF_W \Big) +m_WH_W + e r_F(1-\alpha)(1+\beta H_W)  \left(1 - \dfrac{F_W}{(1-\alpha)K_F} \right) F_W.
\end{multline*}
Using that $\left(1-\dfrac{F_W}{(1-\alpha)K_F}\right) F_W \leq \dfrac{(1-\alpha)K_F}{4}$ for $0\leq F_W \leq (1-\alpha)K_F$ and $H_W \leq H_W^{max}$,

\begin{multline*}
(\nabla G_1 \cdot g)|_{S^{max}} \leq \cI + (f_D - \mu_D - m_D) S^{max} + m_W H_W^{max} - \Big(f_D - \mu_D \Big) e(1-\alpha)K_F + \\ er_F (1-\alpha)(1+\beta H_W^{max}) \dfrac{(1-\alpha)K_F}{4}.
\end{multline*}

Using $H_W^{max} = \dfrac{m_D}{m_W}S^{max}$,

\begin{align*}
(\nabla G_1 \cdot g)|_{S^{max}} &\leq \cI + (f_D - \mu_D) S^{max} - \Big(f_D - \mu_D \Big) e(1-\alpha)K_F +  er_F (1-\alpha)(1+\beta m S^{max}) \dfrac{(1-\alpha)K_F}{4}, \\
& \leq \cI + \Big(f_D - \mu_D + \dfrac{er_F (1-\alpha)^2K_F m}{4} \beta\Big) S^{max} + \left( \mu_D - f_D   +  \dfrac{(1-\alpha)r_F}{4}\right)e(1-\alpha)K_F, \\
& \leq \cI - \Big(\mu_D -f_D - \dfrac{er_F (1-\alpha)^2K_F m}{4} \beta\Big) S^{max} + \left( \mu_D - f_D   +  \dfrac{(1-\alpha)r_F}{4}\right)e(1-\alpha)K_F.
\end{align*}

Using the definition of $S^{max}$, we finally obtain $(\nabla G_1 \cdot g)|_{S = S^{max}} \leq 0$. 

The two others inequalities are straightforward to obtain. We have:
\begin{align*}
(\nabla G_2 \cdot g)|_{F_W = F_W^{max}} &= r_F (1+\beta H_W) \left(1 - \dfrac{K_F (1-\alpha)}{K_F (1-\alpha)}\right)K_F (1-\alpha)  - \lfw H_W K_F (1-\alpha), \\
(\nabla G_2 \cdot g)|_{F_W = F_W^{max}} & = - \lfw H_W K_F (1-\alpha), \\
(\nabla G_2 \cdot g)|_{F_W = F_W^{max}} & \leq 0.
\end{align*}

The computations for $(\nabla G_3 \cdot g)|_{H_W = H_W^{max}}$ give:

\begin{align*}
(\nabla G_3 \cdot g)|_{H_W = H_W^{max}} &= m_D (S - eF_W) - m_W H_W^{max}, \\
(\nabla G_3 \cdot g)|_{H_W = H_W^{max}} &= m_D (S - S^{max} -  eF_W), \\
(\nabla G_3 \cdot g)|_{H_W = H_W^{max}} & \leq 0.
\end{align*}

We have shown that $(\nabla G_1 \cdot g)|_{S = S^{max}} \leq 0$, $(\nabla G_2 \cdot g)|_{F_W = F_W^{max}} \leq 0$ and $(\nabla G_3 \cdot g)|_{H_W = H_W^{max}} \leq 0$ in  $\Omega_S$.  According to \cite{smoller_shock_1994}, this proves that $\Omega_S$ is an invariant region for system \eqref{equation:SFWHW}.

This also shows that the set  $$\Big\{\Big(H_D, F_W, H_W \Big) \in \mathbb{R}^3  \Big|H_D + eF_W \leq S^{max}, F_W \leq F_W^{max}, H_W \leq H_W^{max} \Big\}$$ is invariant for system \eqref{equation:HDFWHW}.

Moreover, for any point $y \in \partial (\mathbb{R}_+^3)$, the vector field defined by $f(y)$ is either tangent or directed inward. Then, $\Omega$ is an invariant region for equations \eqref{equation:HDFWHW}. 

\end{proof}

The next proposition shows that system \eqref{equation:HDFWHW} defines a dynamical system on $\Omega$.

\begin{prop}
\vquatre{When $0\leq \beta < \beta^*$,} system \eqref{equation:HDFWHW} defines a dynamical system on $\Omega$, that is, for any initial condition $(t_0, y_0)$ with $t_0 \in \mathbb{R}$ and $y_0 \in \Omega$, it exists a unique solution of equations \eqref{equation:HDFWHW}, and this solution is defined for all $t \geq t_0$.
\end{prop}

\begin{proof}
The right hand side of equations \eqref{equation:HDFWHW} defines a function $f(y)$ which is of class $\mathcal{C}^1$ on $\mathbb{R}^3$. The theorem of Cauchy-Lipschitz ensures that model \eqref{equation:HDFWHW} admits a unique solution, at least locally, for any given initial condition, see \cite{walter_ordinary_1998}.
Since \vquatre{$0\leq \beta < \beta^*$,  $\Omega$ is an invariant and compact region according to proposition \ref{prop:invariantRegion}}. Therefore, the solutions with initial condition in $\Omega$ are bounded. Based on uniform boundedness, we deduce that solutions of system \eqref{equation:HDFWHW} with initial condition on $\Omega$ exist globally, for all $t\geq t_0$. Therefore, \eqref{equation:HDFWHW} defines a forward dynamical system on $\Omega$.
\end{proof}

\section{Analysis Using a Quasi-Steady State Approach} \label{sec:qssa}
The model we propose combines demographic and migration processes, which have different timescales, one being much faster than the other. Based on this observation, we develop a quasi-steady state approach. We use the Tikhonov theorem (see theorem \ref{theorem:Tikhonov} page \pageref{theorem:Tikhonov}) to derive a reduced system. The reduced system is analyzed and then compared to the non-reduced one.

\subsection{Derivation of a Reduced System}

Considering the characteristic times of humans and wildlife, the time scale used throughout the study is year. However, the human displacements between the domestic and wild area are daily processes, and are faster than the demographic processes. In consequence, we write $$m_W = \dfrac{\mW}{\epsilon}, \quad m_D = \dfrac{\mD}{\epsilon}$$ where $\epsilon$ is the conversion parameter between year and day. 

Taking into consideration this difference of temporality of demographic and migration processes, one could used a quasi steady state approach to analyze the system \eqref{equation:HDFWHW}. In this purpose, we rewrite the system under the form:

\begin{equation}
\def\arraystretch{2}
\left\{ 
\begin{array}{l}
\dfrac{dH_D}{dt}= \cI + e\lfw F_WH_W + (f_D - \mu_D) H_D + \dfrac{1}{\epsilon}(\mW H_W - \mD H_D), \\
\dfrac{dF_W}{dt} = r_F(1- \alpha) (1+ \beta H_W) \left(1 - \dfrac{F_W}{K_F(1-\alpha)} \right) F_W - \lfw F_WH_W, \\
\dfrac{dH_W}{dt}= - \dfrac{1}{\epsilon}(\mW H_W - \mD H_D).
\end{array} \right.
\label{equation:HDFWHW, fast}
\end{equation}

%
%
%

\begin{prop}
For any initial condition $(H_{D}^0, F_{W}^0, H_{W}^0)$, the solution $\Big(H_{D, \epsilon}, F_{W,  \epsilon}, H_{W,  \epsilon} \Big)$ of system \eqref{equation:HDFWHW, fast} is such that for any time $T > 0$, 
\begin{equation*}
\def\arraystretch{2}
\left\lbrace \begin{array}{l}
\lim\limits_{ \epsilon \rightarrow 0}{H_{D, \epsilon}(t)} = H_D(t), \quad t \in [0,T], \\
\lim\limits_{ \epsilon \rightarrow 0} F_{W,  \epsilon}(t) = F_W(t),\quad t \in [0,T], \\
 \lim\limits_{ \epsilon \rightarrow 0} H_{W,  \epsilon}(t) = \dfrac{\mD}{\mW} H_D(t), \quad  t\in (0,T] \\
\end{array} \right.
\end{equation*}
where $(H_D, F_W)$ is the solution of the following reduced system
\begin{equation}
\def\arraystretch{2}
\left\lbrace \begin{array}{l}
\dfrac{dH_D}{dt} = \dfrac{\cI}{1+ \dfrac{\mD}{\mW}} + \dfrac{f_D - \mu_D}{1+\dfrac{\mD}{\mW}} H_D + \dfrac{e}{1+\dfrac{\mD}{\mW}} \dfrac{\mD}{\mW} \lfw F_WH_D, \\
\dfrac{dF_W}{dt} = (1-\alpha) (1+\beta \dfrac{\mD}{\mW} H_D) r_F \left(1 - \dfrac{F_W}{(1-\alpha)K_F} \right) F_W - \dfrac{\mD}{\mW} \lfw F_W H_D.
\end{array} \right.
\label{equation:HDFW}
\end{equation} 
\end{prop}

\begin{proof}
We will apply the Tikhonov theorem (see theorem \ref{theorem:Tikhonov}, page \pageref{theorem:Tikhonov}) on system \eqref{equation:HDFWHW, fast}. Following \cite{banasiak_methods_2014} and to avoid division by zero, we rewrite the system with the variable $H = H_D + H_W$. It gives:

\begin{equation*} 
\def\arraystretch{2}
\left\{ 
\begin{array}{l}
\dfrac{dH}{dt}= \cI + e\lfw F_WH_W + (f_D - \mu_D)(H - H_W), \\
\dfrac{dF_W}{dt} = r_F(1- \alpha) (1+ \beta H_W) \left(1 - \dfrac{F_W}{K_F(1-\alpha)} \right) F_W - \lfw F_W H_W, \\
\epsilon \dfrac{dH_W}{dt}= m_D H - (m_D + m_W) H_W.
\end{array} \right.
\label{equation:HFWHW} 
\end{equation*}

This system 
is under the appropriate form to apply the Tikhonov's theorem with $x = (H, F_W)$, $y = H_W$,  

$$f(t,x,y,\epsilon) = \begin{bmatrix}
cI + e\lfw F_W H_W + (f_D - \mu_D) (H - H_W), \\
r_F(1- \alpha) (1+ \beta H_W) \left(1 - \dfrac{F_W}{K_F(1-\alpha)} \right) F_W - \lfw F_W H_W
\end{bmatrix}  $$
and $g(t,x,y,\epsilon) = \mW \Big(m H - (1 + m)H_W \Big) $.

The functions $f(t, \cdot, \cdot, \epsilon)$ and $g(t, \cdot, \cdot, \epsilon)$ are continuously differentiable in arbitrary $\mathcal{\bar{U}}$ and $\mathcal{V}$ intervals. Therefore, the first assumption of the Tikhonov's theorem is satisfied. The equation $ g(t,x,y,0) =  0$ admits for solution $H_W = \dfrac{m}{1+m}H$. It is clearly continuous as a function of $t, H$ and isolated. The second assumption of Tikhonov's theorem is also satisfied.

The auxiliary equation is given by
\begin{equation*}
\dfrac{d \tilde{H_W}}{d \tau} = \mW \Big(m H(t) - (1 + m)\tilde{H_W} \Big)
\end{equation*}

for which the fixed point $\tilde{H_W} = \dfrac{m}{1+m}H(t)$ is globally uniformly asymptotically stable with respect to $H$ and $t$. Therefore, the third and fifth assumptions of Tikhonov's theorem hold true.

The reduced equation writes:
\begin{equation} \label{equationHFW}
\def\arraystretch{2}
\left\lbrace \begin{array}{l}
\dfrac{dH}{dt}= \cI + e\lfw F_W\dfrac{m}{1+m}H + (f_D - \mu_D)\dfrac{H}{1+m}, \\
\dfrac{dF_W}{dt} = r_F(1- \alpha) (1+ \beta \dfrac{m}{1+m}H ) \left(1 - \dfrac{F_W}{K_F(1-\alpha)} \right) F_W - \lfw F_W \dfrac{m}{1+m}H.  \\
\end{array} \right.
\end{equation}

Using the derivability of the right hand side on $\Omega$, and the fact that those derivatives are bounded on $\Omega$, it is straightforward to show that the fourth assumption of Tikhonov's theorem is satisfied.

Consequently, we can claim that the solutions 
$\Big(H_{ \epsilon}, F_{W,  \epsilon}, H_{W,  \epsilon} \Big)$  of system \eqref{degenerateSystem} satisfy:
\begin{equation*}
\def\arraystretch{2}
\left\lbrace \begin{array}{l}
\lim\limits_{ \epsilon \rightarrow 0}{H_{ \epsilon}(t)} = H(t), \quad t \in [0,T], \\
\lim\limits_{ \epsilon \rightarrow 0} F_{W,  \epsilon}(t) = F_W(t),\quad t \in [0,T], \\
 \lim\limits_{ \epsilon \rightarrow 0} H_{W,  \epsilon}(t) = \dfrac{m}{1+m}H(t), \quad  t\in (0,T], \\
\end{array} \right.
\end{equation*}
for any $T > 0$, where $(H, F_W)$ are the solutions to \eqref{equationHFW}. 
Finally, returning to the original variables $H_D, F_W$ and $H_W$, allows us to conclude the proof.
\end{proof}

The following parts are dedicated to the study of the long term dynamics of system \eqref{equation:HDFW}. The study is separated in two cases: a first one when there is no immigration ($\cI = 0$), and a second one when $\cI > 0$.

\subsection{Reduced System: Local Analysis Without Immigration}

On the following, we study the existence and stability of equilibria of model \eqref{equation:HDFW}, when $\cI = 0$.

\begin{prop}
\label{prop:equilibre, cI=0}
When $\cI = 0$, the following results hold true.
\begin{itemize}
\item System \eqref{equation:HDFW} admits a trivial equilibrium $TE = \Big(0,0\Big)$ and a fauna-only equilibrium $EE^{F_W} = \Big(0, (1-\alpha)K_F \Big)$ that always exist.

\item When
$$
\mathcal{N}_{\cI = 0} := \dfrac{m e \lfw (1-\alpha)K_F}{\mu_D - f_D} >1,
$$ 
then system \eqref{equation:HDFW} admits a unique coexistence equilibrium $EE^{HF_W}_{\cI = 0} = \Big(H^*_{D, \cI = 0}, F^*_{W, \cI = 0}\Big)$ \\ 
where

$$F^*_{W, \cI = 0} = \dfrac{\mu_D - f_D}{\lfw m e},
\quad 
H^*_{D, \cI = 0} = \dfrac{(1-\alpha)r_F\Big(1 - \dfrac{F^*_{W, \cI = 0}}{K_F(1-\alpha)} \Big)}{m\left(\lfw - \beta (1-\alpha) r_F + \beta r_F  \dfrac{F^*_{W, \cI = 0}}{K_F}\right)}.
$$
\end{itemize}
\end{prop}

\begin{remark}\label{remark:existence1}
Straightforward observations show that the equilibrium of the non-reduced system \eqref{equation:HDFWHW} are the same than the ones of system \eqref{equation:HDFW}, with the same condition of existence and with $H_W^* = m H_D^*$.
\end{remark}

\begin{proof}
To derive the equilibrium of system \eqref{equation:HDFW}, we equate the rates of change of the variables to 0. Therefore, an equilibrium satisfies the system of equations:
\begin{equation}\label{equation:system-equilibre, cI=0}
\def\arraystretch{2}
\left\lbrace \begin{array}{cll}
 e \lfw m F_W^* + f_D - \mu_D = 0& \mbox{or} & H_D^* = 0,\\
m H_D^*\Big(\lfw - (1-\alpha)r_F \beta + \dfrac{r_F \beta}{K_F}F_W^* \Big) + r_F \dfrac{F_W^*}{K_F} - (1-\alpha)r_F= 0& \mbox{or} & F^*_W = 0.
\end{array} \right.
\end{equation}
When $H_D^*=0$ and $F_W^*=0$, we recover the trivial equilibrium $TE = \Big(0,0\Big)$. When $H_D^*=0$ and $F_W^*\neq0$, we obtain the fauna-only equilibrium $EE^{F_W} = \Big(0, K_F(1-\alpha)\Big)$. Finally, when $H_D^*\neq0$ and $F_W^*\neq0$, direct computations lead to a unique set of values given by:
$$F^*_{W} = \dfrac{\mu_D - f_D}{\lfw m e},
\quad 
H^*_{D} = \dfrac{(1-\alpha)r_F\Big(1 - \dfrac{F^*_{W}}{K_F(1-\alpha)} \Big)}{m\left(\lfw - \beta (1-\alpha) r_F + \beta r_F  \dfrac{F^*_{W}}{K_F}\right)}.$$

Those values are biologically meaningful if $F_W^* \leq (1-\alpha) K_F$ and if $H_D^*$ is positive. The first inequality gives the constraint $\dfrac{\mu_D - f_D}{\lfw m e} \leq (1-\alpha)K_F$. From now, we assume it. The numerator of $H^*_{D}$ is non negative, and positive when $\dfrac{\mu_D - f_D}{\lfw m e} < (1-\alpha)K_F$. We need to check the sign of its denominator, which has to be positive. We have:

\begin{align*}
\lfw - \beta (1-\alpha) r_F + \beta r_F  \dfrac{F^*_{W}}{K_F} &= \lfw\Big(1 - \dfrac{\beta (1-\alpha) r_F}{\lfw} + \beta r_F  \dfrac{\mu_D - f_D}{\lfw^2 m e K_F} \Big), \\
&= \lfw\left(1 - \dfrac{\beta (1-\alpha) r_F}{\lfw}\Big(1 -\dfrac{\mu_D - f_D}{\lfw m e K_F(1-\alpha)} \Big) \right), \\
&\geq 0.
\end{align*}
The last inequality is obtained using $\beta < \beta^*$ (see proposition \ref{propBeta}, page \pageref{propBeta}). Therefore, the equilibrium of coexistence is biologically meaningful if $\dfrac{\mu_D - f_D}{\lfw m e} < (1-\alpha)K_F$, that is $1 < \dfrac{\lfw (1-\alpha)K_F m e}{\mu_D - f_D}$.

\end{proof}

The following proposition assesses the local asymptotic stability of the equilibrium.

\begin{prop}\label{prop:stab 2D, cI=0}
When $\cI =0$, the following results are valid.
\begin{itemize}
\item The trivial equilibrium $TE$ is unstable.
\item The fauna equilibrium $EE^{HF_W}_{\cI = 0}$ is Locally Asymptotically Stable (LAS) if $\mathcal{N}_{\cI = 0} < 1$.
\item When $\mathcal{N}_{\cI = 0} > 1$, the coexistence equilibrium $EE^{HF_W}$ exists and is LAS.
\end{itemize}
\end{prop}

\begin{proof}
To investigate the local stability of the equilibria of the system \eqref{equation:HDFW}, we rely on the Jacobian matrix of the system. It is given by:

\begin{multline}
\mathcal{J}_{QSSA}(H_D, F_W) = \\ \begin{bmatrix}
- \dfrac{\mu_D - f_D}{1+m} + \dfrac{e \lfw m}{1+m}  F_W & \dfrac{e \lfw m}{1+m}  H_D \\
m\left(-\lfw + \beta (1-\alpha) r_F \Big(1- \dfrac{F_W}{(1-\alpha)K_F} \Big) \right) F_W & (1-\alpha) (1+\beta m H_D) r_F \left(1 - \dfrac{2F_W}{K_F} \right) - \lfw m H_D
\end{bmatrix}.
\label{equation:Jqssa}
\end{multline}

\begin{itemize}
\item At equilibrium $TE$, the Jacobian is given by:
\begin{equation*}
\mathcal{J}_{QSSA}(TE) = \begin{bmatrix}
- \dfrac{\mu_D - f_D}{1+m} &0 \\
0 & (1-\alpha)  r_F 
\end{bmatrix}.
\end{equation*}
$(1-\alpha) r_F > 0$ is a positive eigenvalue, an therefore $TE$ is unstable.

\item At equilibrium $EE^{F_W}$, the Jacobian is given by: 
\begin{equation*}
\mathcal{J}_{QSSA}(EE^{F_W}) = \begin{bmatrix}
- \dfrac{\mu_D - f_D}{1+m} + \dfrac{e}{1+m}\lfw m K_F(1-\alpha) &0 \\
- m \lfw (1-\alpha)K_F & -(1-\alpha)  r_F 
\end{bmatrix}.
\end{equation*}
The eigenvalues are $-(1-\alpha)  r_F$ and $- \dfrac{\mu_D - f_D}{1+m} + \dfrac{e}{1+m}\lfw m K_F(1-\alpha)$. They are both negative if $-(\mu_D - f_D) + e\lfw m K_F(1-\alpha) <0 \Leftrightarrow \mathcal{N}_{\cI = 0} < 1$.

\item At equilibrium $EE^{HF_W}_{\cI = 0}$, the Jacobian is given by:

\begin{multline*}
\mathcal{J}_{QSSA}(EE^{HF_W}_{\cI = 0}) = \\
\begin{bmatrix}
0 & \dfrac{e}{1+m} \lfw m H^*_{D, \cI = 0} \\
m\left(-\lfw + \beta (1-\alpha) r_F \Big(1- \dfrac{F^*_{W, \cI = 0}}{(1-\alpha)K_F} \Big) \right) F^*_{W, \cI = 0} & -(1+\beta m H^*_{D, \cI = 0}) r_F \dfrac{F^*_{W, \cI = 0}}{K_F} 
\end{bmatrix}
\end{multline*}

using equilibrium conditions. According to the Routh-Hurwitz criterion (see theorem \ref{theorem:Routh-Hurwitz}, page \pageref{theorem:Routh-Hurwitz} ),  $EE^{HF_W}_{\cI = 0}$ is LAS if the trace of $\mathcal{J}_{QSSA}(EE^{HF_W}_{\cI = 0}) $ is negative and its determinant positive. We have:

\begin{equation*}
\Tr(\mathcal{J}_{QSSA}(EE^{HF_W}_{\cI = 0})) = -(1+\beta m H^*_{D, \cI = 0}) r_F \dfrac{F^*_{W, \cI = 0}}{K_F}
\end{equation*}
that is $\Tr(\mathcal{J}_{QSSA}(EE^{HF_W}_{\cI = 0})) < 0$.

The determinant is given by:

\begin{align*}
\det(\mathcal{J}_{QSSA}(EE^{HF_W}_{\cI = 0})) &=  \dfrac{- m^2 e}{1+m} \lfw \left(-\lfw + \beta (1-\alpha) r_F \Big(1- \dfrac{F^*_{W, \cI = 0}}{(1-\alpha)K_F} \Big) \right) F^*_{W, \cI = 0} H^*_{D, \cI = 0}, \\
&= \dfrac{m^2 e}{1+m} \lfw^2 \left(1 + \dfrac{\beta (1-\alpha) r_F}{\lfw} \Big(1- \dfrac{\mu_D - f_D}{me \lfw(1-\alpha)K_F} \Big) \right) F^*_{W, \cI = 0} H^*_{D, \cI = 0},
\end{align*}

using the equilibrium value. Using proposition \ref{propBeta}, page \pageref{propBeta}, we show that this last expression is positive. Consequently $\det(\mathcal{J}_{QSSA}(EE^{HF_W}_{\cI = 0})) > 0$ and $EE^{HF_W}_{\cI = 0}$ is LAS without condition.
\end{itemize}
\end{proof}

\subsection{Reduced System: Local Analysis With Immigration}
On the following, we study the existence and local stability of equilibria of model \eqref{equation:HDFW}, when $\cI > 0$.

\begin{prop}\label{prop:eq, cI>0}
When $\cI > 0$, the following results hold true.
\begin{itemize}
\item System \eqref{equation:HDFW} admits a human-only equilibrium $EE^{H} =\Big(\dfrac{\cI}{\mu_D - f_D}, 0 \Big)$ that always exists.
\item When 
$$ \mathcal{N}_{\cI >0} :=  \dfrac{r_F(1-\alpha)\Big({\dfrac{\mu_D - f_D}{m\cI}+\beta\Big)}}{\lfw}  > 1,$$
system \eqref{equation:HDFW} has a unique coexistence equilibrium $EE^{HF_W}_{\cI > 0} = \Big(H^*_{D, \cI > 0}, F^*_{W, \cI > 0}\Big)$
where
$$F^*_{W, \cI > 0} = \dfrac{(1-\alpha)K_F}{2}\left(1 - \dfrac{\sqrt{\Delta_F}}{e(1-\alpha)r_F}\right) + \dfrac{\mu_D - f_D + \cI \beta m}{2\lfw m e},$$
$$
H^*_{D, \cI > 0} = \dfrac{(1-\alpha)r_F\Big(1 - \dfrac{F^*_{W, \cI > 0}}{(1-\alpha)K_F} \Big)}{m\left(\lfw - \beta (1-\alpha) r_F + \beta r_F  \dfrac{F^*_{W, \cI > 0}}{K_F}\right)}
$$
and
\begin{multline*}
\Delta_F = \left(e(1-\alpha)r_F - \dfrac{(\mu_D - f_D) r_F}{\lfw m K_F}\right)^2 + \dfrac{\cI \beta r_F}{\lfw K_F} \left(\dfrac{\cI \beta r_F}{\lfw K_F} + 2\dfrac{(\mu_D - f_D) r_F}{\lfw m K_F} + 2e(1-\alpha)r_F \right) + \\ 4\dfrac{er_F}{K_F}  \cI\Big(1 - \dfrac{(1-\alpha)\beta r_F}{\lfw} \Big).
\end{multline*}
\end{itemize} 
\end{prop}

\begin{proof}
An equilibrium of system \eqref{equation:HDFW} satisfies the system of equations:
\begin{equation}\label{systemEquilibre}
\left\lbrace \begin{array}{cll}
\cI + e \lfw m F_W^* H_D^* + (f_D - \mu_D) H_D^* = 0,&&\\
F_W^* - \dfrac{(1-\alpha)K_F}{1 + \beta m H_D^*} \Big(1 - \dfrac{m(\lfw - (1-\alpha)\beta r_F) H^*_D}{(1-\alpha)r_F} \Big) = 0& \mbox{or} & F^*_W = 0.
\end{array} \right.
\end{equation}

The solution of system \eqref{systemEquilibre} when $F_W^* = 0$ is the Human-only equilibrium $EE^{H} = \Big(\dfrac{\cI}{\mu_D - f_D}, 0\Big)$.
In the sequel, we assume that $F_W^* > 0$. In this case, $F^*_W$ is solution of the quadratic equation
\begin{multline}
P_F(X) := X^2 \left(\dfrac{er_F}{K_F} \right) - X \left(e(1-\alpha)r_F + \dfrac{(\mu_D - f_D) r_F}{\lfw m K_F} + \dfrac{\cI \beta r_F}{\lfw K_F} \right) + \\ \left(\dfrac{(\mu_D - f_D)(1-\alpha) r_F}{\lfw m} - \cI\Big(1 - \dfrac{(1-\alpha)\beta r_F}{\lfw} \Big) \right) = 0.
\end{multline}

The polynomial $P_F$ is studied in appendix \ref{section:study of PF}. In particular, we show that $P_F$ admits two real roots $F_1^* $ and $F_2^*$, with $F_1^* < K_F(1- \alpha) < F_2^*$.

To define an equilibrium, $F^*_W$ must be biologically meaningful, that is be positive and lower than $(1-\alpha) K_F$. Therefore, $F_2^*$ is not biologically meaningful and $F_1^*$ defines an equilibrium only if it is positive, \textit{i.e.} only if  $\dfrac{(\mu_D - f_D) r_F}{\lfw m } > \cI\Big(1 - \dfrac{(1-\alpha)\beta r_F}{\lfw} \Big)$ (see proposition \ref{prop:study of PF}, in appendix \ref{section:study of PF}).  In this case, $F_W^*$ is given by:
$$F^*_W = \dfrac{(1-\alpha)K_F}{2}\left(1 - \dfrac{\sqrt{\Delta_F}}{e(1-\alpha)r_F}\right) + \dfrac{\mu_D - f_D + \cI \beta m}{2\lfw m e}.$$

According to the first equation of system \eqref{systemEquilibre}, the value of $H_D^*$ at equilibrium is given by:

$$
H_D^* = \dfrac{\cI}{\mu_D - f_D - e \lfw m F_1^*}.
$$

It is biologically meaningful if it is positive. Since $F_W^* < \dfrac{\mu_D - f_D}{e \lfw m}$ (see proposition \ref{prop:study of PF}), it is the case. Finally, the coexistence equilibrium exists if $\cI\Big(1 - \dfrac{(1-\alpha)\beta r_F}{\lfw} \Big) < \dfrac{(\mu_D - f_D) r_F}{\lfw m } $.
\end{proof}

The following proposition assesses the local asymptotic stability of the equilibrium.

\begin{prop} \label{prop:stab 2D, cI>0}
When $\cI > 0$, the following results are valid.
\begin{itemize}
\item The human-only equilibrium $EE^{H}$ is LAS if $\mathcal{N}_{\cI > 0} < 1$.
\item When $\mathcal{N}_{\cI > 0} > 1$, the coexistence equilibrium $EE^{HF_W}_{\cI > 0}$ exists and is LAS.
\end{itemize}
\end{prop}

\begin{proof}
To investigate the local stability of the equilibria of the system \eqref{equation:HDFW}, we rely on the Jacobian matrix of the system, computed at equation \eqref{equation:Jqssa}.

\begin{itemize}
\item At equilibrium $EE^{H}$, the Jacobian is given by:
\begin{equation*}
\mathcal{J}_{QSSA}(EE^{H}) = \begin{bmatrix}
-\dfrac{\mu_D - f_D}{1 + m} &  \dfrac{e \lfw m \cI}{(1+m)\mu_D - f_D} \\
0 & (1-\alpha)\Big(1+ \dfrac{m \beta \cI}{\mu_D - f_D}\Big)  r_F -  \dfrac{\lfw m  \cI}{\mu_D - f_D}
\end{bmatrix}.
\end{equation*}
The eigenvalues are $-\dfrac{\mu_D - f_D}{1 + m} $ and $(1-\alpha)\Big(1+ \dfrac{\cI}{\mu_D - f_D}\Big)  r_F -  \dfrac{\lfw m  \cI}{\mu_D - f_D}$. They are both negative if $(1-\alpha)\Big(1+ \dfrac{m \beta \cI}{\mu_D - f_D}\Big)  r_F -  \dfrac{\lfw m  \cI}{\mu_D - f_D} <0 \Leftrightarrow \mathcal{N}_{\cI > 0} < 1$.

\item At equilibrium $EE^{HF_W}_{\cI > 0}$, the Jacobian is given by:

\begin{multline*}
\mathcal{J}_{QSSA}(EE^{HF_W}_{\cI > 0}) \\= \begin{bmatrix}
-\dfrac{\cI}{H^*_{D, \cI > 0}} & e \lfw m H^*_{D, \cI > 0} \\
m\left(-\lfw + \beta (1-\alpha) r_F \Big(1- \dfrac{F^*_{W, \cI > 0}}{(1-\alpha)K_F} \Big) \right) F^*_{W, \cI > 0} & -(1+\beta m H_D) r_F \dfrac{F^*_{W, \cI > 0}}{K_F} 
\end{bmatrix}
\end{multline*}

using equilibrium conditions.  According to the Routh-Hurwitz criterion (theorem \ref{theorem:Routh-Hurwitz}, page \pageref{theorem:Routh-Hurwitz}),  $EE^{HF_W}_{\cI > 0}$ is LAS if the trace of $\mathcal{J}_{QSSA}(EE^{HF_W}_{\cI > 0}) $ is negative and its determinant positive. We have:

\begin{equation*}
\Tr(\mathcal{J}_{QSSA}(EE^{HF_W}_{\cI > 0})) = -\dfrac{\cI}{H^*_{D, \cI > 0}} -(1+\beta m H^*_{D, \cI > 0}) r_F \dfrac{F^*_{W, \cI > 0}}{K_F}, 
\end{equation*}
and therefore $\Tr(\mathcal{J}_{QSSA}(EE^{HF_W}_{\cI > 0})) < 0$.
The determinant is given by:
\begin{align*}
\det(\mathcal{J}_{QSSA}(EE^{HF_W}_{\cI > 0})) &= \dfrac{\cI}{H_D^*} \dfrac{1 + \beta m H_D^*}{K_F} r_F F_W^* + m^2 e \lfw \left(\lfw - \beta(1-\alpha)r_F + r_F \beta\dfrac{F_W^*}{ K_F} \Big) \right) F_W^* H_D^*.
\end{align*}

According to proposition \ref{prop:study of PF}, page \pageref{prop:study of PF}, we know that $(1-\alpha)K_F - \dfrac{\lfw K_F}{\beta r_F} < F_W^*$. Consequently $\det(\mathcal{J}_{QSSA}(EE^{HF_W}_{\cI > 0})) > 0$ and $EE^{HF_W}_{\cI > 0}$ is LAS without condition.
\end{itemize}
\end{proof}

\subsection{Reduced System: Global Stability Results}
This subsection is valid for $\cI \geq 0$. On it, we show that the system \eqref{equation:HDFW} does not admit any limit cycle and that the equilibrium which are LAS are Globally Asymptotically Stable (GAS).

\begin{prop} \label{prop:no limit cycle, 2D}
System \eqref{equation:HDFW} admits no limit cycle on $\Omega$.
\end{prop}

\begin{proof}
We will use the Bendixson-Dulac criterion (see theorem \ref{theorem:Dulac}, page \pageref{theorem:Dulac}) with the function $\phi(H_D, F_W) = \dfrac{1}{H_D F_W}$. We note $f$ the right hand side of equations \eqref{equation:HDFW}. We have:

\begin{equation*}
(\phi \times f_1)(H_D, F_W) = \dfrac{\cI}{(1+m)H_D F_W} -\dfrac{\mu_D - f_D}{F_W} + \dfrac{e}{1+m}\lfw m
\end{equation*} and therefore

\begin{equation*}
\dfrac{\partial (\phi f_1)}{\partial H_D}(H_D, F_W) = - \dfrac{\cI}{(1+m)F_W H_D^2} \leq 0.
\end{equation*}

On the other hand, we also have:
\begin{equation*}
(\phi \times f_2)(H_D, F_W) = - m \lfw + \dfrac{(1-\alpha) (1+ \beta m H_D) r_F}{H_D} - \dfrac{(1+\beta m H_D) r_F}{H_D} \dfrac{F_W}{K_F}
\end{equation*} and therefore

\begin{equation*}
\dfrac{\partial (\phi f_2)}{\partial F_W}(H_D, F_W) = - \dfrac{(1+\beta m H_D) r_F}{H_D K_F} <0.
\end{equation*}

Consequently, $\dfrac{\partial (\phi f_1)}{\partial H_D} + \dfrac{\partial (\phi f_2)}{\partial F_W} < 0$ and according to the Bendixson-Dulac criterion, the system \eqref{equation:HDFW} does not admit any limit cycle.

\end{proof}

\begin{prop} \label{prop:GAS, 2D}
The equilibrium of system \eqref{equation:HDFW} which are LAS are GAS in $\Omega$.
\end{prop}

\begin{proof}
According to proposition \ref{prop:no limit cycle, 2D}, the reduced system \eqref{equation:HDFW} admits no limit cycle in $\Omega$. Therefore, according to the Poincaré-Bendixon's theorem (see theorem \ref{theorem:Poincaré-Bendixson}, page \pageref{theorem:Poincaré-Bendixson}), the only possible dynamic for the system in $\Omega$, is to converge towards a stable equilibrium. 
\end{proof}

The long term dynamics of system \eqref{equation:HDFW} is summarized in table \ref{table: reduced long term dynamics, I = 0}.

\begin{table}[!ht]
\centering
\def\arraystretch{2}
\begin{tabular}{c|c|c}
$\cI$  & $\mathcal{N}_{I = 0}$ & \\
\hline
\multirow{2}{*}{$=0$}& $ < 1$ &$EE^{F_W}$ exists and is GAS on $\Omega$.  \\
\cline{2-3}
 &  \multirow{1}{*}{$> 1$}  &$EE^{F_W}$ exists and is unstable; $EE^{HF_W}_{\cI=0}$ exists and is GAS on $\Omega$.\\
\hline
\hline
$\cI$  & $\mathcal{N}_{I > 0}$ & \\
\hline
\multirow{2}{*}{$>0$} & $< 1$ &$EE^{H}$ exists and is GAS on $\Omega$. \\
\cline{2-3}
 & \multirow{1}{*}{$> 1$}  &$EE^{H}$ exists and is unstable; $EE^{HF_W}_{\cI>0}$ exists and is GAS on $\Omega$. \\
\end{tabular}
\caption{\centering Long term dynamics of system \eqref{equation:HDFW}.}
\label{table: reduced long term dynamics, I = 0}
\end{table}

\subsection{Numerical Comparison Between the Non-Reduced and the Reduced Systems}
\begin{figure}[!ht]
\centering
\begin{subfigure}{0.45\textwidth}
\centering
\includegraphics[width=1\textwidth]{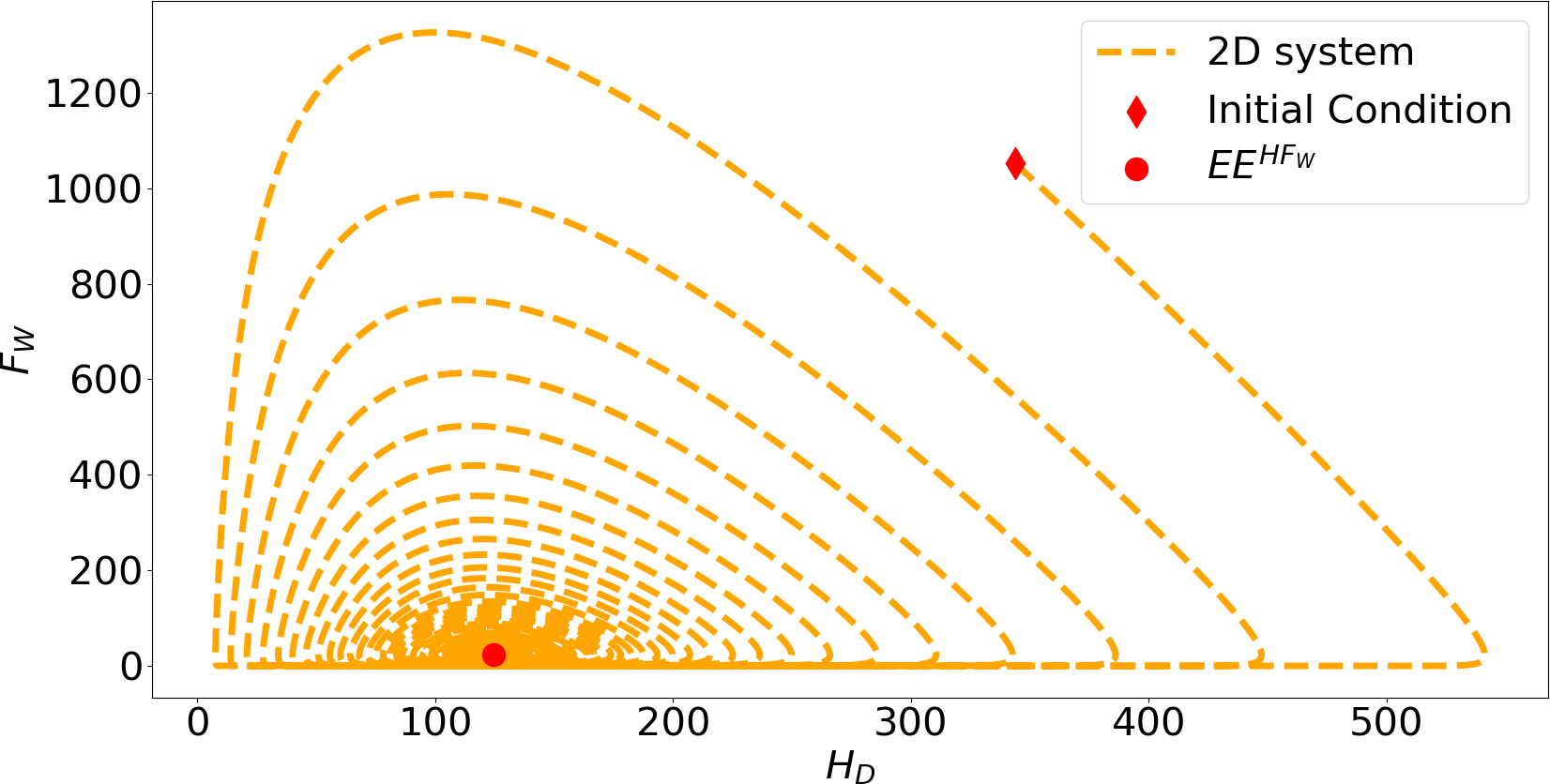}
\caption{}
\label{fig:comparison 2D-3D, a}
\end{subfigure}
\begin{subfigure}{0.45\textwidth}
\centering
\includegraphics[width=1\textwidth]{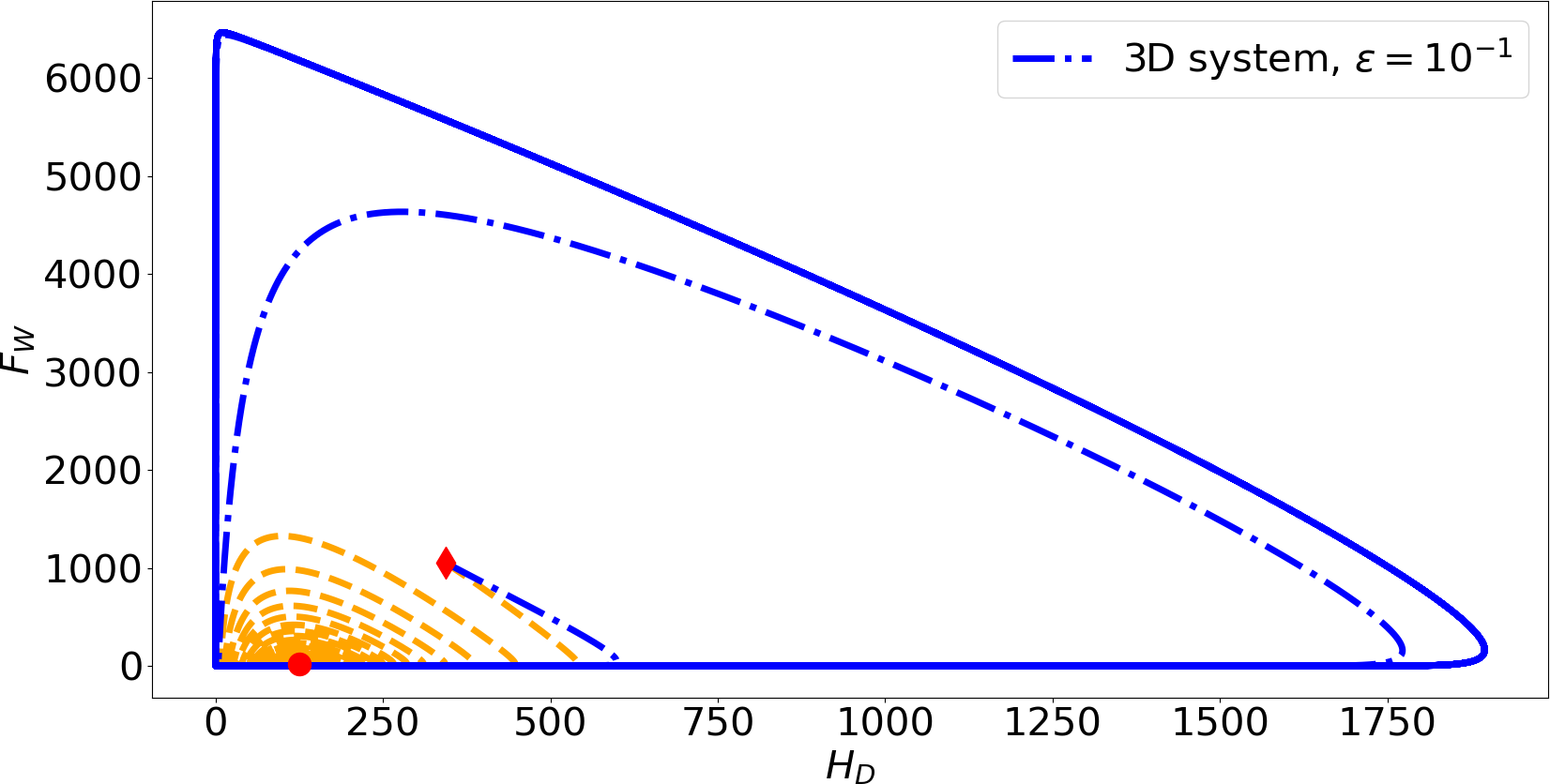}
\caption{}
\label{fig:comparison 2D-3D, b}
\end{subfigure}
\begin{subfigure}{0.45\textwidth}
\centering
\includegraphics[width=1\textwidth]{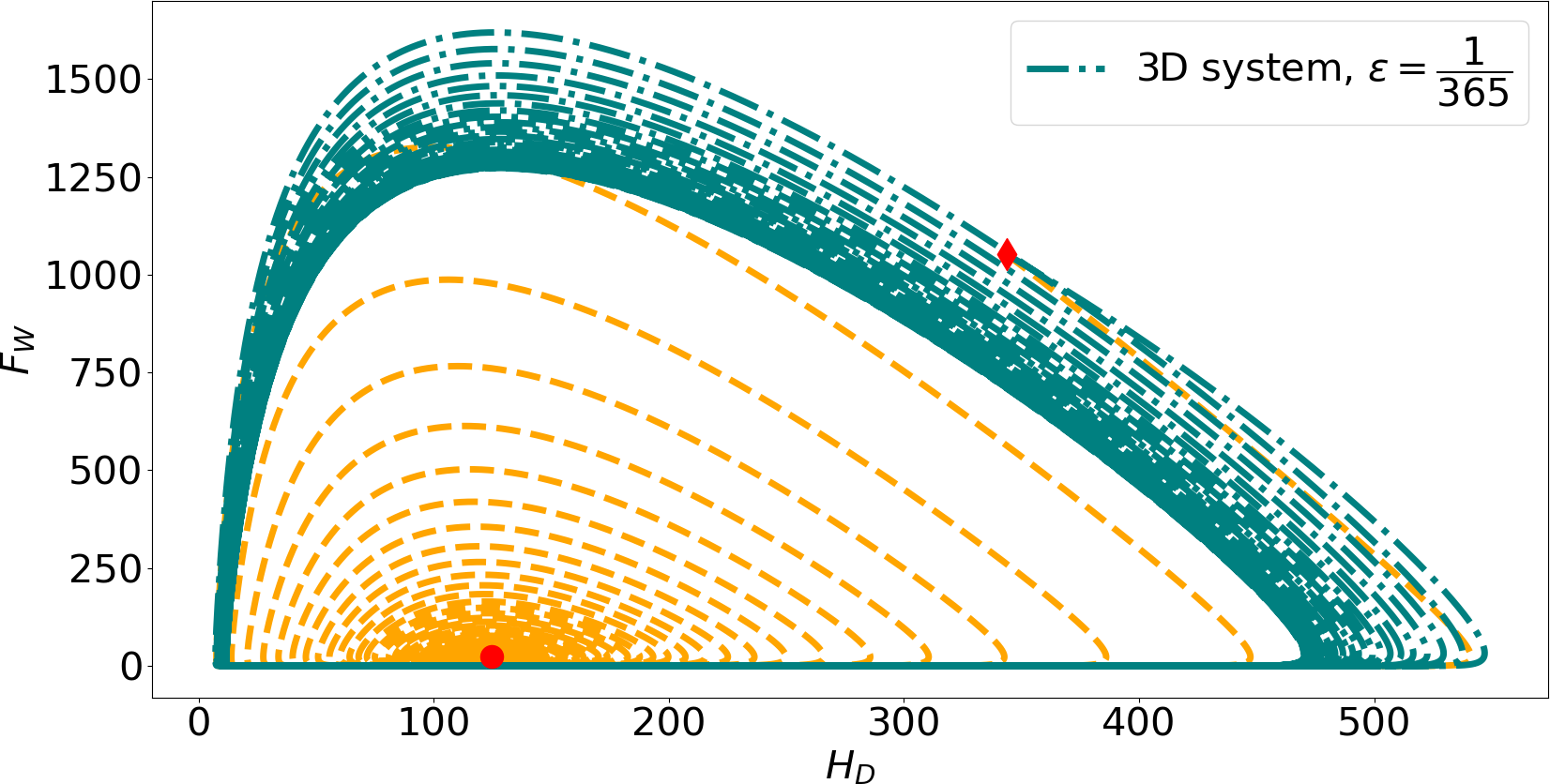}
\caption{}
\label{fig:comparison 2D-3D, c}
\end{subfigure}
\begin{subfigure}{0.45\textwidth}
\centering
\includegraphics[width=1\textwidth]{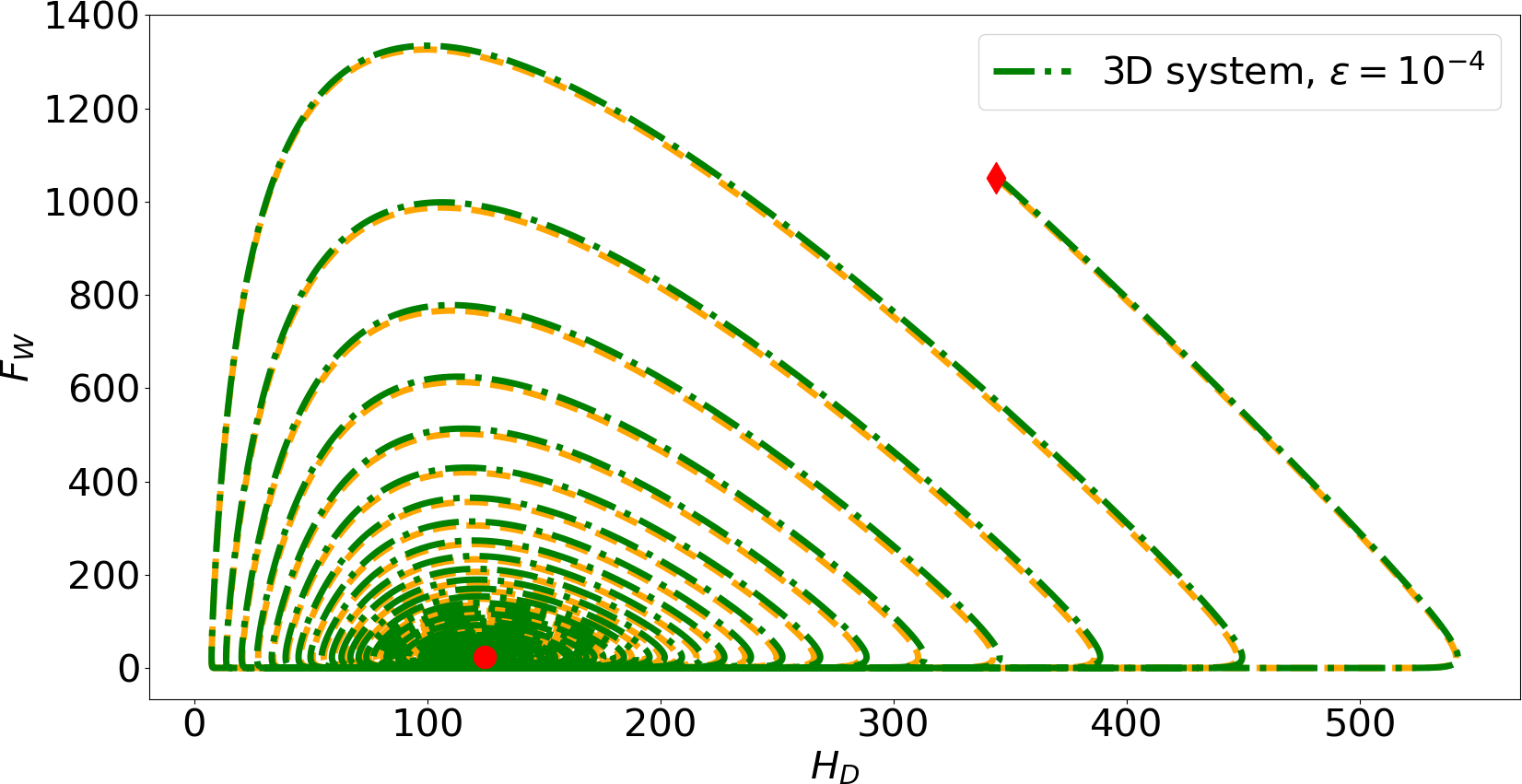}
\caption{}
\label{fig:comparison 2D-3D, d}
\end{subfigure}
\caption{Orbits of systems \eqref{equation:HDFW} and \eqref{equation:HDFWHW} for $\epsilon = 0.1$, $\epsilon = 1/365$ and $\epsilon = 10^{-4}$ in the $H_D - F_W$ plane. In the two first cases, the non-reduced system \eqref{equation:HDFWHW} converges towards a limit cycle while for $\epsilon = 10^{-4}$ it converges towards $EE^{HF_W}$, as does the reduced system \eqref{equation:HDFWHW}. \\
Parameters values: $\cI = 0$, $\beta = 0$, $r_F = 0.6$, $K_F = 7250$, $\alpha = 0.1$, $\lfw = 0.015$, $e = 0.2$, $\mu_D = 0.02$, $f_D = 0.0$, $\mD = 0.0019$, $\mW = 0.0066$.}
\label{fig:comparison 2D-3D}
\end{figure}

In this section, we seek to evaluate the quality of the quasi-steady state approximation. To do this, we compare numerical simulations of the original system \eqref{equation:HDFWHW} and of its approximation, the system \eqref{equation:HDFW}. The simulations are shown in figure \ref{fig:comparison 2D-3D}.

\vspace{0.1cm}

In figure \ref{fig:comparison 2D-3D}, we used parameter values such that $\mathcal{N}_{I=0}=287.17 > 1$. Therefore, according to proposition \ref{prop:GAS, 2D}, the coexistence equilibrium $EE^{HF_W}$ is GAS for system \eqref{equation:HDFW}. This is illustrated in figure \ref{fig:comparison 2D-3D, a}, where the solution of the reduced system \eqref{equation:HDFW}, represented by the dashed orange line, converges towards $EE^{HF_W}$.

In figures \ref{fig:comparison 2D-3D, b},\ref{fig:comparison 2D-3D, c} and \ref{fig:comparison 2D-3D, d} we compare this solution with the solution of the non-reduced system \eqref{equation:HDFWHW} for different values of $\epsilon$ (the conversion parameter between the demographic time scale (year) and the forest-village travel time scale (day)). When $\epsilon$ is small enough ($\epsilon = 10^{-4}$), the original system also converges towards $EE^{HF_W}$, and the orbits are close to each others.
However, for highest values of $\epsilon$ ($\epsilon = 10^{-1}, 1/365 $), the solution of the non-reduced system converges towards a limit cycle around $EE^{HF_W}$.

This difference of behavior does not contradict proposition \ref{prop: equivalentSystem} (which is valid for $\epsilon \rightarrow 0$) but raises questions about the pertinence of the quasi steady state approximation, since $\epsilon = \dfrac{1}{365}$ is precisely the value we should use. 

Therefore, in the following section we propose an analysis of the non-reduced system \eqref{equation:HDFWHW}, using the theory of monotone systems.

\section{Analysis of an Equivalent Competitive System} \label{sec:competitive}
The quasi-steady state approach was a tentative to provide a full and simple analysis of our model. However, as seen above, it does not recover the complete dynamic of the original system \eqref{equation:HDFWHW}. To provide a full analysis of the long term dynamics of our model, we will work on an equivalent system which is competitive. This characteristic enables us to use the theory of monotone systems (see definition \ref{def:monotone}, page \pageref{def:monotone}), and especially, theorem \ref{theorem:Zhu}, page \pageref{theorem:Zhu} (see \cite{zhu_stable_1994}).

\subsection{Derivation of the Equivalent Competitive System}
On this section we derive an equivalent system by a change of variable, and state some of its properties.
\begin{prop} \label{prop: equivalentSystem}
System \eqref{equation:HDFWHW} is equivalent to the system:
\begin{equation}
\def\arraystretch{2}
\left\{ \begin{array}{l}
\dfrac{dh_D}{dt}= \cI + e\lfw h_W f_W + (f_D - \mu_D) h_D - m_D h_D - m_W h_W, \\
\dfrac{df_W}{dt} = (1-\alpha)(1 - \beta h_W) r_F \left(1 + \dfrac{f_W}{K_F(1-\alpha)} \right) f_W + \lfw f_W h_W, \\
\dfrac{dh_W}{dt}= -m_D h_D - m_W h_W. 
\end{array} \right.
\label{equation:hdfwhw}
\end{equation}
We note $y_{compet} = (h_D, f_W, h_W)$ the equivalent variable, $f_{compet}(y_{compet})$ the right hand side of \eqref{equation:hdfwhw} and 
$$
\Omega_{compet} = \Big\{\Big(h_D, f_W, h_W \Big) \in \R_+ \times \R_-^2  \Big|h_D - ef_W \leq S^{max},  -F_W^{max} \leq f_W, -H_W^{max} \leq h_W\Big\}
$$ the domain of interest. 
This equivalent system is irreducible and dissipative on $\Omega_{compet}$, and
\begin{itemize}
\item if $\lfw - (1-\alpha)\beta r_F \geq 0$, it is competitive on $\mathcal{D} = \R_+ \times \R_-^2$,
\item if $\lfw - (1-\alpha)\beta r_F < 0$, it is competitive on $$\Big\{(h_D, f_W, h_W)\in \R_+ \times \R_-^2  | f_W \leq K_F\big(\dfrac{\lfw}{\beta r_F}-(1-\alpha)\big)\Big\}.$$ 
\end{itemize}

\end{prop}

\begin{proof}
Following \cite{wang_predator-prey_1997}, we consider the following change of variable: 

\begin{equation}\label{equation:change of variable}
h_D =  H_D, \quad f_W = -F_W, \quad h_W = -H_W.
\end{equation} The system \eqref{equation:HDFWHW} is transformed into:

\begin{equation}
\def\arraystretch{2}
\left\{ \begin{array}{l}
\dfrac{dh_D}{dt}= \cI + e\lfw h_W f_W + (f_D - \mu_D) h_D - m_D h_D - m_W h_W, \\
\dfrac{df_W}{dt} = (1-\alpha)(1 - \beta h_W) r_F \left(1 + \dfrac{f_W}{K_F(1-\alpha)} \right) f_W + \lfw f_W h_W, \\
\dfrac{dh_W}{dt}= -m_D h_D - m_W h_W.
\end{array} \right.
\end{equation}
It is clear that $\mathcal{D}$ is a $p$-convex set, in which $f_{compet}$ is analytic. According to proposition \ref{prop:invariantRegion}, the region $\Omega_{compet}$ is an invariant region for system \eqref{equation:hdfwhw}, and a compact subset of $\mathcal{D}$. This shows that the system \eqref{equation:hdfwhw} is dissipative for initial condition in  $\Omega_{compet}$.

The Jacobian of $f_{compet}$ is given by:

{
\begin{multline}
\mathcal{J}_{f_{compet}}(y_{compet}) = \\ \begin{bmatrix}
f_D -\mu_D - m_D & e \lfw h_W & e \lfw f_W - m_W \\
0 & r_F (1-\alpha)(1-\beta h_W) \Big(1 + \dfrac{2 f_W}{K_F(1-\alpha)}\Big) + \lfw  h_W & \left(\lfw- (1-\alpha)\beta r_F - \beta r_F \dfrac{f_W}{K_F} \right)f_W\\
-m_D & 0 & -m_W
\end{bmatrix}.
\label{equation:jacobianMatrix compet}
\end{multline}}
Therefore, it is clear that system \eqref{equation:hdfwhw} is irreducible. The system is competitive if the non-diagonal term $\lfw f_W - (1-\alpha)\beta r_F f_W - \beta r_F \dfrac{f_W^2}{K_F}$ is non-positive. For $f_W \in (-\infty, 0]$, we have:

\begin{align*}
&\lfw f_W - (1-\alpha)\beta r_F f_W - \beta r_F \dfrac{f_W^2}{K_F} \leq 0, \\
&\Leftrightarrow \beta r_F \dfrac{f_W}{K_F} \leq \lfw - (1-\alpha)\beta r_F .
\end{align*}

Therefore, the system is competitive on $(-\infty, 0]$ if $0 \leq \lfw - (1-\alpha)\beta r_F$. When $\lfw - (1-\alpha)\beta r_F<0$, the previous computations show that the system is competitive only for $f_W \in \Big(-\infty, -K_F(1-\alpha) + \dfrac{K_F \lfw}{\beta r_F}\Big]$. 
\end{proof}

As stated by the previous proposition, when $\lfw - (1-\alpha)\beta r_F < 0$, the equivalent system is not competitive on the whole domain of interest, $\Omega_{compet}$, but only on a subdomain, which we note $\Omega_{compet}^1$. However, this is not so important since $\Omega_{compet}^1$ is an invariant and absorbing set. That is the subject of the next proposition.

\begin{prop} \label{prop:absorbing set}

When $\lfw < r_F(1-\alpha) \beta$, we define 
\begin{multline*}
\Omega_{compet}^1 = \\ \Big\{\Big(h_D, f_W, h_W \Big) \in \R_+ \times \R_-^2  \Big|h_D - ef_W \leq S^{max}, -F_W^{max} \leq f_W \leq -K_F\Big(1-\alpha - \dfrac{ \lfw}{\beta r_F}\Big), -H_W^{max}\leq  h_W\Big\}.
\end{multline*}
The set $\Omega_{compet}^1$ is an invariant region for system \eqref{equation:hdfwhw}, on which it is competitive. Moreover, any solution of system \eqref{equation:hdfwhw} with initial condition in $\Omega \setminus \Omega_{compet}^1$ will enter in $\Omega_{compet}^1$.
\end{prop}

\begin{proof}
We assume $\lfw < {r_F(1-\alpha) \beta}$ and we note $f_W^{max} = -K_F\Big(1-\alpha - \dfrac{ \lfw}{\beta r_F}\Big)$.
We start by showing that $\Omega_{eq, 1}$ is an invariant region. In fact, since we already prove that $\Omega_{compet}$ is an invariant region, we only need to show that 
$\nabla (G \cdot f_{compet})_{|f_W = f_W^{\max}}(y_{compet}) < 0$, for $y_{compet} \in \Omega_{compet}^1$ where $G = f_W - f_W^{\max}$. We have:

\begin{align*}
(\nabla G \cdot y_{compet})_{|f_W = f_W^{max}} &= r_F(1-\alpha)(1-\beta h_W) \left(1 + \dfrac{f_W^{max}}{(1-\alpha) K_F} \right)f_W^{max} + \lfw h_W f^{max}_W, \\
&= \left(r_F(1-\alpha)(1-\beta h_W) \left(1 + \dfrac{-(1-\alpha) K_F + \dfrac{K_F \lfw}{r_F \beta}}{(1-\alpha) K_F}\right) + \lfw h_W \right) f^{max}_W \\
&= \left((1-\beta h_W)\dfrac{\lfw}{\beta} + \lfw h_W \right) f^{max}_W, \\
&= \dfrac{\lfw}{\beta} f_W^{max}, \\
&< 0.
\end{align*}

Therefore, $\Omega_{compet}^1$ is an invariant region. Now, we show that any solution with initial condition in $\Omega \setminus \Omega_{compet}^1$ enter in $\Omega_{compet}^1$. We consider $f_W \in (f_W^{max}, 0)$, and using the same computations than before, we obtain:

\begin{align*}
\dfrac{df_W}{dt} = &r_F(1-\alpha)(1-\beta h_W) \left(1 + \dfrac{f_W}{(1-\alpha) K_F}\right)f_W + \lfw h_W  f_W, \\
& \leq \left(r_F(1-\alpha)(1-\beta h_W) \left(1 + \dfrac{f_W^{max}}{(1-\alpha) K_F}\right) + \lfw h_W  \right) f_W, \\
& \leq \dfrac{\lfw}{\beta} f_W,\\
&< 0
\end{align*}
This means that any solution with initial condition in $\Omega_{compet} \setminus \Omega_{compet}^1$ will enter in $\Omega_{compet}^1$.
\end{proof}

\begin{remark} \label{remark:competitivity}
This proposition allows us to use the competitiveness of the system on the whole domain $\Omega_{compet}$ without consideration for the sign of $\lfw - (1-\alpha)\beta r_F$.
\end{remark}

\subsection{Equilibria Existence}
As state in remark \ref{remark:existence1}, the equilibrium of the system \eqref{equation:HDFWHW} are the same than the ones of system \eqref{equation:HDFW}. Therefore, the equilibrium of system \eqref{equation:hdfwhw} are obtained after applying the change of variable described by equations \eqref{equation:change of variable}.

\begin{prop}
When $\cI = 0$, the following results hold true.
\begin{itemize}
\item System \eqref{equation:hdfwhw} admits a trivial equilibrium $TE = \Big(0,0, 0\Big)$ and an equilibrium $EE^{f_W} = \Big(0, -(1-\alpha)K_F,0 \Big)$ that always exist.

\item When
$$
\mathcal{N}_{\cI = 0} = \dfrac{m e \lfw (1-\alpha)K_F}{\mu_D - f_D} > 1,
$$ 
then system \eqref{equation:hdfwhw} admits a coexistence equilibrium $EE^{hf_W}_{\cI = 0} = \Big(h^*_{D, \cI = 0}, f^*_{W, \cI = 0}, h^*_{W, \cI = 0}\Big)$ \\ 
where $h^*_{D, \cI = 0} = H^*_{D, \cI = 0}$ and $f^*_{W, \cI = 0} = -F^*_{W, \cI = 0}$ are given in proposition \ref{prop:equilibre, cI=0}, and
$$ 
h^*_{W, \cI = 0} = -m h^*_{D, \cI = 0}.
$$
\end{itemize}

When $\cI > 0$, the following results hold true.
\begin{itemize}
\item System \eqref{equation:hdfwhw} admits an equilibrium $EE^{h} = \Big(\dfrac{\cI}{\mu_D - f_D}, 0, -m\dfrac{\cI}{\mu_D - f_D} \Big)$ that always exists.
\item When 
$$ \mathcal{N}_{\cI >0} =\dfrac{r_F(1-\alpha)\Big({\dfrac{\mu_D - f_D}{m\cI}+\beta\Big)}}{\lfw}  > 1,$$
system \eqref{equation:hdfwhw} admits a coexistence equilibrium $EE^{hf_W}_{\cI > 0} = \Big(h^*_{D, \cI > 0}, f^*_{W, \cI > 0}, h^*_{W, \cI > 0}\Big)$
where $h^*_{D, \cI > 0} = H^*_{D, \cI > 0}$ and $f_{W, \cI > 0} = - F^*_{W, \cI > 0}$ are given in proposition \ref{prop:eq, cI>0}, and
$$ 
h^*_{W, \cI > 0} = -m h^*_{D, \cI > 0}.
$$
\end{itemize} 
\end{prop}  

We continue by studying the local stability of the equilibrium.

\subsection{Competitive System: Local Stability}

\subsubsection{Stability Analysis Without Immigration}
\begin{prop}\label{prop:stab, cI=0}

 When $\cI = 0$, the following results are valid.
\begin{itemize}
\item The trivial equilibrium $TE$ is unstable.
\item When $\mathcal{N}_{\cI = 0} < 1$, the equilibrium $EE^{f_W}$ is LAS.
\item When $\mathcal{N}_{\cI = 0} > 1$, the equilibrium $EE^{hf_W}_{\cI =0}$ exists. It is LAS if $\Delta_{Stab} > 0$ where 
\begin{multline}
\Delta_{Stab, \cI =0} = \Big(\mu_D - f_D + m_D - r_F (1-\beta h_W^*) \dfrac{f_W^*}{K_F} + m_W\Big) \times \\ \big( \mu_D  -f_D + m_D + m_W \big) r_F(1 - \beta h_W^*) \dfrac{-f^*_W}{K_F} + 
m_D e \lfw (1- \alpha) r_F \left(1 + \dfrac{f_W^*}{(1- \alpha)K_F}\right) f_W^*,
\label{eq:deltaStab, I=0}
\end{multline}
and unstable when $\Delta_{Stab, \cI =0} < 0.$
\end{itemize}
\end{prop}

\begin{proof}
To prove this proposition, we look at the Jacobian of system \eqref{equation:hdfwhw}, given by equation \eqref{equation:jacobianMatrix compet}.

\begin{itemize}
\item At equilibrium $TE$, we have:
\begin{equation*}
\mathcal{J}(TE) = \begin{bmatrix}
f_D-\mu_D - m_D & 0 &  -m_W \\
0 & r_F(1-\alpha)  &  0\\
-m_D & 0 & -m_W
\end{bmatrix}.
\end{equation*}
$r_F > 0$ is a positive eigenvalue of $\mathcal{J}(TE)$. So, $TE$ is unstable.
\item At equilibrium $EE^{f_W}$, we have
\begin{equation*}
\mathcal{J}(EE^{f_W}) = \begin{bmatrix}
f_D-\mu_D - m_D & 0 & -e\lfw K_F(1-\alpha) - m_W \\
0 & -(1-\alpha)r_F  & -\lfw(1-\alpha)K_F  \\
-m_D & 0 & -m_W
\end{bmatrix}.
\end{equation*}

The characteristic polynomial of $\mathcal{J}(EE^{F_W})$ is given by:
\begin{equation*}
\chi(X) = \big(X +(1-\alpha)r_F\big) \times \left(X^2 - X\Big(f_D - \mu_D - m_D - m_W \Big) + m_W(\mu_D - f_D) - m_D e \lfw K_F(1-\alpha) \right).
\end{equation*}
The equilibrium is LAS if the real part of the roots of $\chi(X)$ are negative. One of the root is $-(1-\alpha)r_F < 0$, and therefore we only need to determine the sign of the real part of the second factor's roots. As the coefficient in $X$ is positive, the sign of the real part of the root is determined by the sign of the constant coefficient.
They have a negative real part if the constant coefficient is positive \textit{i.e.} if $\dfrac{m e \lfw K_F(1-\alpha)}{\mu_D - f_D} < 1 $, and a positive real part if $\dfrac{m e \lfw K_F(1-\alpha)}{\mu_D - f_D} > 1 $. The local stability of $EE^{f_W}$ follows.

\item Now, we look for the asymptotic stability of the coexistence equilibrium $EE^{hf_W}_{\cI=0}$.
We have

\begin{equation*}
\mathcal{J}(EE^{hf_W}_{\cI \geq 0}) = \begin{bmatrix}
f_D -\mu_D - m_D & e \lfw h_W^* & e \lfw f^*_W - m_W \\
0 & r_F(1 - \beta h_W^*)\dfrac{f_W^*}{K_F} & \big( \lfw - (1-\alpha)\beta r_F\big) f_W^* -  \dfrac{r_F\beta}{K_F} (f_W^*)^2 \\
-m_D & 0 & -m_W
\end{bmatrix}.
\end{equation*} 

We will use the Routh-Hurwitz criterion (see theorem \ref{theorem:Routh-Hurwitz}, page \pageref{theorem:Routh-Hurwitz}) to determine the local stability of $EE^{hf_W}$. First, it requires to compute the coefficients $a_i, i =1,2,3$, which have to be positive. Note that the first part of the computations are common with the ones for proving the local stability of equilibrium $EE^{hf_W}_{\cI >0}$. We have:

\begin{subequations}
\begin{align}
a_2 &= - \Tr\Big(\mathcal{J}(EE^{hf_W}_{\cI \geq 0})\Big), \\
 &= -\Big(f_D - \mu_D - m_D + r_F (1 - \beta h_W^*) \dfrac{h_W^*}{K_F} - m_W\Big), \\
 &= \mu_D - f_D + m_D - r_F (1 - \beta h_W^*) \dfrac{f_W^*}{K_F} + m_W
 \label{equation:coefficient a2}
\end{align}
\end{subequations}
that is $a_2>0$. Coefficient $a_0$ is given by:

\begin{subequations}
\begin{align}
a_0 &= -\det\Big(\mathcal{J}(EE^{hf_W}_{\cI \geq 0})\Big), \\
a_0 &= -\Big(\mu_D + m_D -f_D \Big) m_W r_F (1 - \beta h_W^*) \dfrac{f^*_W}{K_F}  - m_D r_F (1 - \beta h_W^*) \dfrac{f_W^*}{K_F}(e\lfw f_W^* - m_W) + \\
\nonumber
&  m_D e \lfw  \left(\lfw - (1-\alpha)\beta r_F  - \dfrac{r_F\beta}{K_F} f_W^* \right)h_W^* f_W^*, \\
a_0 &= -\Big(\mu_D -f_D \Big) r_F m_W (1 - \beta h_W^*) \dfrac{f^*_W}{K_F}  - m_D e\lfw (1 - \beta h_W^*) r_F \dfrac{(f_W^*)^2}{K_F} + \\
\nonumber
&  m_D e \lfw \left(\lfw - (1-\alpha)\beta r_F  - \dfrac{r_F\beta}{K_F} f_W^* \right)h_W^*f_W^*, \\
a_0 &= -\Big(\mu_D -f_D \Big) r_F m_W (1 - \beta h_W^*) \dfrac{f^*_W}{K_F}  - m_D e\lfw (1 - \beta h_W^*) r_F \dfrac{(f_W^*)^2}{K_F} - \\
\nonumber
&  m_D e \lfw (1- \alpha) r_F \left(1 + \dfrac{f_W^*}{(1- \alpha)K_F}\right) f_W^*, \\
a_0 &= -m_D \lfw e r_F (1 - \beta h_W^*) \left(\dfrac{\mu_D -f_D }{e \lfw m} + f_W^*\right) \dfrac{f_W^*}{K_F} - \\
\nonumber
& m_D e \lfw (1- \alpha) r_F \left(1 + \dfrac{f_W^*}{(1- \alpha)K_F}\right) f_W^*,  \\
a_0 &= -m_D \lfw e r_F \left(\dfrac{\mu_D -f_D }{e \lfw m K_F} + 2\dfrac{f_W^*}{K_F} + (1-\alpha) - \dfrac{\beta h_W^*}{K_F} \left(\dfrac{\mu_D -f_D }{e \lfw m} + f_W^*\right) \right) f_W^*.   \label{equation:coefficient a0}
\end{align}
\end{subequations}

When $\cI = 0$, we have:

\begin{equation*}
f_W^* = - \dfrac{\mu_D - f_D}{\lfw m e}.
\end{equation*} 

Using this expression in equation \eqref{equation:coefficient a0}, we obtain:

\begin{equation*}
a_{0, \cI=0} = -m_D \lfw e r_F  (1- \alpha) \left(1 + \dfrac{f_W^*}{(1 - \alpha)K_F}\right) f_W^* 
\end{equation*}
that is $a_{0, \cI=0}>0$. The coefficient $a_1$ is given by:
\begin{subequations}
\begin{align}
a_1 &= - \big( \mu_D  -f_D + m_D) r_F(1 - \beta h_W^*) \dfrac{f^*_W}{K_F} + (\mu_D -f_D + m_D) m_W - r_F(1 - \beta h_W^*) \dfrac{f_W^*}{K_F} m_W + \\ \nonumber &m_D (e\lfw f^*_W - m_W), \\
a_1 &= - \big( \mu_D-f_D + m_D + m_W) r_F(1 - \beta h_W^*) \dfrac{f^*_W}{K_F} + (\mu_D -f_D) m_W  + m_D e\lfw f^*_W, \\
a_1 &= -\big( \mu_D  -f_D + m_D + m_W) r_F(1 - \beta h_W^*) \dfrac{f^*_W}{K_F} + \left(\dfrac{\mu_D -f_D}{e\lfw m} + f_W^*\right) e \lfw m_D . \label{equation:coefficient a1}
\end{align}
\end{subequations}

Again, using the expression of $f^*_W$ in the case where $\cI = 0$, we have:

\begin{equation*}
a_{1, \cI =0} = -\big( \mu_D  -f_D + m_D + m_W) r_F(1 - \beta h_W^*) \dfrac{f^*_W}{K_F}
\end{equation*}
and we do have $a_{1, \cI =0} > 0$.

The first assumption of the Rough-Hurwitz criterion is verified, $a_{i, \cI =0} > 0$ for $i=1,2,3$. Therefore, the asymptotic stability of $EE^{HF_W}_{\cI =0}$ only depends on the sign of $\Delta_{Stab}= a_2 a_1 - a_0$, which has to be positive. 
\end{itemize}
\end{proof}

\begin{prop} \label{prop:stab, cI=beta=0}

When $\beta = 0$, we have the following equivalence:
$$ \Delta_{Stab, \cI = \beta = 0} > 0 \Leftrightarrow \lfw < \lfw^*$$
where 
\begin{multline*}
\lfw^* := \left[m_{W}(\mu_{D}-f_{D})+\big(\mu_{D}-f_{D}+m_{D}+m_{W})^{2}\right] \times \\
 \dfrac{\left(1+\sqrt{1+4\dfrac{(1-\alpha)m_{W}r_{F}\left(\mu_{D}-f_{D}\right)\big(\mu_{D}-f_{D}+m_{D}+m_{W})}{\left[m_{W}\dfrac{\mu_{D}-f_{D}}{e}+\big(\mu_{D}-f_{D}+m_{D}+m_{W})^{2}\right]^{2}}}\right)}{2em_D (1-\alpha) K_F }.
\end{multline*}
\end{prop}

\begin{proof}
See appendix \ref{sec:stab, cI = beta = 0}, page \pageref{sec:stab, cI = beta = 0}.
\end{proof}

\subsubsection{Stability Analysis With Immigration}
Now, we look for the local stability of the equilibrium.
\begin{prop}\label{prop:stab, cI>0} 

When $\cI > 0$, the following results are valid.
\begin{itemize}
\item When $\mathcal{N}_{\cI > 0} < 1$, the equilibrium $EE^{h}$ is LAS.
\item When $\mathcal{N}_{\cI > 0} > 1$, the equilibrium $EE^{hf_W}$  exists. It is LAS if 
$$\Delta_{Stab, \cI > 0} > 0,$$  where 

\begin{multline}
\Delta_{Stab, \cI > 0} = \left(\mu_D -f_D + m_D - (1 - \beta h_W^*)r_F \dfrac{f_W^*}{K_F} + m_W  \right) \times \\ \left(- \big( \mu_D  -f_D + m_D + m_W) r_F(1 - \beta h_W^*) \dfrac{f^*_W}{K_F} + \left(\dfrac{\mu_D -f_D}{e\lfw m} + f_W^*\right) e \lfw m_D \right) + \\
m_D \lfw e r_F \left(\dfrac{\sqrt{\Delta_F}}{er_F} - \dfrac{\cI \beta}{\lfw K_F e} - \dfrac{\beta h_W^*}{K_F} \left(\dfrac{\mu_D -f_D }{e \lfw m} + f_W^*\right)\right)  f^*_{W}.
\label{eq:deltaStab, I>0}
\end{multline}
\end{itemize}
\end{prop}

\begin{proof}
To assess the local stability or instability of the equilibrium, we look at the Jacobian matrix. The Jacobian of the system \eqref{equation:hdfwhw}, noted $\mathcal{J}$, was computed in equation \eqref{equation:jacobianMatrix compet}.

\begin{itemize}
\item At equilibrium $EE^{h}$, we have
\begin{equation*}
\mathcal{J}(EE^{h}) = \begin{bmatrix}
f_D-\mu_D - m_D & -e \lfw \dfrac{m \cI}{\mu_D - f_D} & -m_W \\
0 & r_F(1-\alpha)\left(1-\beta\dfrac{m\cI}{\mu_D - f_D}\right) - \lfw\dfrac{m\cI}{\mu_D - f_D} & 0 \\
-m_D & 0 & -m_W
\end{bmatrix}.
\end{equation*}

The characteristic polynomial of $\mathcal{J}(EE^{h})$ is given by:
\begin{multline*}
\chi(X) = \left(X - r_F(1-\alpha)\left(1+\beta\dfrac{m\cI}{\mu_D - f_D}\right) + \lfw\dfrac{m\cI}{\mu_D - f_D} \right) \times \\ \left(X^2 - X\Big(f_D - \mu_D - m_D - m_W \Big) + m_W(\mu_D - f_D)\right).
\end{multline*}

The constant coefficient of the second factor and its coefficient in $X$ are positive. So, the roots of the second factor have a negative real part. Therefore, only the sign of $r_F(1-\alpha)\left(1+\beta\dfrac{m\cI}{\mu_D - f_D}\right) - \lfw\dfrac{m\cI}{\mu_D - f_D}$ determines the stability of $EE^{h}$. If it is negative, $EE^{h}$ is LAS.

\item In order to determine the asymptotic stability of the coexistence equilibrium $EE^{hf_W}_{\cI > 0}$, we will use the Routh-Hurwitz criterion (see theorem \ref{theorem:Routh-Hurwitz}, page \pageref{theorem:Routh-Hurwitz}). General expressions of the Jacobian matrix and of coefficients $a_i, i = 1,2,3$ were computed in equations \eqref{equation:jacobianMatrix compet}, \eqref{equation:coefficient a2}, \eqref{equation:coefficient a0} and \eqref{equation:coefficient a1}. 

According to \eqref{equation:coefficient a2}, we have:
\begin{equation*}
a_2 =  \mu_D - f_D + m_D - r_F (1 - \beta h_W^*) \dfrac{f_W^*}{K_F} + m_W
\end{equation*}
and we have $a_2>0$. According to \eqref{equation:coefficient a0}, we have:
\begin{equation*}
a_0 = -m_D \lfw e r_F \left(\dfrac{\mu_D -f_D }{e \lfw m K_F} + 2\dfrac{f_W^*}{K_F} + (1-\alpha) - \dfrac{\beta h_W^*}{K_F} \left(\dfrac{\mu_D -f_D }{e \lfw m} + f_W^*\right) \right) f_W^*.
\end{equation*}

Using
\begin{equation*}
f_W^* = -\dfrac{(1-\alpha)K_F}{2}\left(1 - \dfrac{\sqrt{\Delta_F}}{e(1-\alpha)r_F}\right) - \dfrac{\mu_D - f_D + \cI \beta m}{2\lfw m e},
\end{equation*}
we obtain
\begin{equation*}
a_0 = - m_D \lfw e r_F \left(\dfrac{\sqrt{\Delta_F}}{er_F} - \dfrac{\cI \beta}{\lfw K_F e} -  \dfrac{\beta h_W^*}{K_F} \left(\dfrac{\mu_D -f_D }{e \lfw m} + f_W^*\right)\right)  f^*_{W}.
\end{equation*}
Let show that the terms in parentheses are positive. This will give $a_0> 0$.

Using proposition \ref{prop:study of PF}, we know that $0 < \dfrac{\mu_D -f_D }{e \lfw m} + f_W^*$. According to proposition \ref{prop:eq, cI>0}, we have:

\begin{multline*}
\Delta_F = \left(e(1-\alpha)r_F - \dfrac{(\mu_D - f_D) r_F}{\lfw m K_F}\right)^2 + \dfrac{\cI \beta r_F}{\lfw K_F} \left(\dfrac{\cI \beta r_F}{\lfw K_F} + 2\dfrac{(\mu_D - f_D) r_F}{\lfw m K_F} + 2e(1-\alpha)r_F \right) + \\ 4\dfrac{er_F}{K_F}  \cI\Big(1 - \dfrac{(1-\alpha)\beta r_F}{\lfw} \Big)
\end{multline*}
 which gives $\Delta_F > \left(\dfrac{\cI \beta r_F}{\lfw K_F}\right)^2$ therefore $
\dfrac{\sqrt{\Delta_F}}{er_F} - \dfrac{\cI \beta}{\lfw K_F e} > 0$.

According to \eqref{equation:coefficient a1}, coefficient $a_1$ is given by:
\begin{equation*}
a_1 = -\big( \mu_D  -f_D + m_D + m_W) r_F(1 - \beta h_W^*) \dfrac{f^*_W}{K_F} + \left(\dfrac{\mu_D -f_D}{e\lfw m} + f_W^*\right) e \lfw m_D .
\end{equation*}

which is positive, since  $0 < \dfrac{\mu_D - f_D}{e \lfw m} +  f^*_{W} $.

The first assumption of the Rough-Hurwitz criteria is verified, $a_i > 0$ for $i=1,2,3$. Therefore, the local asymptotic stability of $EE^{HF_W}$ only depends on the sign of $\Delta_{Stab}= a_2 a_1 - a_0$, which has to be positive.
\end{itemize}
\end{proof}

\subsection{Competitive System: Global Stability and Existence of Limit Cycle}

In this section, we are interested in the global stability of the equilibrium. The proof of propositions \ref{prop:limitCycle, cI=0} and \ref{prop:limitCycle, cI>0} concerning the equilibrium $EE^{hf_W}_{\cI \geq 0}$ are done using the theorem \ref{theorem:Zhu}, page \pageref{theorem:Zhu}. To prove the global stability of equilibrium $EE^{f_W}$ and $EE^{h}$ (propositions \ref{prop:EEFGAS} and \ref{prop:EEHGAS}) we use the theorem \ref{theorem: monotone GAS}, page \pageref{theorem: monotone GAS}, which involves inequalities between vectors of $\R^3$. These inequalities are considered in the following sense: for $\mathbf{a}, \mathbf{b} \in \R^3$,
\begin{itemize}
\item $\mathbf{a} \leq \mathbf{b} \Leftrightarrow \left\lbrace \begin{array}{l} 
a_1 \leq b_1 \\ -a_2 \leq - b_2 \\ -a_3 \leq - b_3
\end{array} \right.$
\item $\mathbf{a} < \mathbf{b} \Leftrightarrow \mathbf{a} \leq \mathbf{b}$ and $\mathbf{a}\neq \mathbf{b}$
\end{itemize}
We also note $[\mathbf{a},\mathbf{b}] = \lbrace \mathbf{x} \in \mathbb{R}^3 | \mathbf{a} \leq \mathbf{x} \leq \mathbf{b} \rbrace$.

\subsubsection{Without Immigration}
We first investigate the global stability of the system in the absence of immigration, \textit{i.e.} when $\cI = 0$.
\begin{prop}\label{prop:EEFGAS}If $\mathcal{N}_{I =0} \leq 1
$, the equilibrium $EE^{f_W}$ is globally asymptotically stable (GAS) on $\Omega_{compet}$ for system \eqref{equation:hdfwhw}.
\end{prop}

\begin{proof}
We will apply Theorem \ref{theorem: monotone GAS}, page \pageref{theorem: monotone GAS}.
\begin{itemize}
\item First, let assume that $\lfw -(1-\alpha) \beta r_F \geq 0$. In this case, according to proposition \ref{prop: equivalentSystem}, system \eqref{equation:hdfwhw} is competitive on $\mathcal{D} = \R_+\times \R_-^2$. Let $\mathbf{a} = \mathbf{0}_{\R^3}$ and $\mathbf{b} = \Big(S^{max}, -K_F(1-\alpha), -mS^{max}\Big)$. We have $\mathbf{a} < \mathbf{b}$ and $[\mathbf{a}, \mathbf{b}] \subset \mathcal{D}$. 
It is immediate to check that $EE^{f_W} \in [\mathbf{a}, \mathbf{b}]$, and since $\N_{\cI = 0} \leq 1$, it is the only equilibrium of the system. Moreover, $f_{compet}(\mathbf{a}) = \mathbf{0} \geq \mathbf{0}$. We note $f_{compet}(\mathbf{b}) = \Big(f_1, f_2, f_3\Big)$ and we have:
\begin{align*}
f_1 &= \Big(e\lfw K_F(1-\alpha) + m_W\Big)mS^{max} - \Big(\mu_D - f_D + m_D\Big) S^{max}, \\
&= \Big(m e\lfw K_F(1-\alpha) - (\mu_D - f_D) \Big)S^{max}, \\
& \leq 0
\end{align*}
using $N_{I= 0} = \dfrac{m e \lfw (1-\alpha)K_F}{\mu_D - f_D} \leq 1$. 
Moreover $f_2 = \lfw m (1-\alpha) K_F S^{max} > 0$ and $f_3 = 0$. All these prove that $f_{compet}(\mathbf{b}) \leq 0$. Therefore, according to Theorem \ref{theorem: monotone GAS}, $EE^{f_W}$ is GAS on $[\mathbf{a}, \mathbf{b}]$ and since $\Omega_{compet} \subset [\mathbf{a}, \mathbf{b}]$, it is GAS on $\Omega_{compet}$.
\item Now, let assume that $\lfw -(1-\alpha) \beta r_F < 0$. In this case, according to proposition \ref{prop: equivalentSystem}, the system \eqref{equation:hdfwhw} is only competitive on $\Big\{(h_D, f_W, h_W) | 0 \leq h_D, f_W \leq K_F\big(\dfrac{\lfw}{\beta r_F}-(1-\alpha)\big), h_W \leq 0 \Big\}$, which does not contain $\mathbf{0}_{\R^3}$. Therefore, instead of using $\mathbf{a}= \mathbf{0}_{\R^3}$, we use $\mathbf{a} = \Big(0, \dfrac{\lfw}{\beta r_F}K_F -	(1-\alpha)K_F, 0 \Big)$  to show that $EE^{f_W}$ is GAS on $\Omega_{compet}^1 \subset[\mathbf{a}, \mathbf{b}]$. Since $\Omega_{compet}\setminus\Omega_{compet}^1$ is absorbed by $\Omega_{compet}^1$ (see proposition \ref{prop:absorbing set}), $EE^{f_W}$ is GAS on $\Omega_{compet}$.
\end{itemize}
\end{proof}

We now consider the case $\mathcal{N}_{I = 0} > 1$ and use the Zhu and Smith's Theorem, recalled in theorem \ref{theorem:Zhu} page \pageref{theorem:Zhu} to complete the characterization of the system.

\begin{prop}
\label{prop:limitCycle, cI=0}
Assume $\mathcal{N}_{I =0} > 1$. If
\begin{itemize}
\item $\Delta_{stab, \cI =0} > 0$, the equilibrium $EE^{hf_W}$ is GAS on $\Omega_{compet}$ for system \eqref{equation:hdfwhw}.
\item $\Delta_{stab, \cI =0} < 0$, system \eqref{equation:hdfwhw} admits an orbitally asymptotically stable periodic solution.
\end{itemize}
\end{prop}

\begin{remark}
When $\beta = 0$, using proposition \ref{prop:stab, cI=beta=0}, we equivalently have:
\begin{itemize}
\item if $\lfw <  \lfw^*$, then $EE^{hf_W}$ is GAS on $\Omega_{compet}$ for system \eqref{equation:hdfwhw}.
\item if $\lfw  > \lfw^*$, system \eqref{equation:hdfwhw} admits an orbitally asymptotically stable periodic solution.
\end{itemize}
\end{remark}

\begin{proof}
Since $\mathcal{N}_{I =0} > 1$, the system \eqref{equation:hdfwhw} admits a unique positive equilibrium, $EE^{hf_W}$. By applying Theorem \ref{theorem:Zhu} and considering remark \ref{remark:competitivity}, we know that either $EE^{hf_W}$ is GAS, or that it exists an asymptotically stable periodic solution. According to proposition \ref{prop:stab, cI=0}, the condition for stability is precisely $0 < \Delta_{stab, \cI =0}$. 
\end{proof}

\subsubsection{With Immigration}

As with the case $\cI = 0$, we supplement the local results obtained when $\cI > 0$ with global stability results.

\begin{prop}
When $\mathcal{N}_{I > 0} < 1$, the equivalent system \eqref{equation:hdfwhw} is competitive on the whole domain $\Omega_{compet}$.
\end{prop}

\begin{proof}
Since $\mathcal{N}_{\cI > 0} = \dfrac{r_F(1-\alpha)}{\lfw}\dfrac{\mu_D - f_D}{m \cI} + \dfrac{r_F(1-\alpha) \beta}{\lfw}$, we have $\dfrac{r_F(1-\alpha) \beta}{\lfw} < \mathcal{N}_{\cI > 0}$. Therefore, when $\mathcal{N}_{\cI > 0} < 1$,we deduce that $\dfrac{r_F(1-\alpha) \beta}{\lfw} < 1$, which implies, according to proposition \ref{prop: equivalentSystem}, that the equivalent system \eqref{equation:hdfwhw} is competitive on $\Omega_{compet}$.
\end{proof}

\begin{prop}\label{prop:EEHGAS}

When $\mathcal{N}_{\cI > 0} \vtrois{<} 1$, the equilibrium $EE^{h}$ is GAS on $\Omega_{compet}$ for system \eqref{equation:hdfwhw}.
\end{prop}

\begin{proof}
We note $f_W^{\min} := \max\Big(-\dfrac{\mu_D - f_D}{e m \lfw}, -K_F(1-\alpha)\Big)$ and 
$$\Omega_{GAS} :=\Big\{\Big(h_D, f_W, h_W \Big)\in \R_+ \times \R_-^2  \Big|h_D - ef_W \leq S^{max},  f_W^{min} \leq f_W, -H_W^{max}\leq h_W \Big\}.$$ 
To show that $EE^{h}$ is GAS on $\Omega_{compet}$, we proceed in two steps: we first show that $EE^{h}$ is GAS on $\Omega_{GAS}$ and then we show that any solution with initial condition in $\Omega_{compet}$ enters in $ \Omega_{GAS}$. 

\begin{itemize}
\item We use theorem \ref{theorem: monotone GAS}, page \pageref{theorem: monotone GAS} to show that $EE^{h}$ is GAS on $\Omega_{GAS}$. Let $\mathbf{a} = \mathbf{0}_{\R^3}$ and $\mathbf{b} = \Big(h_D^{\max}, f_W^{\min}, -mh_D^{\max} \Big)$, where $h_D^{max} = \max\Big(S^{max}, \dfrac{\cI}{\mu_D - f_D - e m \lfw f^{\min}_W}\Big)$. We have $\mathbf{a} < \mathbf{b}$, $\mathbf{0}_{\R^3} \leq f_{compet}(\mathbf{a})$ and $EE^{h}$ is the only equilibrium in $[\mathbf{a}, \mathbf{b}]$. We note $f(\mathbf{b}) = (f_1, f_2, f_3)$. We have:
\begin{align*}
f_1 &= I + \Big(-e\lfw f_W^{min} + m_W\Big)mh_D^{\max} - \Big(\mu_D - f_D + m_D\Big) h_D^{\max}, \\
&= I + \Big(-e m \lfw f_W^{min} - (\mu_D - f_D) \Big)h_D^{\max}, \\
& \leq 0
\end{align*}
since $\dfrac{\cI}{\mu_D - f_D - e m \lfw f^{\min}_W} \leq h_D^{\max}$. We immediately have $f_3 = 0$, and since $\N_{I>0} \leq 1$ and $-m h_D^{\max} \leq -mS^{max} \leq -m\dfrac{I}{\mu_D -f_D}$, it is straightforward to obtain $-f_2 \leq 0$. All this prove that $f_{compet}(\mathbf{b}) \leq 0$, and according to theorem \ref{theorem: monotone GAS}, $EE^{h}$ is GAS on $ \Omega_{GAS} \subset [\mathbf{a}, \mathbf{b}]$ for system \eqref{equation:hdfwhw}.
\item Now we show that any solution of \eqref{equation:hdfwhw} with initial condition in $\Omega_{compet}$ enters in $\Omega_{GAS}$. 
Let $\Big(h_D^{f_{compet}},f_W^{f_{compet}}, h_W^{f_{compet}}\Big)$ be a solution of system \eqref{equation:hdfwhw}. We introduce the following system:
\begin{equation}
\def\arraystretch{2}
\left\lbrace \begin{array}{l}
\dfrac{dh_D}{dt}= \cI + (f_D - \mu_D) h_D - m_D h_D - m_W h_W, \\
\dfrac{df_W}{dt} = (1-\alpha)(1 - \beta h_W) r_F \left(1 + \dfrac{f_W}{K_F(1-\alpha)} \right) f_W + \lfw f_W h_W, \\
\dfrac{dh_W}{dt}= -m_D h_D - m_W h_W. 
\end{array} \right.
\label{equation:model H GAS}
\end{equation}

We note $g(y)$ its right hand side. For all $y \in \Omega_{compet}$, we have $g(y) \leq f_{compet}(y)$. Since $f_{compet}$ is competitive, we have the following inequality between the solutions of the two systems, with the same initial condition:
\begin{equation}
\Big(h_D^g,f_W^g, h_W^g \Big) \leq \Big(h_D^{f_{compet}},f_W^{f_{compet}}, h_W^{f_{compet}}\Big),
\label{equation:inequality H GAS}
\end{equation}

where we note $\Big(h_D^g,f_W^g, h_W^g \Big)$ the solution of system \eqref{equation:model H GAS}.

System \eqref{equation:model H GAS} can be written like system \eqref{equation: eqVidyasagar} in the Vidyasagar's Theorem (see theorem \ref{theorem:Vidyasagar}, page \pageref{theorem:Vidyasagar}), with $y_1 = (h_D, h_W)$ and $y_2 = (f_W)$. Using this theorem, and since $\N_{I>0} \leq 1$, we obtain that $EE^{h}$ is GAS on $\Omega_{compet}$ for system \eqref{equation:model H GAS} (see appendix \ref{sec:lemma H GAS}, page \pageref{sec:lemma H GAS} for a proof).

In particular this means that $h_W^g$ converges, when $t\rightarrow +\infty$, towards $h_W^* = -m\dfrac{\mu_D - f_D}{\cI}$. 
Therefore, since $\N_{\cI > 0} = r_F(1-\alpha)\dfrac{1 - \beta h_W^* }{-\lfw h_W^* } < 1$, it exists $T>0$ such that for all $t\geq T$, $r_F(1-\alpha) \dfrac{1-\beta h_W^g}{-\lfw h_W^g} < 1$.

Moreover, using inequality \eqref{equation:inequality H GAS} (which holds in the sense given in the beginning of this section), we have that for all $t\geq T$,
$$
r_F(1-\alpha) \dfrac{1-\beta h_W^{f_{compet}}(t)}{-\lfw h_W^{f_{compet}}(t)} < 1.
$$
This gives that for all $t\geq T$, $\dfrac{df_W^{f_{compet}}}{dt}(t) \geq 0$, and $f_W^{f_{compet}}$ converges towards 0. Therefore, it exists $T_1 > 0$, such that $f_W^{min} \leq f_W^{f_{compet}}(T_1) \leq 0 $, and this means that  $\Big(h_D^{f_{compet}},f_W^{f_{compet}}, h_W^{f_{compet}}\Big)(T_1) \in \Omega_{GAS}$.
\end{itemize}
\end{proof}

\vquatre{
\begin{remark}
When $(1-\alpha)K_F \leq \dfrac{\mu_D - f_D}{e m \lfw} $, the first part of the proof shows that $EE^{h}$ is GAS on $\Omega_{compet}$ for $\N_{I>0} \leq 1$, and not only for $\N_{I>0} < 1$.
\end{remark}}

\begin{prop} \label{prop:limitCycle, cI>0}
Assume $\mathcal{N}_{\cI > 0} > 1$. If
\begin{itemize}
\item $\Delta_{Stab} > 0$, the equilibrium $EE^{hf_w}_{\cI >0}$ is GAS on $\Omega_{compet}$ for system \eqref{equation:hdfwhw}.
\item $\Delta_{Stab} < 0$, the system \eqref{equation:hdfwhw} admits an orbitally asymptotically stable periodic solution on $\Omega_{compet}$.
\end{itemize}
\end{prop}

\begin{proof}
When $\mathcal{N}_{I =0} > 1$, system \eqref{equation:hdfwhw} admits a unique positive equilibrium, $EE^{hf_W}$. By applying Theorem \ref{theorem:Zhu} to this system, we know that either $EE^{hf_W}$ is GAS, or it exists a asymptotically stable periodic solution. Since the condition for stability is precisely $0 < \Delta_{stab, \cI > 0}$, the  proposition is proven. 
\end{proof}

\section{Analysis and Numerical Simulation of the Model} \label{sec:ecological}

\subsection{Long Term Dynamics of the Model}

In this section, we presents the long term dynamics of our model (equations \eqref{equation:HDFWHW}). It is obtained by inversing the change of variable done in the previous section (section \ref{sec:competitive}). Consequently, the results are left without any proof. We still briefly recall the following notations:
\begin{itemize}
\item $EE^{F_W}$ stands for a fauna-only equilibrium
\item $EE^{H}$ stands for a human-only equilibrium
\item $EE^{HF_W}_{\cI \geq 0}$ stands for a coexistence equilibrium between human and wild fauna populations. It exists if $\N_{I = 0} > 1$ (respectively $\N_{I > 0} > 1$) when there is no immigration (respectively when there is immigration). The thresholds $\N_{I \geq 0}$ are given by:
$$
\mathcal{N}_{\cI = 0} = \dfrac{m e \lfw (1-\alpha)K_F}{\mu_D - f_D}, \quad \mathcal{N}_{\cI >0} =\dfrac{r_F(1-\alpha)\Big({1 + \dfrac{m\cI}{\mu_D - f_D}\beta\Big)}}{\lfw\dfrac{m\cI}{\mu_D - f_D}}.
$$
\item Two other thresholds, $\Delta_{Stab, \cI = 0}$ and $\Delta_{Stab, \cI > 0}$ are defined by equations \eqref{eq:deltaStab, I=0} and \eqref{eq:deltaStab, I>0}. They are necessary to assess the stability of the coexistence equilibrium.
\end{itemize}

The long term dynamics of \eqref{equation:HDFWHW} is presented in table \ref{table:long term dynamics, I = 0}.
\vtrois{
\begin{table}[!ht]
\centering
\def\arraystretch{2}
\begin{tabular}{c|c|c|c}
$\cI$  & $\mathcal{N}_{I = 0}$ &  $\Delta_{Stab, \cI = 0}$ & \\
\hline
\multirow{4}{*}{$=0$}& $ \leq 1$ & &$EE^{F_W}$ exists and is GAS.  \\
\cline{2-4}
 &  \multirow{3}{*}{$> 1$} & $ >0$ &$EE^{F_W}$ exists and is unstable; $EE^{HF_W}_{\cI=0}$ exists and is GAS.\\
 \cline{3-4}
 &  &\multirow{2}{*}{$ <0 $} & $EE^{F_W}$ and $EE^{HF_W}_{\cI=0}$ exist and are unstable ; there is an asymptotically \\
&  & &  stable periodic solution around $EE^{HF_W}_{\cI=0}$. \\
\hline
\hline
$\cI$  & $\mathcal{N}_{I > 0}$ &  $\Delta_{Stab, \cI > 0}$ & \\
\hline
\multirow{3}{*}{$>0$} & $<1$ & &$EE^{H}$ exists and is GAS on $\Omega$ \\
\cline{2-4}
 & \multirow{3}{*}{$> 1$}  & $>0$ &$EE^{H}$ exists and is unstable; $EE^{HF_W}_{\cI>0}$ exists and is GAS on $\Omega$ \\
 \cline{3-4}
 & & \multirow{2}{*}{$ < 0$} &$EE^{H}$ and $EE^{HF_W}_{\cI>0}$ exist and are unstable ; there is an asymptotically \\
 & & &  stable periodic solution around $EE^{HF_W}_{\cI>0}$.
\end{tabular}
\caption{\centering Long term dynamics of system \eqref{equation:HDFWHW}.}
\label{table:long term dynamics, I = 0}
\end{table}
}

The analytical results and the long term dynamics of our model are interpreted in section \ref{subsec:ecological}. Before doing so, we found it interesting to revisit the quasi-steady state approximation that we presented in section \ref{sec:qssa}. In the sequel, a comparison between the reduced (system \eqref{equation:HDFW}) and the non-reduced system (system \eqref{equation:HDFWHW}) is carried out in terms of their long term dynamics.

\subsubsection{Long Term Dynamics when $\epsilon \rightarrow 0$} \label{subsec:ecological}
In this section, we deal with system \eqref{equation:HDFWHW} when $\epsilon \rightarrow 0$. It is not obvious that the system has the same long term behavior than its approximation \eqref{equation:HDFW} since proposition \ref{prop: equivalentSystem} is only valid for finite time interval (see \cite{auger_aggregation_2008} or section 4.4.3 in \cite{banasiak_methods_2014} for a counter example). However, in the case of this model, the dynamics are the same, as stated by the following proposition.

\begin{prop}
Since 
$$
\lim\limits_{ \epsilon \rightarrow 0}{\Delta_{Stab, \cI \geq 0}} > 0,
$$
the non-reduced system \eqref{equation:HDFWHW} and the reduced system \eqref{equation:HDFW} have the same long term dynamics when $\epsilon \rightarrow 0$.
\end{prop}
\begin{proof}
We first prove that $\lim\limits_{ \epsilon \rightarrow 0}{\Delta_{Stab, \cI \geq 0}} > 0$. We can note that for $\cI \geq 0$, the value at equilibrium $EE^{HF_W}$ does not depend on $\epsilon$. However, according to \eqref{equation:coefficient a2}, \eqref{equation:coefficient a0} and \eqref{equation:coefficient a1} the Routh-Hurwitz coefficient depend on it. They are given by (the expression are valid for $\cI \geq 0$):

\begin{align*}
a_2  &= \mu_D - f_D + \dfrac{\tilde m_D}{\epsilon} + (1+\beta H_W^*)r_F \dfrac{F_W^*}{K_F} + \dfrac{\tilde m_W}{\epsilon} \\
a_0 &= e \lfw \dfrac{\tilde m_D}{\epsilon} r_F \left(\dfrac{\mu_D -f_D }{e \lfw m K_F} - 2\dfrac{F_W^*}{K_F} + (1-\alpha) + \dfrac{\beta H_W^*}{K_F} \left(\dfrac{\mu_D -f_D }{e \lfw m} - F_W^*\right) \right) F_W^* \\
a_1 &= \Big( \mu_D  -f_D + \dfrac{\tilde m_D}{\epsilon} + \dfrac{\tilde m_W}{\epsilon} \Big) r_F(1+ \beta H_W^*) \dfrac{F^*_W}{K_F} + \left(\dfrac{\mu_D -f_D}{e\lfw m} - F_W^*\right) e \lfw \dfrac{\tilde m_D}{\epsilon}.
\end{align*}

Therefore, direct computations lead to:

\begin{equation*}
\Delta_{Stab} := a_2 a_1 - a_0 = \dfrac{C_2}{\epsilon^2} + \dfrac{C_1}{\epsilon} + C_0,
\end{equation*}
with $C_2 > 0$.

Therefore, when $\epsilon \rightarrow 0$, $\Delta_{Stab} > 0$. The rest of the proposition is immediate (see tables \ref{table: reduced long term dynamics, I = 0} and \ref{table:long term dynamics, I = 0}).
\end{proof}

\subsection{Ecological Interpretation of the Analytical Results} 
In this section, we give an ecological interpretation of the different results we obtained. We are specially interested in the role played by the parameters of hunting rate $\lfw$ and the level of anthropisation, $\alpha$, which traduce the possibility of over-hunting and the deforestation. 

The mathematical analysis shows that the long term dynamics is totally different with or without immigration.

\begin{itemize}
\item When no immigration is considered, the system always admits a fauna only equilibrium, $EE^{F_W}$, and may admit a coexistence equilibrium between wild fauna and human $EE^{HF_W}$. This equilibrium exists only if $1 < \N_{\cI = 0}$, which can be rewritten as a condition on the hunting rate; the coexistence equilibrium exists only if 
$$
\lambda_{F, \cI = 0}^{\min} = \dfrac{\mu_D - f_D}{(1-\alpha) e m K_F} < \lfw,
$$ that is if the hunting rate is sufficiently high. 
This is because the human population relies on bush meat for subsistence. We can notice that the minimum hunting rate, $\lambda_{F, \cI = 0}^{\min}$ increases with the anthropisation level. The more anthropised the regions, the fewer wild animals there are, which means that humans have to hunt more to survive.

When $\lambda_{F, \cI = 0}^{\min} < \lfw $, the system converges towards the coexistence equilibrium or towards a periodic cycle around this equilibrium. When the growth of wild fauna is not affected by human feedback (\textit{i.e.} when $\beta = 0$), we are able to characterize this difference in behavior using a second hunting rate threshold. If 
$$
\lambda_{F, \cI=\beta=0}^{\max} = \lfw^* < \lfw,
$$ (where $\lfw^*$ is defined in proposition \ref{prop:stab, cI=beta=0}) then the system converges towards the limit cycle. This means that a too large hunting rate can reduce the wild fauna below a level sufficient to feed the human population. Consequently, the human population decreases before increasing again, once the wild fauna population has grown sufficiently. In this case, numerical simulations helps to assess the period and amplitude of the cycles, in order to determine the extent and duration of the populations variations.

\item When we consider immigration, the persistence of wild fauna  is no longer guaranteed. Indeed, while a human-only equilibrium $EE^{H}$ always exist, the wildlife exists only if the coexistence equilibrium $EE^{HF_W}$ exists, that is only if $1 < \N_{I>0}$. We can rewrite this condition as:
$$
\lfw < \lambda_{F, \cI > 0}^{max} = (1-\alpha) r_F\left(\dfrac{\mu_D - f_D}{m \cI} + \beta \right).
$$
Therefore, the wild fauna exists only if the hunting rate is not too large. Moreover, $\lambda_{F, \cI > 0}^{max}$ is a decreasing function of both immigration and anthropisation level. This indicates that the wild fauna is endangered by these two phenomenons, which must be controlled if a coexistence between the two populations is desired.

We have also proved the possibility of periodic solutions around the coexistence equilibrium. Due to the complexity of the expressions, this possibility was not described with a hunting rate threshold, so instead, we have studied it numerically.

\end{itemize}

\subsection{Numerical Simulations}
To complete the theoretical work, we provide numerical simulations. The parameter values used are based on ecological and anthropological reviews concerning the South-Cameroon.
\begin{table}[ht]
\centering
\begin{tabular}{|c|c|c|c|c|}
\hline 
Parameter & Unit & Range of value & \vquatre{Baseline for simulations} & Reference \\ 
\hline 
$\mu_D$ & Year$^{-1}$  & $1/50$ & $1/50$& \cite{ins_demographie}\\
$f_D$ & Year$^{-1}$ &$[0, 0.0164]$ & $0.005$& \cite{koppert_consommation_1996}\\
$m_W$ &Year$^{-1}$  &$[2.28, 73.0]$ & $4.13$ &\cite{duda_hunting_2017, avila_interpreting_2019}\\
$m_D$ & Year$^{-1}$  &$[0.17, 0.52] \times m_W$ & $0.826$&\cite{avila_interpreting_2019}\\
$r_F$ & Year$^{-1}$ & $[0.44, 0.84]$ & $0.8$ & \cite{duda_hunting_2017, robinson_intrinsic_1986}\\
$K_F$ & Ind& $[900, 34000]$ & $7200$ & \cite{janson_ecological_1990} \\
$e$ & - & $(0, 1]$ &$0.4$ & Assumed\\
$\alpha$ &-&  $[0, 1)$ &  \YD{varying} & -  \\
$\beta$ & Ind$^{-1}$ & $[0, \beta^*)$ & \YD{varying} & - \\
$\lfw$ & Ind$^{-1}\times$ Year$^{-1}$ & & \YD{varying} & - \\
$\mathcal{I}$ &  Ind$\times$Year$^{-1}$ &  & \YD{varying} & -  \\
\hline
\end{tabular}
\caption{Value range for parameters.}
\label{table:param values}
\end{table}

\YD{In the forthcoming simulations we will consider the parameter values given in Table \ref{table:param values}(column $4$), page \pageref{table:param values} with given values for the parameters $\beta$ and $\mathcal{I}$.}

First, we present bifurcation diagrams (figures \ref{fig:bifurcation, I=0}, \ref{fig:bifurcation, I=0.1}, \ref{fig:bifurcation, I=1} and \vquatre{\ref{fig:bifurcation, I=1, betaMax}}). They are plotted in the $\lfw-\alpha$ plane, and they illustrate the long term dynamics when these parameters varies.
The first one, presented in figure \ref{fig:bifurcation, I=0} is drawn for $\cI = 0$. On it, we see that the range of values of $\lfw$ and $\alpha$ for which the wildlife persists without human is very low. This means that the human population finds easily enough prey to survive, except for high level of anthropisation.

\begin{figure}[!ht]
\centering
\includegraphics[width=1\textwidth]{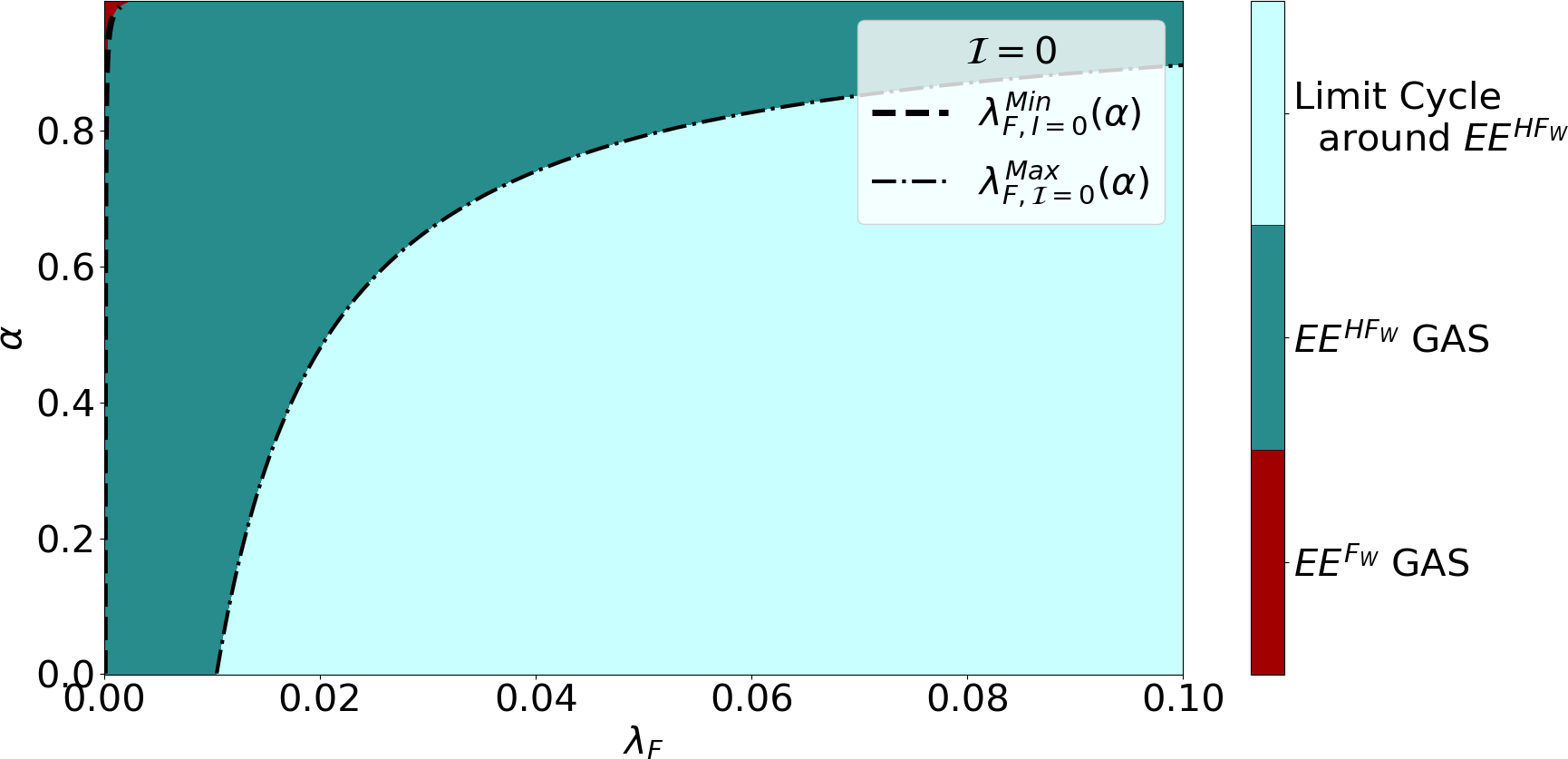}
\caption{Bifurcation diagram when $\cI = 0$ and $\beta = 0$.}
\label{fig:bifurcation, I=0}
\end{figure}

The second and third bifurcation diagrams (figures \ref{fig:bifurcation, I=0.1} and \ref{fig:bifurcation, I=1}) are respectively drawn for $\cI = 0.1$ and $\cI = 1$. As mentioned above, if the hunting rate or the anthropisation is too high, only the human population subsists. Moreover, the range of values for which coexistence is possible decreases as $\cI$ increases, as does the possibility of limit cycles. This last point is important. Indeed, the value of $F_W^*$ at equilibrium rapidly decreases when the hunting rate, the anthropisation or the immigration increase. Therefore, in order to ensure the presence of a significant wild fauna population, the system must be in a limit cycle configuration.

\vquatre{
The \YD{first three} bifurcation diagram\YD{s} are drawn for $\beta = 0$. The last one, in figure \ref{fig:bifurcation, I=1, betaMax}, is drawn for $\beta = 0.9 \beta^*$, \YD{where $\beta^*$, defined in \eqref{betastar}, is an increasing function of $\alpha$}. With the parameter values used, $\beta^*_{\alpha = 0} = 1.3 10^{-4}$ and $\beta^*_{\alpha = 0.99} = 1.3$. A comparison of this diagram with the one in figure \ref{fig:bifurcation, I=1} (same parameter values except that $\beta = 0$) shows that having $\beta > 0$ slightly increases the range of values for which the system converges towards the coexistence equilibrium. This effect is notable when $\alpha$ is close to 1, but becomes quickly imperceptible when $\alpha$ decreases.
}

All these remarks show that, according to our model and our parametrization, the situation in South Cameroon is such that wild fauna is in great danger in the long term. It also underlines that immigration of people complicates the possibility of coexistence between human and wild fauna populations, which becomes possible only at a low level of anthropisation and hunting.

\begin{figure}[!ht]
\centering
\includegraphics[width=1\textwidth]{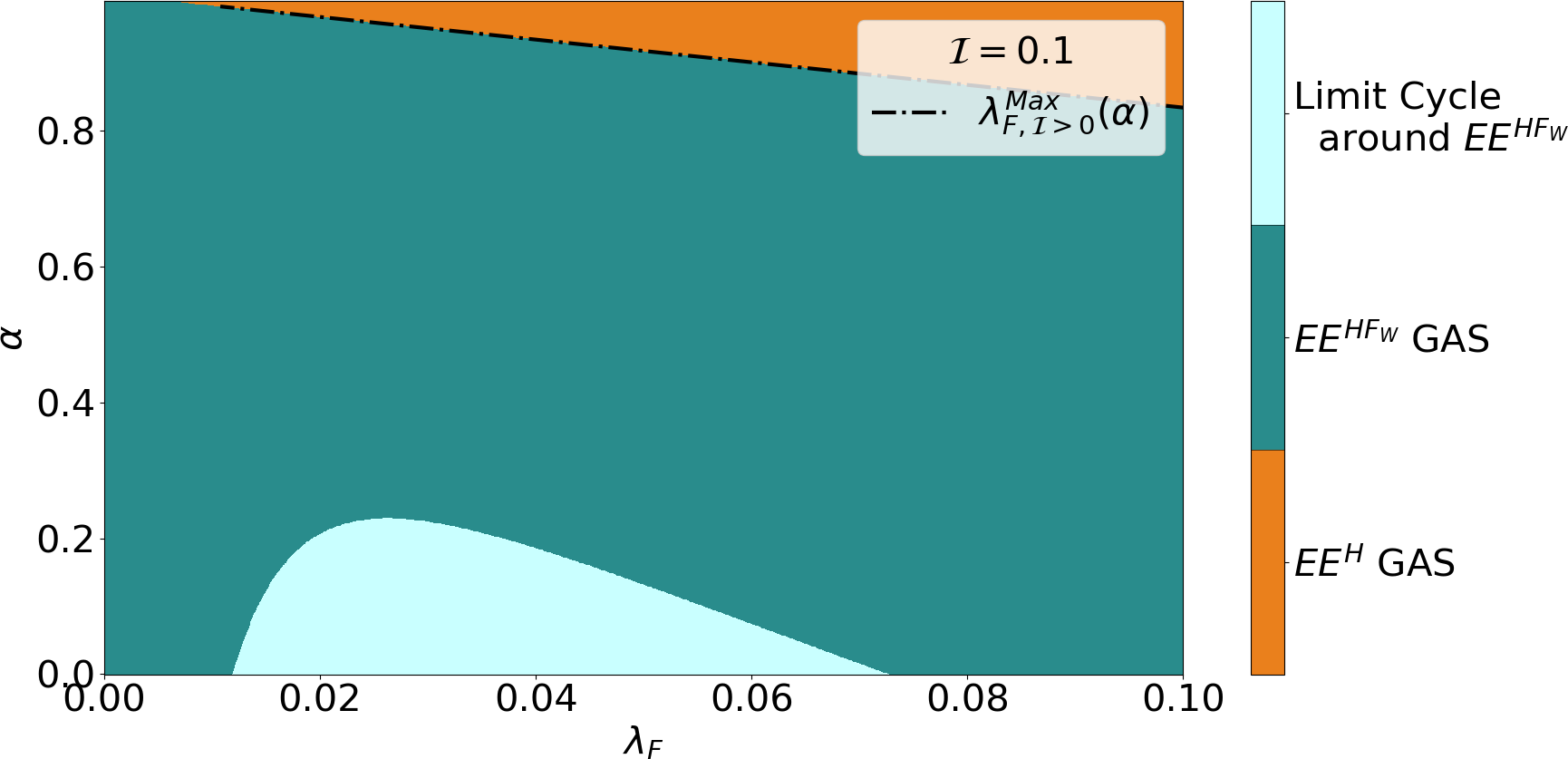}
\caption{Bifurcation diagram when $\cI = 0.1$ and $\beta = 0$.}
\label{fig:bifurcation, I=0.1}
\end{figure}
\begin{figure}[!ht]
\centering
\includegraphics[width=1\textwidth]{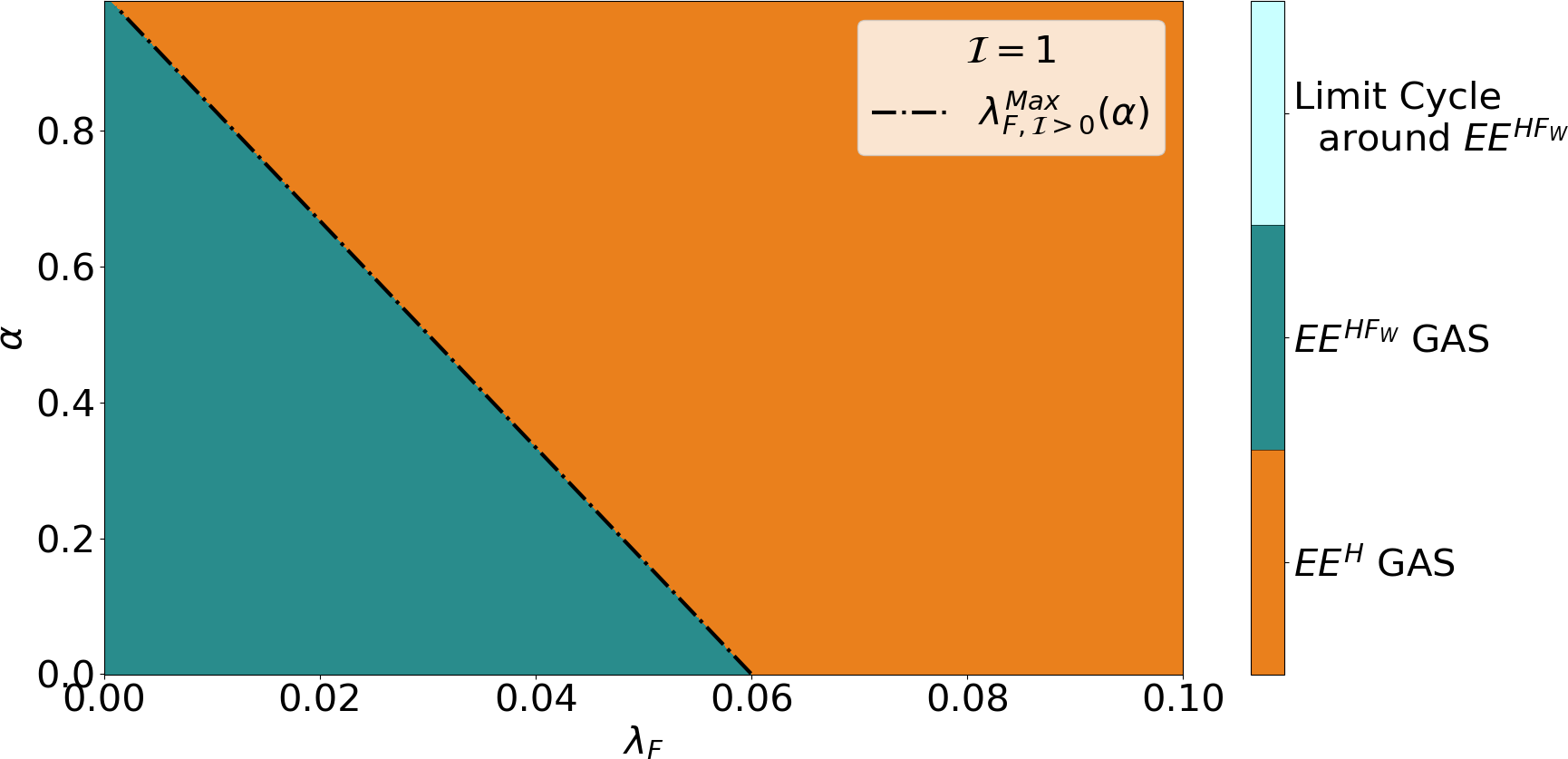}
\caption{Bifurcation diagram when $\cI = 1$ and $\beta = 0$.}
\label{fig:bifurcation, I=1}
\end{figure}
\begin{figure}[!ht]
\centering
\includegraphics[width=1\textwidth]{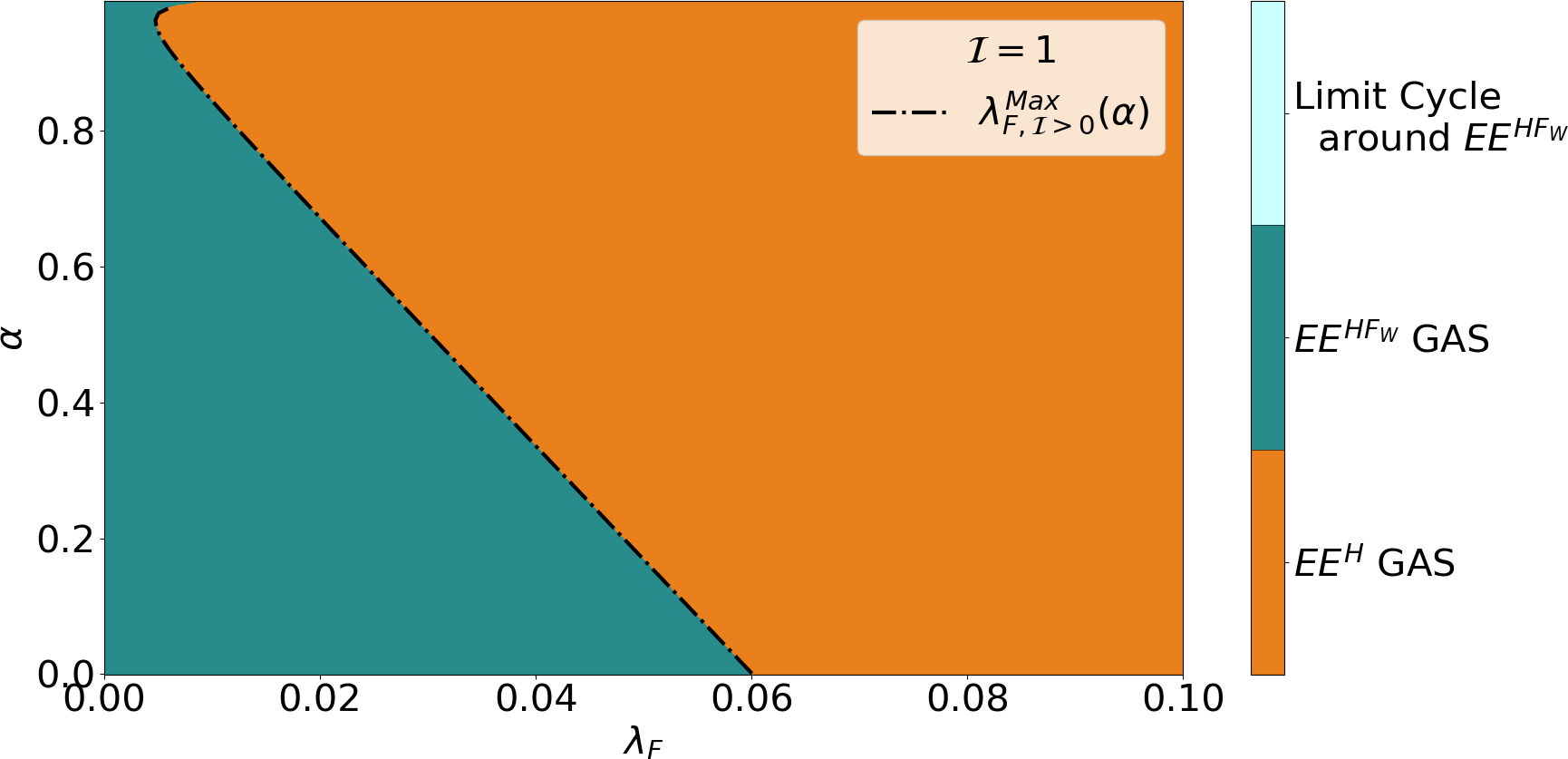}
\caption{ Bifurcation diagram when $\cI = 1$ and $\beta = 0.9\beta^*_\alpha$.} 
\label{fig:bifurcation, I=1, betaMax}
\end{figure}

\begin{figure}[!ht]
\begin{subfigure}{1\textwidth}
\centering
\includegraphics[width=1\textwidth]{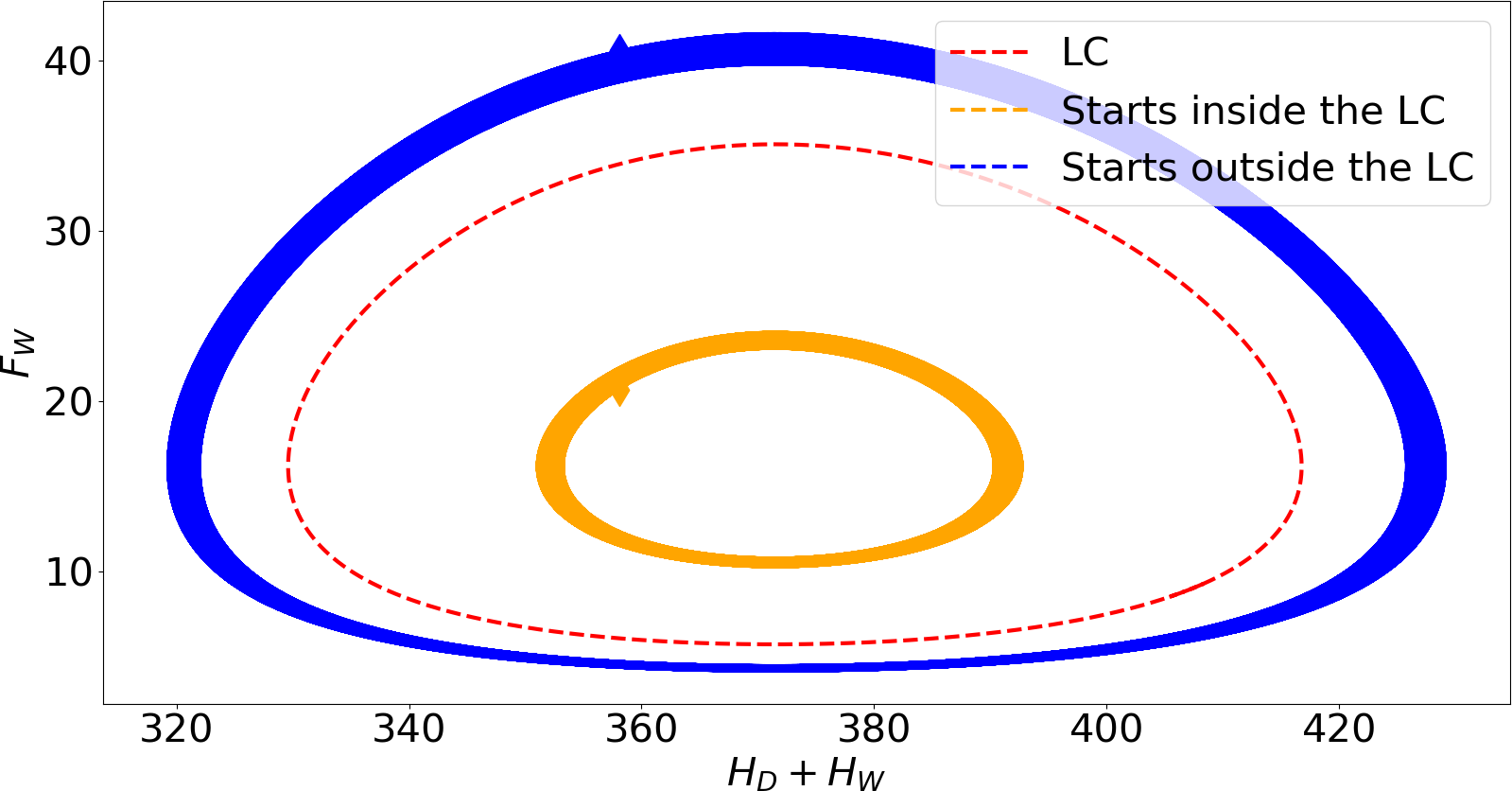}
\caption{\centering Orbits in the $\Big(H_D + H_W ; F_W\Big)$ plane converging towards the LC.}
\label{fig:LCI0, 1}
\end{subfigure}
\begin{subfigure}{1\textwidth}
\centering
\includegraphics[width=1\textwidth]{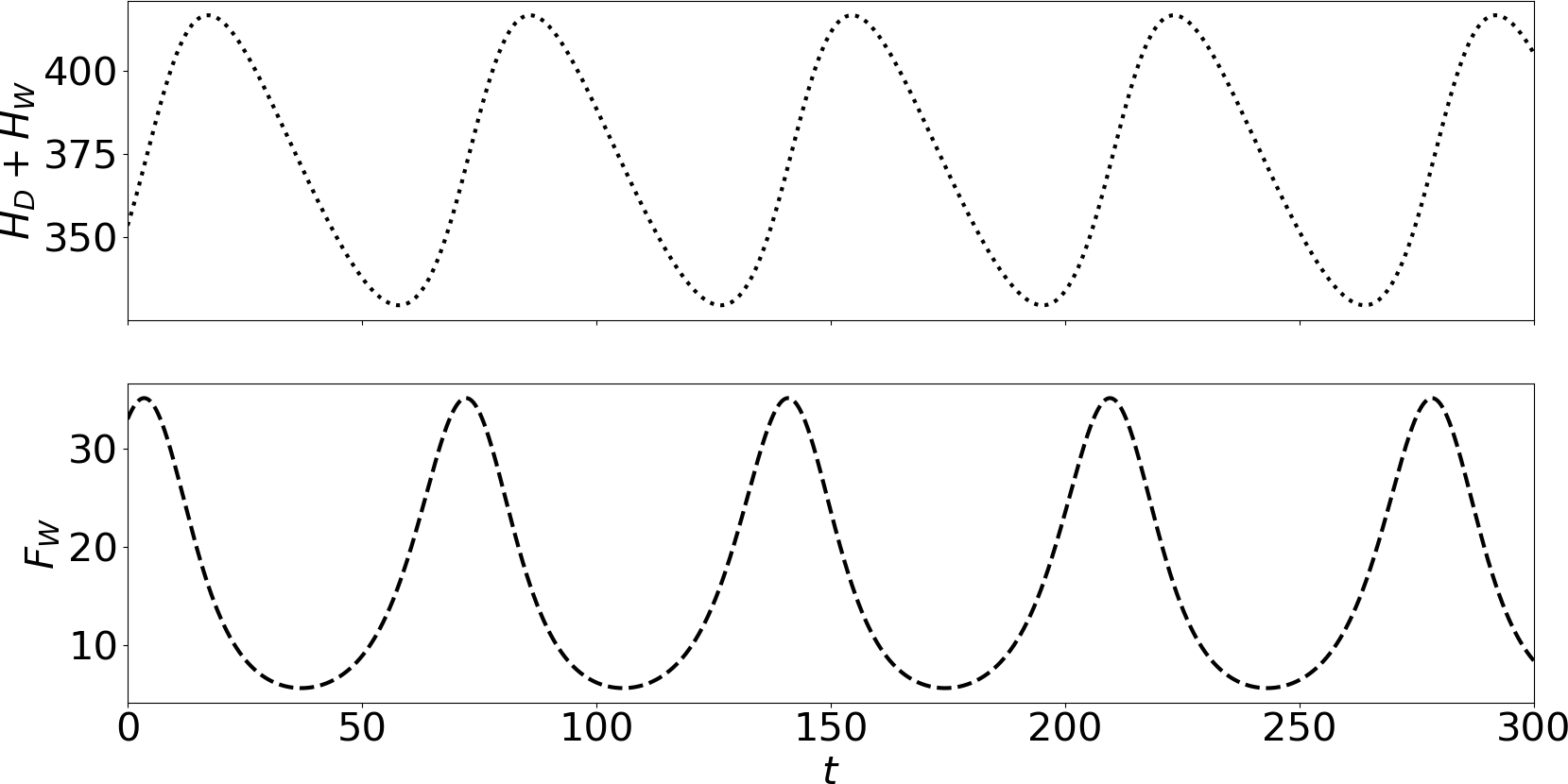}
\caption{\centering Values of $H_D + H_W$ and $F_W$ as function of time on the limit cycle.}
\label{fig:LCI0, 2}
\end{subfigure}
\caption{Illustration of the system's convergence toward a stable Limit Cycle (LC), when $\cI =0$ and $\lfw = 0.0116$. \\
Other parameter values: $\mu_D = 0.02$, $f_D = 0.005$, $m_W = 4.13$, $m_D = 0.826$, $r_F = 0.8$, $K_F = 7200$, $e=0.4$, $\alpha = 0.1$, $\beta = 0$.}
\label{fig:LCI0}
\end{figure}

\begin{figure}[!ht]
\begin{subfigure}{1\textwidth}
\centering
\includegraphics[width=1\textwidth]{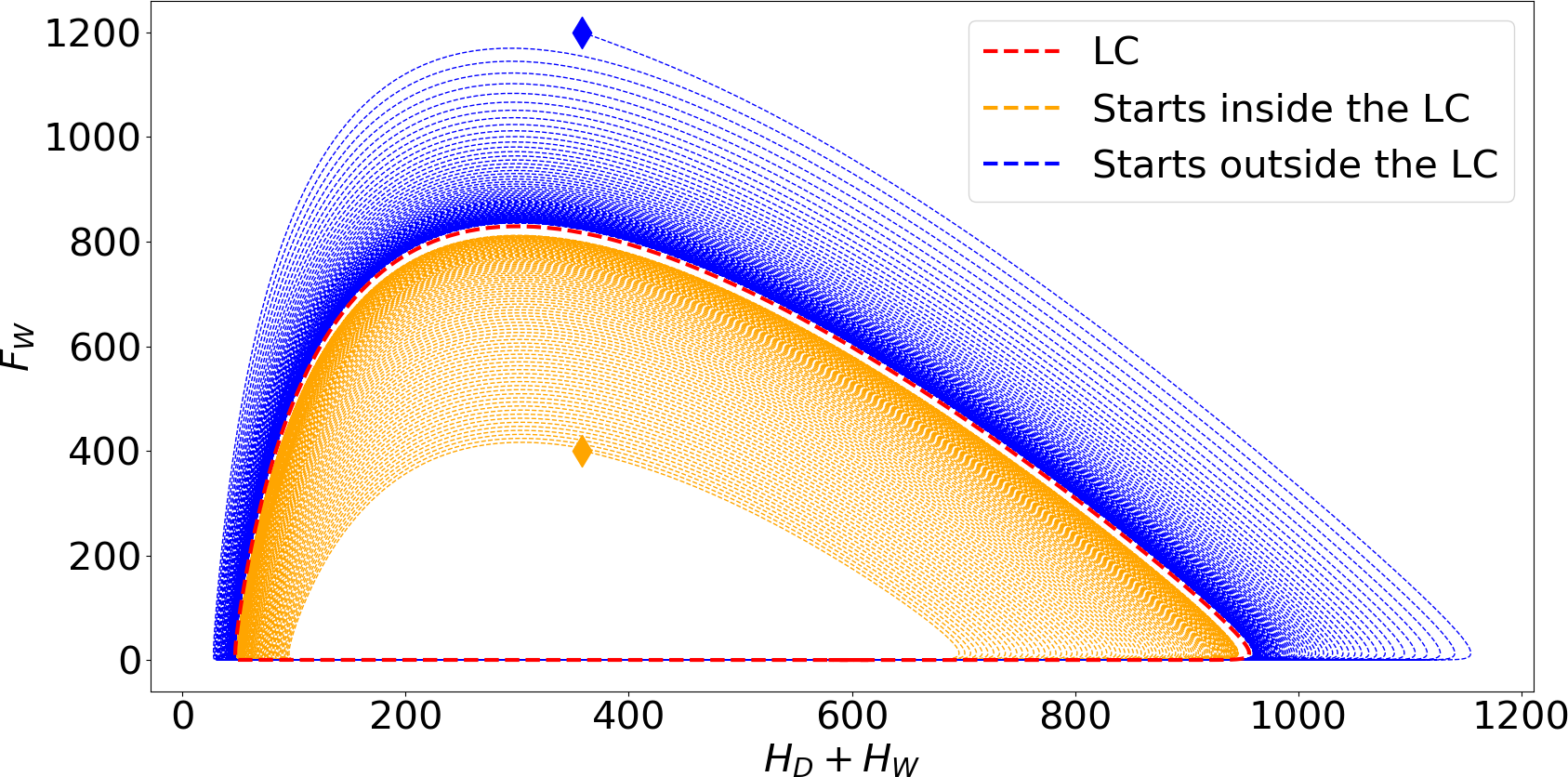}
\caption{\centering Orbits in the $\Big(H_D + H_W ; F_W\Big)$ plane converging towards the LC.}
\label{fig:LCAttoFox, 1}
\end{subfigure}
\begin{subfigure}{1\textwidth}
\centering
\includegraphics[width=1\textwidth]{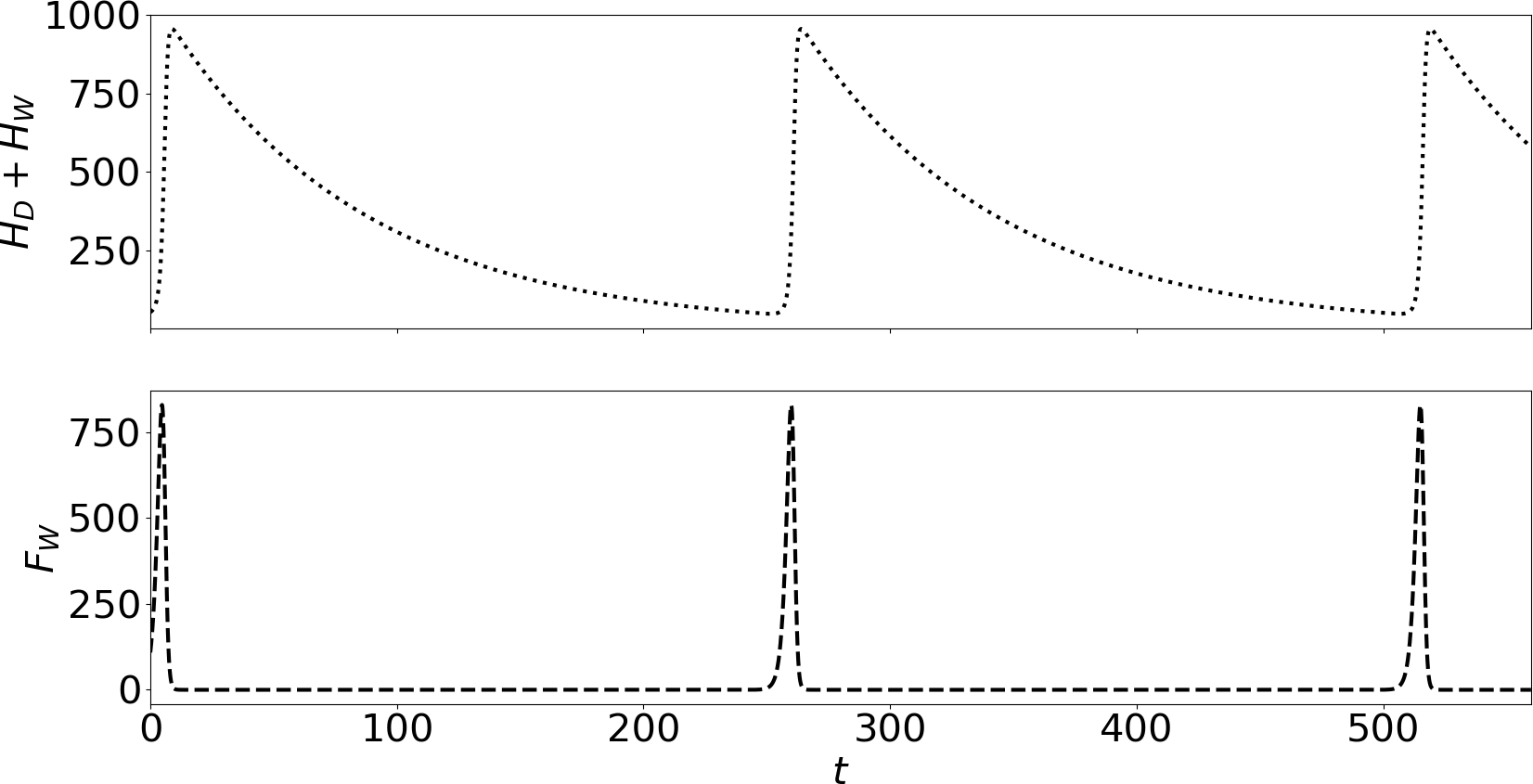}
\caption{\centering Values of $H_D + H_W$ and $F_W$ as function of time on the limit cycle.}
\label{fig:LCAttoFox, 2}
\end{subfigure}
\caption{Illustration of the system's convergence toward a stable Limit Cycle (LC), when $\cI =0$ and $\lfw = 0.01425$. \\
Other parameter values: $\mu_D = 0.02$, $f_D = 0.005$, $m_W = 4.13$, $m_D = 0.826$, $r_F =0.8$, $K_F = 7200$, $e=0.4$, $\alpha = 0.1$, $\beta = 0$.}
\label{fig:LCIAttoFox}
\end{figure}

\begin{figure}[!ht]
\begin{subfigure}{1\textwidth}
\centering
\includegraphics[width=1\textwidth]{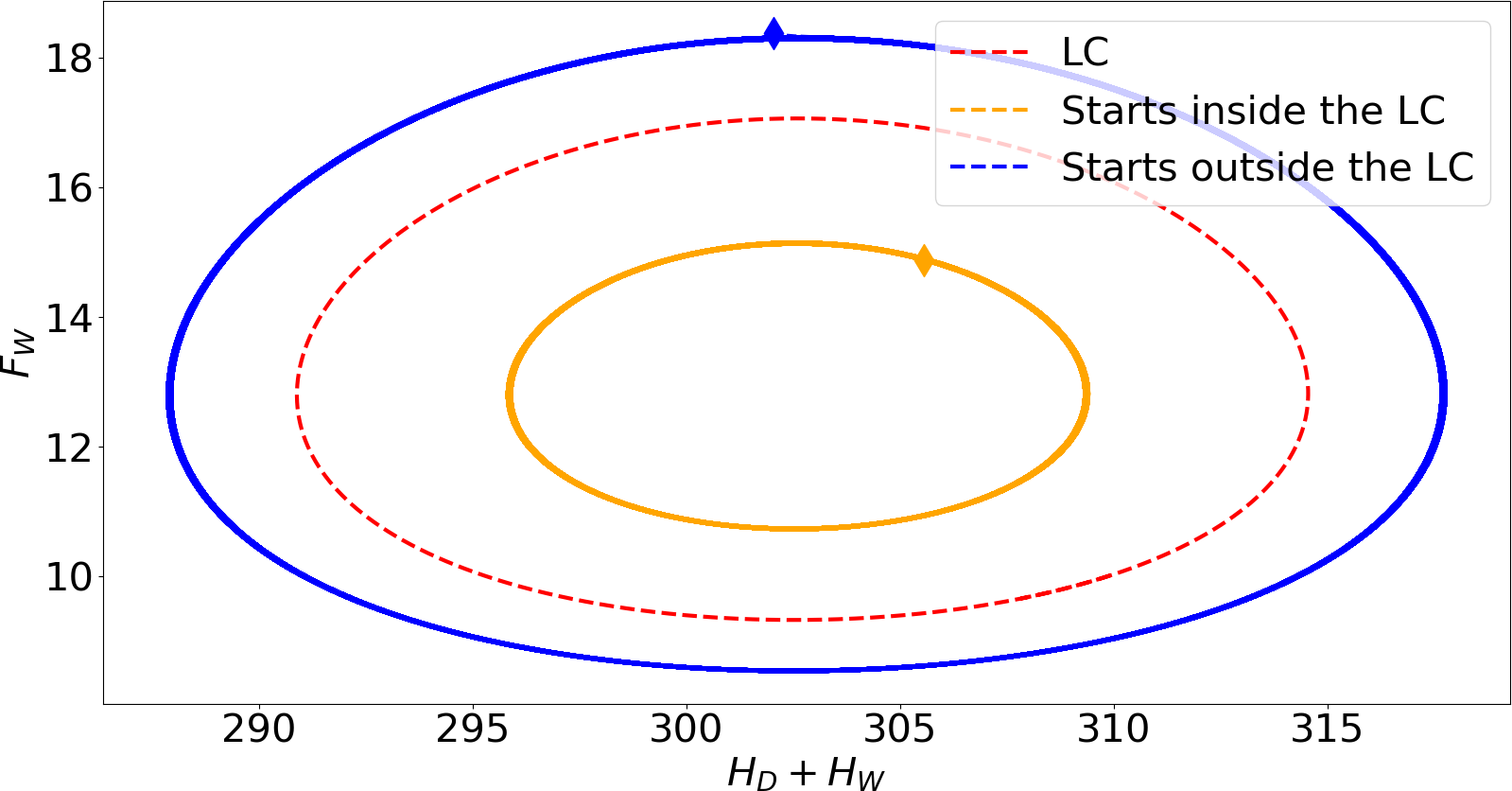}
\caption{\centering Orbits in the $\Big(H_D + H_W ; F_W\Big)$ plane converging towards the LC.}
\label{fig:LCI, 1}
\end{subfigure}
\begin{subfigure}{1\textwidth}
\centering
\includegraphics[width=1\textwidth]{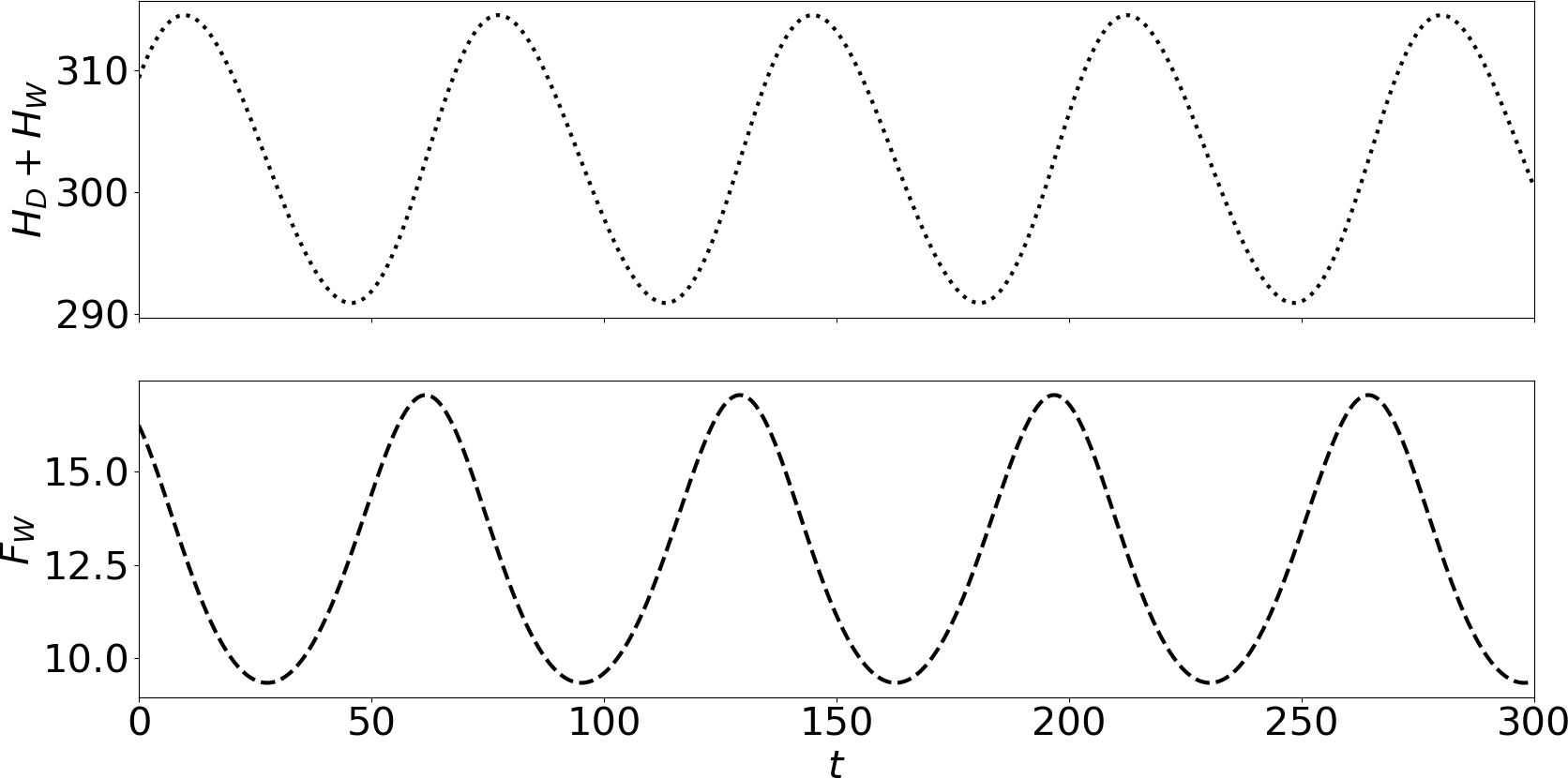}
\caption{\centering Values of $H_D + H_W$ and $F_W$ as function of time on the limit cycle.}
\label{fig:LCI, 2}
\end{subfigure}
\caption{Illustration of the system's convergence toward a stable Limit Cycle (LC), when $\cI =0.1$ and $\lfw = 0.01425$. \\
Other parameter values: $\mu_D = 0.02$, $f_D = 0.005$, $m_W = 4.13$, $m_D = 0.826$, $r_F =0.8$, $K_F = 7200$, $e=0.4$, $\alpha = 0.1$, $\beta = 0$.}
\label{fig:LCI}
\end{figure}

\begin{figure}[!ht]
\centering
\begin{subfigure}{1\textwidth}
\centering
\includegraphics[width=1\textwidth]{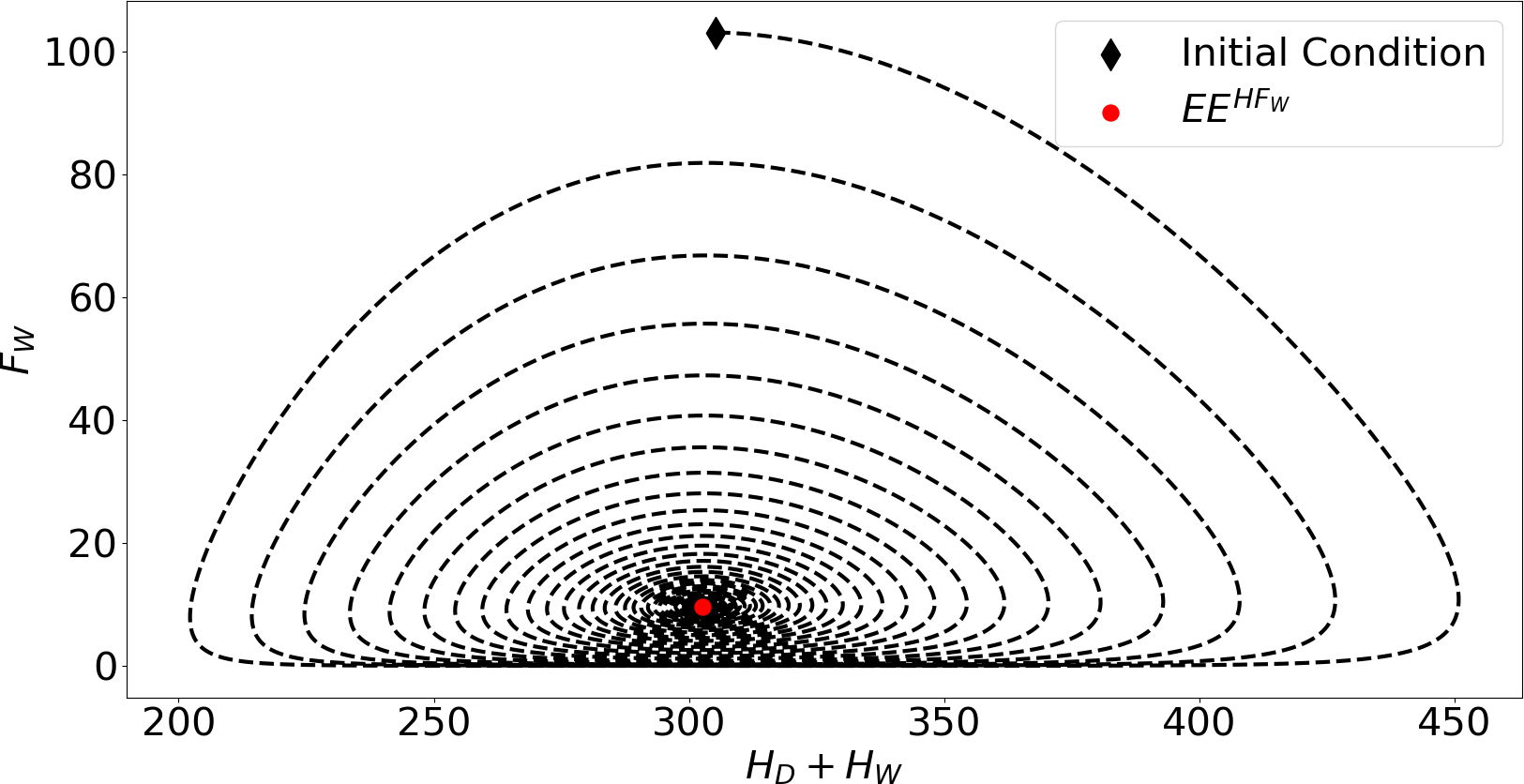}
\caption{\centering Orbits in the $\Big(H_D + H_W ; F_W\Big)$ plane converging towards the coexistence equilibrium.}
\end{subfigure}
\begin{subfigure}{1\textwidth}
\centering
\includegraphics[width=1\textwidth]{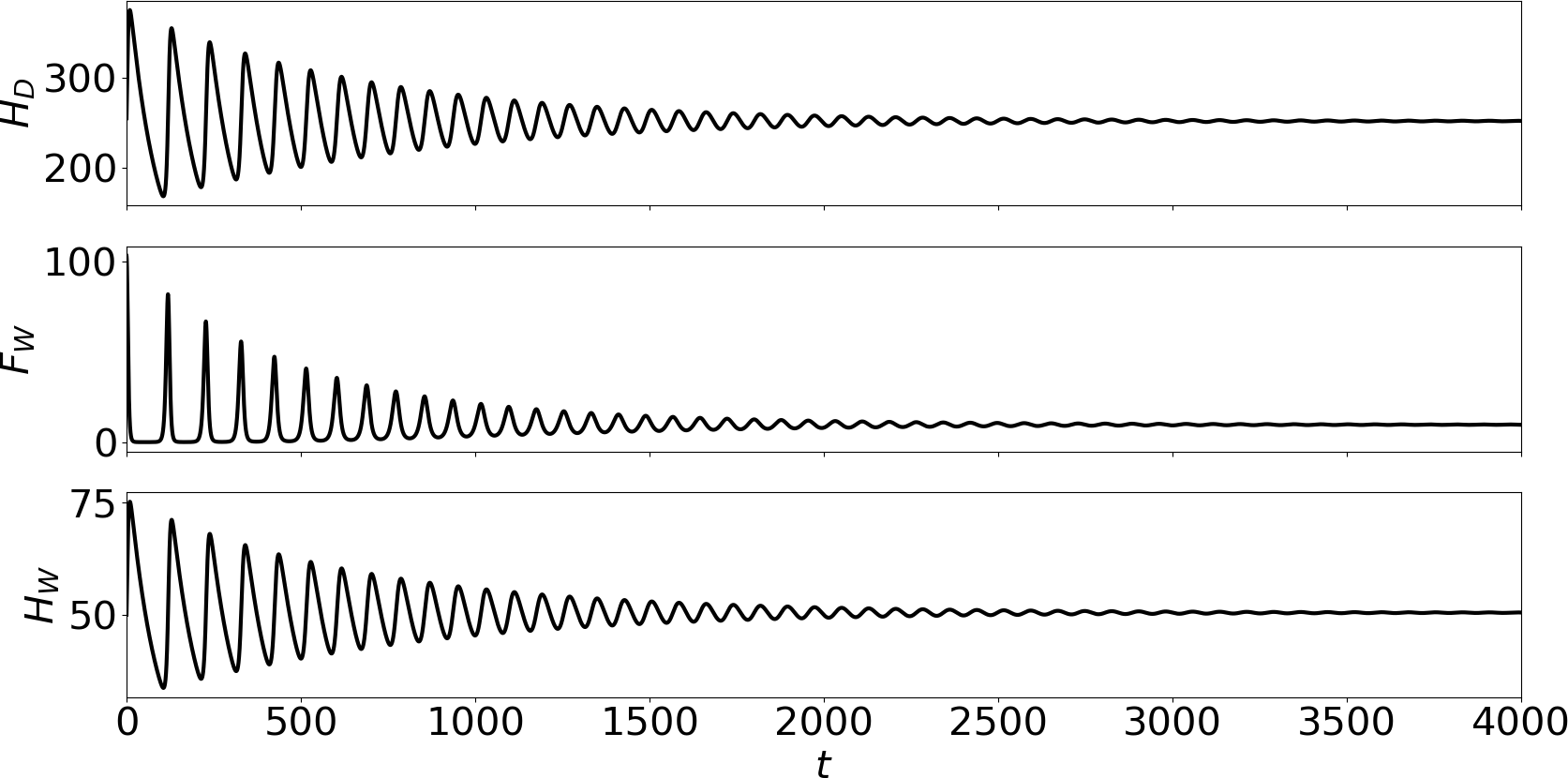}
\caption{\centering Values of $H_D$, $F_W$ and $H_W$ as function of time.}
\end{subfigure}
\caption{Illustration of the system's convergence towards the coexistence-equilibrium, when $\cI = 1$ and $\lfw = 0.01425$. \\
Other parameter values: $\mu_D = 0.02$, $f_D = 0.005$, $m_W = 4.13$, $m_D = 0.826$, $r_F =0.8$, $K_F = 7200$, $e=0.4$, $\alpha = 0.1$, $\beta = 0$.}
\label{fig:HFWI}
\end{figure}

We also present numerical solutions. They all share the same parameters values except for $\cI$ and $\lfw$. The values are respectively
\begin{itemize}
\item $\cI = 0$ and $\lfw = 0.0116$ for figure \ref{fig:LCI0},
\item $\cI = 0$ and $\lfw = 0.01425$ for figure \ref{fig:LCIAttoFox},
\item $\cI = 0.1$ and $\lfw = 0.01425$ for figure \ref{fig:LCI},
\item $\cI = 1$ and $\lfw = 0.01425$ for figure \ref{fig:HFWI}.
\end{itemize}
In the three first cases, the solutions converge towards a limit cycle (see the bifurcation diagrams). For the last configuration, the system converges towards the coexistence equilibrium.

A comparison between figures \ref{fig:LCI0} and \ref{fig:LCIAttoFox} shows that increasing the hunting rate (from $0.0116$ to $0.01425$) increases the amplitude of the limit cycle, but also the period for which the wild fauna population is low. A small increase of the immigration (from $\cI = 0$ to $\cI = 0.1$) is sufficient for the period and the oscillations to reduce. This is shown in figure \ref{fig:LCI}.
\YD{A slight} increase of the immigration parameter (from $\cI = 0.1$ to $\cI = 1$) changes the long term dynamics of the system. It is now converging, still by oscillating, towards the coexistence equilibrium $EE^{HF_W}$. This is shown on figure \ref{fig:HFWI}. \YD{In fact it is clear that the immigration of only a few individuals seems to be enough to eliminate the co-existence equilibrium, and, thus, shows that the immigration process is maybe not well modeled.}

\section{Conclusion} \label{sec:conclusion}
To better understand the risks of overhunting and habitat anthropisation to the coexistence between a forest-dwelling human population and wildlife, we developed a model based on ordinary differential equations. Initially, we studied the long-term behavior of the model using a quasi-stationary approach. However, numerical simulations showed that this approach yielded incomplete results. To address this issue, we conducted a second analysis using the theory of monotone systems. Depending on the values of the different thresholds, the system can converge towards the extinction of one population, the coexistence of both populations, or a limit cycle. We provided an ecological interpretation of the thresholds.

Our bifurcation diagrams underlined the importance of maintaining low level of hunting and habitat anthropisation in order to preserve coexistence, especially in areas subject to human immigration. The numerical simulations showed that due to the hunting activities, the wild fauna population quickly declines when the human population becomes large. 

We opted for a Holling-type 1 functional to model hunting activities in order to facilitate our mathematical analysis, which ultimately proved to be not so simple. This functional, because it is unbounded, causes the wild fauna population to decrease rapidly when the human population increases. 
In our future works, we intend to consider a Holling-type 2 functional (with $\Lambda_F(F_W) = \dfrac{\lfw F_W}{1+ \lfw \theta F_W}$, see \cite{feng_mathematical_2018}), which introduces a decelerating and bounded intake rate, and also includes our model, when $\theta = 0$. However, if the Holling-type 2 functional can be more realistic, it also introduces another non-linearity in the model, which will certainly complicate the mathematical analysis.

In the absence of recent literature, the parameter values used can be outdated. For example, the upper-bound for parameter $f_D$ is probably underestimated. A larger value of $f_D$ means that the human population will be larger, further threatening the persistence of wildlife.
\vquatre{
For some parameter values (for example when $\mu_D - f_D \simeq 0$, or for large value of $\lfw$), an atto-fox problem may arise. In \cite{lobry_migrations_2015}, the authors propose several solutions to resolve this problem. One of them is to introduce a constant immigration term for the wild animal population, suggesting that there is a reservoir of wild fauna unaffected by human activities somewhere. This can be, for example, a National Park with strict legislation \YD{of non-hunting}. However, modeling the immigration by a constant term is not satisfying, whether for human population or wildlife: the model is highly sensitive to this parameter, which is difficult to estimate. Another solution to the atto-fox problem proposed in \cite{lobry_migrations_2015} is to develop a multi-patch approach for the resource. This \YD{c}ould be an interesting extension of the model.
}

\section*{Acknowledgments}
This work takes place within the I-CARE project, funded by ExposUM Institute of the University of Montpellier, Occitanie Region, and the investment programme France 2030, ANR-21-EXES-0005. YD is (partially) supported by the DST/NRF SARChI Chair in Mathematical Models and Methods in Biosciences and Bioengineering at the University of Pretoria (Grant 82770). YD and VY-D acknowledges the support of the AFRICAM project funded by the French Agency for Development (AFD). YD acknowledges the support of the Conseil R\'egional de la R\'eunion, the Conseil d\'epartemental de la R\'eunion, the European Regional Development Fund (ERDF), and the Centre de Coop\'eration Internationale en Recherche Agronomique pour le D\'eveloppement (CIRAD).

\bibliographystyle{plain}
\bibliography{Biblio/Math, Biblio/Context, Biblio/interactionsHumanEnvironmentModel}

\newpage
\section*{Appendix}
\begin{appendix}
\section{Useful Theorems Coming from the Literature} \label{sec:litterature theorems}

On the following, we consider the system of differential equation
\begin{equation}
\dfrac{dy}{dt} = f(y), \quad y(0) = y_0,
\label{equation:generic system}
\end{equation}
defined on an open subset $\mathcal{D}$ of $\mathbb{R}^n$. We assume that $f \in C^1(\mathcal{D})$, and we note $\mathcal{J}(y)$ the Jacobian of $f$ at $y$.

\begin{theorem}  \label{theorem:Poincaré-Bendixson} [Poincaré-Bendixson's Theorem  \cite{wiggins_introduction_2003}]

The following theorem holds true when $\mathcal{D}$ is an open subset of $\mathbb R ^2$.

Let $\mathcal{M}$ be a positively invariant region for $f$, containing a finite number of fixed points. Let $y \in \mathcal{M}$ and consider its $\omega$-limit set, $\omega(y)$. Then one of the following possibilities holds:
\begin{enumerate}
\item $\omega(y)$ is a fixed point;
\item $\omega(y)$ is a closed orbit;
\item $\omega(y)$ consists of a finite number of fixed points $y_1, \ldots, y_n$ and regular heteroclinic or homoclinic orbits joining them.
\end{enumerate}
\end{theorem}

\begin{theorem} \label{theorem:Dulac} [Dulac's Criteria \cite{perko_differential_1996}]

The following theorem holds true when $\mathcal{D}$ is a simply connected region in $\mathbb R ^2$. 

If there exists a function $B \in \mathbb{C}^1(\mathcal{D})$ such that $\nabla \cdot (Bf)$ is not identically zero and does not change sign in $\mathcal{D}$, then the system has no closed orbit lying entirely in $\mathcal{D}$.
\end{theorem}

\begin{theorem} \label{theorem:Routh-Hurwitz} [Routh-Hurwitz's Criterion \cite{wiggins_introduction_2003}]

Let $y^*$ be a critical point of the system. We note $J^*$ the Jacobian matrix at this point, and $\chi_{J^*}$ its characteristic polynomial.

\begin{itemize}
\item When $\mathcal{D}\subset \mathbb{R}^2$, we have $\chi_{J^*} = X^2 - \Tr(J^*) X + \det(J^*)$ and $y^*$ is Locally Asymptotically Stable (LAS) if $\Tr(J^*) < 0$ and $\det(J^*) > 0$.
\item When $\mathcal{D}\subset \mathbb{R}^3$, we have $\chi_{J^*} = X^3 - \Tr(J^*) X^2 + a_1 X - \det(J^*)$, where $$a_1 = J^*_{11}J^*_{22} + J^*_{11} J^*_{33} + J^*_{22}J^*_{33} - J^*_{31}J^*_{13} - J^*_{32}J^*_{23} - J^*_{12}J^*_{21}.$$  $y^*$ is LAS if $\Tr(J^*) < 0$, $a_1 > 0$, $\det(J^*) < 0$ and $-\Tr(J^*) a_1 + \det(J^*) > 0$.
\end{itemize}
\end{theorem}

\begin{theorem} \label{theorem:Vidyasagar} [Vidyasagar's Theorem  \cite{vidyasagar_decomposition_1980, dumont_mathematical_2012}]

Assume that the system can written as
\begin{equation}
\def\arraystretch{2}
\left\{ \begin{array}{l}
\dfrac{dy_1}{dt} = f_1(y_1), \\
\dfrac{dy_2}{dt} = f_2(y_1, y_2) 
\end{array} \right.
\label{equation: eqVidyasagar}
\end{equation}

with $(y_1, y_2) \in \mathbb{R}^{n_1} \times\mathbb{R}^{n_2}$. Let $(y_1^*, y^*_2)$ be an equilibrium point.
If $y^*_1$ is Globally Asymptotically Stable (GAS) in $\mathbb{R}^{n_1}$ for the system $\dfrac{dy_1}{dt} = f_1(y_1)$, and if $y^*_2$ is GAS in $\mathbb{R}^{n_2}$ for the system $\dfrac{dy_2}{dt} = g(y_1^*, y_2)$, then $(y_1^*, y_2^*)$ is (locally) asymptotically stable for system \eqref{equation: eqVidyasagar}. Moreover, if all trajectories of \eqref{equation: eqVidyasagar} are forward bounded, then $(y_1^*, y_2^*)$ is GAS for \eqref{equation: eqVidyasagar}.
\end{theorem}

\begin{theorem} \label{theorem:Tikhonov} [Thikonov's Theorem \cite{banasiak_methods_2014}]

We consider the following systems of ODEs,
\begin{equation}\label{orginalProbelm}
\def\arraystretch{2}
\left\lbrace \begin{array}{l}
\dfrac{d x}{dt} = f(t,x,y,\epsilon), \quad x(0) = x_0, \\
\epsilon \dfrac{d y}{dt} = g(t,x,y,\epsilon), \quad y(0) = y_0, \\
\end{array} \right.
\end{equation}and the following assumptions:
\begin{enumerate}
\item \textbf{Assumption 1:} Assume that the functions $f, g$:
$$
f : [0, T]\times \mathcal{\bar{U}} \times \mathcal{V} \times [0, \epsilon_0] \rightarrow \mathbb{R}^n,
$$
$$
g : [0, T]\times \mathcal{\bar{U}} \times \mathcal{V} \times [0, \epsilon_0] \rightarrow \mathbb{R}^m
$$
are continuous and satisfy the Lipschitz condition with respect to the variables $x$ and $y$ in $[0, T]\times \mathcal{\bar{U}} \times \mathcal{V}$, where $\mathcal{\bar{U}}$ is a compact set in $\mathbb{R}^n$, $\mathcal{V}$ is a bounded open set in $\mathbb{R}^m$ and $T, \epsilon_0 > 0$.
\item \textbf{Assumption 2:} The corresponding degenerate system reads
\begin{equation} \label{degenerateSystem}
\def\arraystretch{1.2}
\left\lbrace \begin{array}{l}
\dfrac{d x}{dt} = f(t,x,y,0), \quad x(0) = x_0, \\
0 =  g(t,x,y,0).
\end{array} \right.
\end{equation}

Assume that there exists a solution $y = \phi(t, x) \in \mathcal{V}$ of the second equation of \eqref{degenerateSystem}, for $(t,x) \in [0, T]\times \mathcal{\bar{U}}$. The solution is such that
$$
\phi \in \mathcal{C}^0([0, T]\times \mathcal{\bar{U}} ; \mathcal{V})
$$
and is isolated in $[0, T]\times \mathcal{\bar{U}}$.
\item \textbf{Assumption 3:} Consider the following auxiliary equation:
\begin{equation}\label{auxiliaryEquation}
\dfrac{d \tilde{y}}{d \tau} =  g(t,x,\tilde{y},0),
\end{equation}
where $t$ and $x$ are treated as parameters.

Assume that the solution $\tilde{y}_0 := \phi(t, x)$ of equation  \eqref{auxiliaryEquation} is asymptocially stable, uniformly with respect to $(t,x) \in [0, T]\times \mathcal{\bar{U}}$.

\item \textbf{Assumption 4:} 
Consider the reduced equation:
\begin{equation}\label{reducedEquation}
\dfrac{d\bar{x}}{dt} = f(t,\bar{x},\phi(t,\bar{x}), 0), \quad \bar{x}(0) = x_0.
\end{equation}
Assume that the function $(t,x) \mapsto f(t,x,\phi(t,x), 0)$ satisfies the Lipschitz condition with respect to $x$ in $[0, T]\times \mathcal{\bar{U}}$. Assume moreover that there exists a unique solution $\bar{x}$ of equation \eqref{reducedEquation} such that $$\bar{x}(t) \in Int \: \mathcal{\bar{U}}, \quad \forall t \in (0,T).$$

\item \textbf{Assumption 5:} We consider the equation \eqref{auxiliaryEquation} in the particular case $t=0$ and $x = x_0$:
\begin{equation}\label{auxiliaryEquation, 0}
\dfrac{d \tilde{y}}{d \tau} =  g(0, x_0,\tilde{y},0), \quad \tilde{y}(0) = y_0.
\end{equation}
Assume that $y_0$ belongs to the region of attraction of the solution $y = \phi(0, x_0)$ of equation $g(0, x_0,\tilde{y},0) = 0$.
\end{enumerate}

Let assumptions 1, 2, 3, 4, 5 be satisfied. There exists $\epsilon_0$ such that for any $\epsilon \in (0, \epsilon_0]$ there exists a unique solution $(x_\epsilon(t), y_\epsilon(t))$ of problem \eqref{orginalProbelm} on $[0,T]$ and
\begin{equation}
\def\arraystretch{1.2}
\left\lbrace \begin{array}{l}
\lim\limits_{\epsilon \rightarrow 0}{x_\epsilon(t)} = \bar{x}(t) \quad t \in [0,T], \\
\lim\limits_{\epsilon \rightarrow 0} y_\epsilon(t) = \phi(t, \bar{x}(t)) := \bar{y}(t) \quad t \in (0,T] \\
\end{array} \right.
\end{equation}
where $\bar{x}(t)$ is the solution of problem \eqref{reducedEquation}.
\end{theorem}

We continue by defining the notions used in theorem \ref{theorem:Zhu}
\begin{definition}
An open set $\mathcal{D} \in \mathbb{R}^n$ is said to be p-convex provided that for every $x, y \in \mathcal{D}$, with $x \leq y$, the line segment joining $x$ and $y$ belongs to $\mathcal{D}$.
\end{definition}

\begin{definition}\cite{kaszkurewicz_matrix_2012}
A square matrix $A \in M_n (\mathbb{R})$ is said reducible if 
for each nonempty proper subset $I$ of $N = \{1, ..., n\}$, there exists $i \in I$ and $j \in N\backslash I$ such that $A_{i,j} \neq 0$
\end{definition}

\begin{definition} \label{def:monotone}\cite{smith_monotone_1995} The system is said to be 
\begin{itemize}
\item irreducible if the Jacobian matrix $\mathcal{J}(y)$ is irreducible.
\item competitive if $\mathcal{J}(y)$ has non positive off-diagonal elements:
$$ \dfrac{\partial f_i}{\partial y_j}(y) \leq 0, \quad i \neq j.
$$
\end{itemize}
\end{definition}

\begin{theorem} \label{theorem: monotone GAS} [Theorem 6 of \cite{anguelov_mathematical_2012}]

Let the system \eqref{equation:generic system} be monotone (and in particular competitive) and let $a,b \in \mathcal{D}$ such that $a <b$, $[a, b] \subset \mathcal{D}$ and $f(b) \leq \mathbf{0} \leq f(a)$. Then the system defines a positive dynamical system on $[a, b] $. Moreover, if $[a, b] $ contains a unique equilibrium then it is GAS on $[a, b] $.
\end{theorem}

\begin{theorem}\label{theorem:Zhu} [Zhu and Smith's Theorem \cite{zhu_stable_1994}]
The following theorem holds true when $\mathcal{D} \subset \R^3$ and verifies the following assumption:
\begin{itemize}
\item $\mathcal{D}$ is an open, $p$-convex subset of $\mathbb{R}^3$,
\item $\mathcal{D}$ contains a unique equilibrium point $y^*$ and $\det(\mathcal{J}(y^*)) < 0$,
\item $f$ is analytic in $\mathcal{D}$,
\item the system is competitive and irreducible in $\mathcal{D}$,
\item the system is dissipative: For each $y_0 \in \mathcal{D}$, the positive semi-orbit through $y_0$, $\phi^+(y_0)$ has a compact closure in $\mathcal{D}$. Moreover, there exists a compact subset $\mathcal{B}$ of $\mathcal{D}$ with the property that for each $y_0 \in \mathcal{D}$, there exists $T(y_0) > 0$ such that $y(t, y_0) \in \mathcal{B}$ for $t \geq T(y_0)$.
\end{itemize}

then either $y^*$ is stable, or there exists at least one non-trivial orbitally asymptotically stable  periodic orbit in $\mathcal{D}$.
\end{theorem}

\section{Preliminary Results}
\begin{prop} \label{propBeta}
Assuming $\beta < \beta^*$, the following inequality holds
$$
\dfrac{\beta (1-\alpha) r_F}{\lfw} \Big(1 - \dfrac{\mu_D - f_D}{\lfw m e (1-\alpha) K_F}\Big) < 1
$$
\end{prop}

\begin{proof}
\vquatre{
When $1 - \dfrac{\mu_D - f_D}{\lfw m e (1-\alpha) K_F} \leq 0$, the result is obvious. Assuming that $0 < 1 - \dfrac{\mu_D - f_D}{\lfw m e (1-\alpha) K_F}$, w}e have
\begin{equation*}
\dfrac{\beta (1-\alpha) r_F}{\lfw} \Big(1 - \dfrac{\mu_D - f_D}{\lfw m e (1-\alpha) K_F}\Big) < \dfrac{4 \mu_D - f_D}{\lfw m e (1-\alpha) K_F} \Big(1 - \dfrac{\mu_D - f_D}{\lfw m e (1-\alpha) K_F}\Big).\\
\end{equation*}

It is straightforward to show that $x(1 - x) \leq \dfrac{1}{4}$ for $x \in \mathbb{R}$. Therefore,
\begin{equation*}
\dfrac{\beta (1-\alpha) r_F}{\lfw} \Big(1 - \dfrac{\mu_D - f_D}{\lfw m e (1-\alpha) K_F}\Big) <  1.
\end{equation*}
\end{proof}
\subsection{Study of $P_F$} \label{section:study of PF}
We define the following polynomial, used to prove propositions \ref{prop:eq, cI>0}, \ref{prop:stab, cI=0} and \ref{prop:stab, cI>0}
\begin{multline}
P_F(X) := X^2 \left(\dfrac{er_F}{K_F} \right) - X \left(e(1-\alpha)r_F + \dfrac{(\mu_D - f_D) r_F}{\lfw m K_F} + \dfrac{\cI \beta r_F}{\lfw K_F} \right) + \\ \left(\dfrac{(\mu_D - f_D)(1-\alpha) r_F}{\lfw m} - \cI\Big(1 - \dfrac{(1-\alpha)\beta r_F}{\lfw} \Big) \right).
\label{polynome-Feq}
\end{multline} 

\begin{prop}\label{prop:study of PF}
When $\beta < \beta^*$ and $\cI > 0$, the following results are true:
\begin{itemize}
\item $P_F\Big((1-\alpha)K_F\Big) = -\cI < 0$.
\item $P_F\Big(\dfrac{\mu_D - f_D}{\lfw m e}\Big) < 0$.
\item $P_F\Big((1-\alpha)K_F - \dfrac{K_F \lfw}{\beta r_F}\Big) > 0$.
\item $P_F$ admits two real roots, $F_1^* \leq F_2^*$. 
\item If $\dfrac{(\mu_D - f_D) r_F}{\lfw m } \leq \cI\Big(1 - \dfrac{(1-\alpha)\beta r_F}{\lfw} \Big)$,  $F_2^*$ is positive and $F_1^*$ is non positive. If $\dfrac{(\mu_D - f_D) r_F}{\lfw m } > \cI\Big(1 - \dfrac{(1-\alpha)\beta r_F}{\lfw} \Big)$, $F^*_1$ and $F^*_2$ are positive. Moreover, $P_F$ is positive on $(-\infty, F_1^*)$, negative on $(F^*_1, F^*_2)$ and positive on $(F^*_2, +\infty)$.
\item From the precedent points, it follows that $$(1-\alpha)K_F - \dfrac{K_F \lfw}{\beta r_F} \leq F^*_1 \leq (1-\alpha)K_F, \dfrac{\mu_D - f_D}{\lfw m e} \leq F_2^* $$
\end{itemize}

\end{prop}

\begin{proof}
We have:
\begin{align*}
P_F((1-\alpha) K_F) &= \Big((1-\alpha) K_F \Big)^2 \left(\dfrac{er_F}{K_F} \right) - (1-\alpha) K_F \left(e(1-\alpha)r_F + \dfrac{(\mu_D - f_D) r_F}{\lfw m K_F} + \dfrac{\cI \beta r_F}{\lfw K_F} \right) + \\ &\left(\dfrac{(\mu_D - f_D)(1-\alpha) r_F}{\lfw m} - \cI\Big(1 - \dfrac{(1-\alpha)\beta r_F}{\lfw} \Big) \right), \\
&=(1-\alpha)^2 e K_F r_F - e(1-\alpha)^2 K_F r_F - \dfrac{(\mu_D - f_D) (1-\alpha) r_F}{\lfw m} - \dfrac{\cI \beta (1-\alpha)r_F}{\lfw}  + \\ &\dfrac{(\mu_D - f_D)(1-\alpha) r_F}{\lfw m} - \cI +\cI \dfrac{(1-\alpha)\beta r_F}{\lfw}, \\
&= -\cI< 0.
\end{align*}
Moreover,

\begin{align*}
P_F\Big(\dfrac{\mu_D - f_D}{\lfw m e}\Big) &= \left(\dfrac{\mu_D - f_D}{\lfw m e}\right)^2 \left(\dfrac{er_F}{K_F} \right) - \dfrac{\mu_D - f_D}{\lfw m e} \left(e(1-\alpha)r_F + \dfrac{(\mu_D - f_D) r_F}{\lfw m K_F} + \dfrac{\cI \beta r_F}{\lfw K_F} \right) + \\ & \left(\dfrac{(\mu_D - f_D)(1-\alpha) r_F}{\lfw m} - \cI\Big(1 - \dfrac{(1-\alpha)\beta r_F}{\lfw} \Big) \right), \\
&= - \dfrac{(\mu_D - f_D) \cI \beta r_F}{e \lfw ^2 K_F} - \cI + \dfrac{\cI (1-\alpha)r_F \beta}{\lfw}, \\
&= -\cI \left( 1 - \dfrac{(1-\alpha)r_F \beta }{\lfw} + \dfrac{(\mu_D - f_D)  \beta r_F}{e \lfw ^2 K_F} \right), \\
&= -\cI \left( 1 - \dfrac{\beta(1-\alpha)r_F  }{\lfw}\Big(1 - \dfrac{(\mu_D - f_D) }{ m e \lfw (1-\alpha) K_F}\Big) \right), \\
& < 0,
\end{align*}
thanks to proposition \ref{propBeta}. We have:

\begin{align*}
P_F\Big((1-\alpha)K_F - \dfrac{K_F \lfw}{\beta r_F}\Big) &= P_F\Big((1-\alpha)K_F\Big) + \Big(\dfrac{K_F \lfw}{\beta r_F}\Big)^2 \dfrac{er_F}{K_F} - 2(1-\alpha)K_F \dfrac{K_F \lfw}{\beta r_F}\dfrac{er_F}{K_F} + \\ &\left(e(1-\alpha)r_F + \dfrac{(\mu_D - f_D) r_F}{\lfw m K_F} + \dfrac{\cI \beta r_F}{\lfw K_F} \right) \dfrac{K_F \lfw}{\beta r_F}, \\
&= -\cI + \dfrac{K_F \lfw^2}{\beta^2 r_F} - 2 \dfrac{(1-\alpha)K_F \lfw e}{\beta} +\dfrac{(1-\alpha)K_F \lfw e}{\beta} + \dfrac{\mu_D - f_D}{\beta m} + \cI, \\
&= \dfrac{K_F \lfw^2}{\beta^2 r_F} -  \dfrac{(1-\alpha)K_F \lfw e}{\beta} + \dfrac{\mu_D - f_D}{\beta m}, \\
&= \dfrac{K_F \lfw^2}{\beta^2 r_F} \left(1 - \dfrac{\beta (1-\alpha) r_F}{\lfw} \Big(1 - \dfrac{\mu_D - f_D}{m e \lfw K_F(1-\alpha)} \Big) \right), \\
&> 0,
\end{align*}
using proposition \ref{propBeta}.

To show the last point of the proposition, we start by computing the discriminant of $P_F$, $\Delta_F$. We have:
\begin{align*}
\Delta_F &= \left(e(1-\alpha)r_F + \dfrac{(\mu_D - f_D) r_F}{\lfw m K_F} + \dfrac{\cI \beta r_F}{\lfw K_F} \right)^2 - 4\dfrac{er_F}{K_F}  \left(\dfrac{(\mu_D - f_D)(1-\alpha) r_F}{\lfw m} - \cI\Big(1 - \dfrac{(1-\alpha)\beta r_F}{\lfw} \Big) \right), \\
\Delta_F &= \left(e(1-\alpha)r_F - \dfrac{(\mu_D - f_D) r_F}{\lfw m K_F}\right)^2 + \dfrac{\cI \beta r_F}{\lfw K_F} \left(\dfrac{\cI \beta r_F}{\lfw K_F} + 2\dfrac{(\mu_D - f_D) r_F}{\lfw m K_F} + 2e(1-\alpha)r_F \right) + 4\dfrac{er_F}{K_F}  \cI\Big(1 - \dfrac{(1-\alpha)\beta r_F}{\lfw} \Big), \\
\Delta_F & > 0.
\end{align*}

Therefore, $P_F$ admits two real roots. Their sign depends on the sign of the constant coefficient. $P_F$ admits:
\begin{itemize}
\item One non positive root $F^*_1$ and one positive root $F^*_2$ if $$\dfrac{(\mu_D - f_D)(1-\alpha) r_F}{\lfw m} - \cI\Big(1 - \dfrac{(1-\alpha)\beta r_F}{\lfw} \Big) \leq 0 \Leftrightarrow \dfrac{(\mu_D - f_D) r_F}{\lfw m } \leq \cI\Big(1 - \dfrac{(1-\alpha)\beta r_F}{\lfw} \Big).$$
\item Two positive roots $F^*_1\leq  F^*_2$ if $\dfrac{(\mu_D - f_D) r_F}{\lfw m } > \cI\Big(1 - \dfrac{(1-\alpha)\beta r_F}{\lfw} \Big)$.
\end{itemize}
They are given by:

\begin{equation*}
F_i^* = \dfrac{K_F(1-\alpha)}{2}\left(1 \pm \dfrac{\sqrt{\Delta_F}}{e(1-\alpha)r_F}\right) + \dfrac{\mu_D - f_D}{2\lfw m e} + \dfrac{\cI \beta}{2\lfw e}, \quad i=1,2.
\end{equation*}
\end{proof}

\subsection{Proof of Proposition \ref{prop:stab, cI=beta=0}} \label{sec:stab, cI = beta = 0}

\begin{proof}
When $\beta = 0$, we have
\begin{multline} \label{DeltaStab, generalCase}
\Delta_{Stab, \cI =\beta = 0} =  \left(\mu_D - f_D + m_D + m_W + r_F\dfrac{F_W^*}{K_F} \right) \\ \times   \left( \mu_D -f_D + m_D + m_W \right) - m_D e \lfw \Big(K_F(1-\alpha) - F_W^* \Big),
\end{multline}

where $F_W^* = \dfrac{\mu_D - f_D}{\lfw m e}$. Therefore,

\begin{multline*}
\Delta_{Stab, \cI=\beta = 0} > 0 \\
\Leftrightarrow \left(\mu_D - f_D + m_D + m_W + r_F \dfrac{\mu_D - f_D}{\lfw K_F m e} \right) \times   \left( \mu_D -f_D + m_D + m_W \right) > \\ m_D e \lfw \left(K_F(1-\alpha) - \dfrac{\mu_D - f_D}{\lfw m e} \right), \\
\Leftrightarrow (\mu_D - f_D + m_D + m_W)^2 + r_F \dfrac{\mu_D - f_D}{\lfw K_F m e}  \times   \left( \mu_D -f_D + m_D + m_W \right) > \\ m_D e \lfw K_F(1-\alpha) - (\mu_D - f_D)m_W , \\
\Leftrightarrow \lfw (\mu_D - f_D + m_D + m_W)^2 + r_F \dfrac{\mu_D - f_D}{K_F m e}  \times   \left( \mu_D -f_D + m_D + m_W \right) > \\ m_D e \lfw^2 K_F(1-\alpha) - \lfw (\mu_D - f_D)m_W , \\
\Leftrightarrow 0 > \lfw^2 (1-\alpha) K_F  m_D e - \lfw \Big((\mu_D - f_D + m_D + m_W)^2 +(\mu_D - f_D)m_W \Big) - \\ \dfrac{r_F (\mu_D - f_D) }{K_F m e}  \big( \mu_D -f_D + m_D + m_W \big).\\
\end{multline*}

We define 
\begin{multline*}
P_{\Delta_{Stab, \cI= \beta = 0}}(X) := X^2 (1-\alpha) K_F  m_D e - X \Big((\mu_D - f_D + m_D + m_W)^2 +(\mu_D - f_D)m_W \Big) - \\ \dfrac{r_F (\mu_D - f_D) m_W}{K_F m_D e}  \big( \mu_D -f_D + m_D + m_W \big),
\end{multline*} 

such that we have 
\begin{equation}
\Delta_{Stab, \cI= \beta = 0} > 0 \Leftrightarrow P_{\Delta_{Stab, \cI= \beta = 0}}(\lfw) < 0.
\label{equation: equivalence DeltaStab}
\end{equation}

$P_{\Delta_{Stab, \cI = \beta = 0}}$ has a positive dominant coefficient, and its other coefficients are negative. So,  $P_{\Delta_{Stab, \cI= \beta = 0}}$ admits a unique positive root, noted $\lfw^*$, given by:
\begin{multline}
\lfw^* = \\
 \dfrac{\left[m_{W}(\mu_{D}-f_{D})+\big(\mu_{D}-f_{D}+m_{D}+m_{W})^{2}\right]\left(1+\sqrt{1+4\dfrac{(1-\alpha)m_{W}r_{F}\left(\mu_{D}-f_{D}\right)\big(\mu_{D}-f_{D}+m_{D}+m_{W})}{\left[m_{W}(\mu_{D}-f_{D})+\big(\mu_{D}-f_{D}+m_{D}+m_{W})^{2}\right]^{2}}}\right)}{2em_D (1-\alpha) K_F }.
\end{multline}

Moreover, $P_{\Delta_{Stab, \cI = \beta = 0}}$ is negative on $\left[0, \lfw^* \right)$ and positive on $\left(\lfw ^*, +\infty \right)$. Using \eqref{equation: equivalence DeltaStab}, we obtain that $EE^{HF_W}_{\cI = \beta = 0}$ is asymptotically stable if $\lfw  < \lfw ^*$.
\end{proof}

\subsection{Proof of Global Stability of $EE^{h}$ For System \eqref{equation:model H GAS}} \label{sec:lemma H GAS}
In this section we proof that $EE^{h}$ is GAS for system \eqref{equation:model H GAS} using the Vidyasagar's theorem (see theorem \ref{theorem:Vidyasagar}, page \pageref{theorem:Vidyasagar}). We note $y_1 = (h_D, h_W)$, $y_2 = f_W$ and $g_1$, $g_2$ the corresponding right hand side of system \eqref{equation:model H GAS}.

\begin{itemize}
\item We have that $y_1^* = \Big(\dfrac{\cI}{\mu_D - f_D}, -m\dfrac{\cI}{\mu_D - f_D} \Big)$ is LAS for system $\dfrac{dy_1}{dt} = g_1(y_1)$. Indeed, the characteristic polynomial of the Jacobian of $g_1$ at $y_1^*$ is $X^2 - X(f_D - \mu_D - m_D - m_W) + (\mu_D - f_D) m_W$, and its roots have a negative real parts. 

The function $B(h_D, h_W) = \dfrac{1}{h_Wh_D}$ is a Dulac's function for system $\dfrac{dy_1}{dt} = g_1(y_1)$. Hence, this system does not admit a closed orbit such that periodic solutions, homoclinic or heteroclinic loop.
Therefore, according to the Poincaré-Bendixson's theorem, $y_1^*$ is GAS for $\dfrac{dy_1}{dt} = g_1(y_1)$.
\item Since $\N_{\cI > 0} \leq 1$, we have $g_2(y_1^*, y_2) \geq 0$ and therefore $y_2^* = 0$ is GAS for $\dfrac{dy_2}{dt} = g_2(y_1^*, y_2)$.
\item The region $\Omega_{eq}$ is an invariant set for system \eqref{equation:model H GAS}, and it is compact. Therefore, all the trajectories with initial condition in $\Omega_{eq}$ are bounded.
\end{itemize}
All the assumptions of the Vidyasagar's theorem hold true, and therefore $EE^{h}$ is GAS for system \eqref{equation:model H GAS}.
\end{appendix}

\end{document}